\input amstex
\documentstyle {amsppt}

\pagewidth{32pc} 
\pageheight{45pc} 
\mag=1200
\baselineskip=15 pt

\hfuzz=5pt
\topmatter
\NoRunningHeads 
\title Finite Localities II 
\endtitle
\author Andrew Chermak
\endauthor
\affil Kansas State University
\endaffil
\address Manhattan Kansas
\endaddress
\email chermak\@math.ksu.edu
\endemail
\date
November 2021
\enddate 

\endtopmatter

\define\w{\widetilde}

\redefine\norm{\trianglelefteq}

\redefine\bar{\overline}

\redefine\maps{\mapsto}
\redefine\i{^{-1}}

\redefine\l{\lambda}
\redefine\s{\sigma}
\redefine\a{\alpha}
\redefine\b{\beta}

\redefine\g{\gamma}

\redefine\r{\rho}

\redefine\G{\Gamma}

\redefine\S{\Sigma}
\redefine\L{\Lambda}

\redefine\<{\langle}
\redefine\>{\rangle}
\redefine\F{\Bbb F}
\redefine\ca{\Cal}

\redefine\D{\Delta}

\redefine\sub{\subseteq} 

\redefine\nsub{\nsubseteq}

\redefine\bX{\bold X}

\redefine\nset{\emptyset} 

\redefine\1{\bold 1} 

\redefine\up{\uparrow}

\redefine\bw{\bold w} 
\redefine\bu{\bold u} 
\redefine\bv{\bold v} 
\redefine\bX{\bold X} 
\redefine\bY{\bold Y}

\document 
\vskip .1in 
\noindent 
{\bf Introduction} 
\vskip .1in 

This is Part II of the series that began with [Ch1] - and the reader is assumed to be familiar with Part I. 
References to results in Part I will be made by prefixing ``I" to the number of the cited result. Thus, 
for example, lemma I.3.12 is Stellmacher's Splitting Lemma. 

Let $(\ca L,\D,S)$ be a locality. There is then a category $\ca F=\ca F_S(\ca L)$ - the ``fusion system" of 
$\ca L$ - whose objects are the subgroups of $S$, and whose morphisms are compositions of conjugation 
maps $c_g:X\to Y$ from one subgroup of $S$ into another, induced by elements $g\in\ca L$. We say that 
$\ca L$ is a locality {\it on} $\ca F$. 

There is an extensive theory of abstract fusion systems and, more particularly, of saturated fusion 
systems. The references by Craven [Cr], and by Aschbacher, Kessar, and Oliver [AKO] provide far more 
material than will be needed here. It is possible to provide a completely self-contained treatment 
of the fusion systems associated with localities; and such a treatment will in fact be necessary when 
generalizing the results here (where $\ca L$ is finite) to {\it discrete localities} 
(where $\ca L$ is ``finite-dimensional"). Our approach here is to compromise, 
by presenting elementary material on fusion systems in a self-contained way in section 1, and to 
avail ourselves of the literature on saturated fusion systems in section 6. Readers who are already 
familiar with fusion systems, and who may be tempted to simply skip over section 1, are advised to 
familiarize themselves with the definition 1.8 of {\it radical} subgroups (not the usual definition, but 
equivalent to the usual definition if the ambient fusion system is saturated), and of {\it subcentric} 
subgroups. 

Throughout this paper, $(\ca L,\D,S)$ will be a {\it proper} locality on $\ca F$ (cf. 2.4). 
In brief, what this means is that $\D$ is ``not too small" (in that $\D$ contains the set $\ca F^{cr}$ of all 
subgroups of $S$ which are both centric and radical in $\ca F$), and that each of 
the groups $G=N_{\ca L}(P)$ for $P\in\D$ has the property that $C_G(O_p(G))\leq O_p(G)$. 

A non-empty collection $\G$ of subgroups of $S$ is {\it $\ca F$-closed} if $Q\in\G$ whenever there exists an 
$\ca F$-homomorphism $\phi:P\to Q$ for some $P\in\G$. For example, the set $\ca F^s$ of $\ca F$-subcentric 
subgroups is $\ca F$-closed by [He2]. 
It is a straightforward exercise with the definitions (see 2.11) to 
show that if $\D_0$ is an $\ca F$-closed subset of $\D$ containing $\ca F^{cr}$, then there is a 
unique proper locality $(\ca L_0,\D_0,S)$ on $\ca F$ such that the partial group $\ca L_0$ is a 
subset of $\ca L$ and such that the inclusion map $\ca L_0\to \ca L$ is a homomorphism of partial groups.  
We shall call $\ca L_0$ the {\it restriction} of $\ca L$ to $\D_0$. The core of this paper, consisting of 
sections 
3 through 5, concerns the opposite sort of operation, by which one expands, rather than restricts, 
the set of objects.

\proclaim {Theorem A1} Let $(\ca L,\D,S)$ be a proper locality on $\ca F$ and let $\D^+$ be an $\ca F$-closed 
collection of subgroups of $S$ such that $\D\sub\D^+\sub\ca F^s$. 
\roster 

\item "{(a)}" There exists a proper locality $(\ca L^+,\D^+,S)$ on $\ca F$ such that $\ca L$ is the 
restriction $\ca L^+\mid_\D$ of $\ca L^+$ to $\D$. Moreover, $\ca L^+$ is generated by $\ca L$ as a 
partial group. 

\item "{(b)}" For any proper locality $(\w{\ca L},\D^+,S)$ on $\ca F$ whose restriction to $\D$ is $\ca L$, 
there is a unique isomorphism $\ca L^+\to\w{\ca L}$ which restricts to the identity map on $\ca L$. 

\endroster 
\endproclaim

Recall from Part I 
that a partial subgroup $\ca N$ of a partial group $\ca L$ is {\it normal} in $\ca L$ (or is a 
{\it partial normal subgroup} of $\ca L$, denoted $\ca N\norm\ca L$) if $x^g:=g\i xg\in\ca N$ for all  
$x\in\ca N$ and $g\in\ca L$ for which the product $g\i x g$ is defined. 
Recall also: for any partial group $\ca L$ and any subset $X$ of $\ca L$, $\<X\>$ is defined to 
be the intersection of the set of partial subgroups of $\ca L$ containing $X$. The 
intersection of partial subgroups is again a partial subgroup by I.1.8, and $\<X\>$ is called 
the partial subgroup generated by $X$.

\proclaim {Theorem A2} Let the hypothesis and notation be as in Theorem A1, let $\ca N$ be a partial 
normal subgroup of $\ca L$, and set $T=S\cap\ca N$. Then there exists a unique partial normal subgroup 
$\ca N^+$ of $\ca L^+$ such that $\ca N=\ca N^+\cap\ca L$. Moreover, 
the mapping $\ca N\maps\ca N^+$ is an inclusion-preserving bijection from the set of partial normal 
subgroups of $\ca L$ to the set of partial normal subgroups of $\ca L^+$. 
\endproclaim 

As mentioned, section 6 concerns saturated fusion systems, and relies on results from the literature.  
We shall obtain the following result. 

\proclaim {Theorem B} Let $\ca F$ be the fusion system of a proper locality $(\ca L,\D,S)$. Then 
$\ca F$ is saturated. 
\endproclaim 
 
Section 6 contains also a treatment of the set $\ca F^s$ of subcentric subgroups of $S$, included also in  
Henke's work [He2] on this subject. 
Section 7 concerns the notions of $O^p_{\ca L}(\ca N)$ and $O^{p'}_{\ca L}(\ca N)$ for $\ca N$ a partial 
normal subgroup of a proper locality $\ca L$.

\vskip .2in 
\noindent 
{\bf Section 1: Fusion systems} 
\vskip .1in

We begin this section by providing a brief summary of some of the terminology and some of the basic results 
pertaining to general fusion systems. Some of the definitions to be given here will be non-standard, but 
will turn out to be equivalent to the standard definitions in the case of 
the fusion system of a proper locality. 
  
\definition {Definition 1.1} Let $S$ be a finite $p$-group. A {\it fusion system} $\ca F$ on $S$ is a 
category, whose set of objects is the set of subgroups of $S$, and whose morphisms satisfy the following 
conditions (in which $P$ and $Q$ are subgroups of $S$).  
\roster 

\item "{(1)}" Each $\ca F$-morphism $P\to Q$ is an injective homomorphism of groups. 

\item "{(2)}" If $g\in S$ and $P^g\leq Q$ then the conjugation map $c_g:P\to Q$ is an $\ca F$-morphism. 

\item "{(3)}" If $\phi:P\to Q$ is an $\ca F$-morphism then the bijection $P\to Im(\phi)$ defined by  
$\phi$ is an $\ca F$-isomorphism. 

\endroster 
\enddefinition  

One most often refers to $\ca F$-morphisms as $\ca F$-homomorphisms, in order to emphasize condition (1). 
Notice that (2) implies that all inclusion maps between subgroups of $S$ are $\ca F$-homomorphisms, and hence 
the restriction of an $\ca F$-homomorphism $P\to Q$ to a subgroup of $P$ is again an $\ca F$-homomorphism.  
\vskip .1in

Let $G$ be a finite group and let $S$ be a $p$-subgroup of $G$. There is then a fusion system 
$\ca F=\ca F_S(G)$ on $S$ in which the $\ca F$-homomorphisms $P\to Q$ are the maps 
$c_g:P\to Q$ given by conjugation by those elements $g\in G$ for which $P^g\leq Q$.

\definition {Definition 1.2} Let $\ca F$ be a fusion system on $S$, and let $\ca F'$ be a fusion system 
on $S'$. A homomorphism $\a:S\to S'$ is a {\it fusion-preseving} (relative to $\ca F$ and $\ca F'$) if, 
for each $\ca F$-homomorphism $\phi:P\to Q$, there exists an $\ca F'$-homomorphism $\psi:P\a\to Q\a$ such 
that $\a\mid_P\circ\psi=\phi\circ\a\mid_Q$. 
\enddefinition 

Notice that each of the $\ca F'$-homomorphisms $\psi$ in the preceding definition is uniquely determined, 
since all $\ca F'$-homomorphisms are injective. Thus, if $\a:S\to S'$ is a fusion-preserving homomorphism  
then $\a$ induces a mapping $Hom_{\ca F}(P,Q)\to Hom_{\ca F'}(P\a,Q\a)$, for each pair $(P,Q)$ of 
subgroups of $S$. The following result is then easily verifed. 

\proclaim {Lemma 1.3} Let $\ca F$ be a fusion system on $S$, let $\ca F'$ be a fusion system on $S'$, 
and let $\a:S\to S'$ be a fusion-preserving homomorphism. Then the mapping $P\maps P\a$ 
from objects of $\ca F$ to objects of $\ca F'$, together with the set of mappings 
$$ 
\a_{P,Q}:Hom_{\ca F}(P,Q)\to Hom_{\ca F'}(P\a,Q\a) \ \ (P,Q\leq S) 
$$ 
defines a functor $\a^*:\ca F\to\ca F'$. 
\qed 
\endproclaim

In view of the preceding result, a fusion-preserving homomorphism $\a$ may also be called a 
{\it homomorphism of fusion systems}. Notice that the inverse of a fusion-preserving isomorphism is 
fusion-preserving, and is therefore an isomorphism of fusion systems.  

In the special case where $S\leq S'$ and the inclusion map $S\to S'$ is fusion-preserving, we say that 
$\ca F$ is a {\it fusion subsystem} of $\ca F'$. Thus, $\ca F_S(S)$ is a fusion subsystem of $\ca F$ for 
each fusion system $\ca F$ on $S$, by 1.1(2). We refer to $\ca F_S(S)$ as the {\it trivial fusion system} 
on $S$. 
\vskip .1in

If $\ca F$ is a fusion system on $S$ and $P$ is a subgroup of $S$, write $P^{\ca F}$ 
for the set of subgroups of $S$ of the form $P\phi$, $\phi\in Hom_{\ca F}(P,S)$. The elements of 
$P^{\ca F}$ are the {\it $\ca F$-conjugates} of $P$. 

\vskip .1in 
For the remainder of this section let $\ca F$ be a fixed fusion system on the finite $p$-group $S$. 

\definition {Definition 1.4} Let $P\leq S$ be a subgroup of $S$. Then $P$ is {\it fully normalized} 
in $\ca F$ if $|N_S(P)|\geq|N_S(P')|$ for all $P'\in P^{\ca F}$. Similarly, $P$ is {\it fully 
centralized} in $\ca F$ if $|C_S(P)|\geq|C_S(P')|$ for all $P'\in P^{\ca F}$. 
\enddefinition

Let $U\leq S$ be a subgroup of $S$. The {\it normalizer} $N_{\ca F}(U)$ of $U$ in $\ca F$ is the category 
whose objects are the subgroups of $N_S(U)$, and whose morphisms $P\to Q$ ($P$ and $Q$ subgroups of 
$N_S(U)$) are restrictions of $\ca F$-homomorphisms $\phi:PU\to QU$ such that $U\phi=U$. Similarly, 
the {\it centralizer} $C_{\ca F}(U)$ of $U$ in $\ca F$ is the category whose objects are the subgroups 
of $C_S(U)$ and whose morphisms $\phi:P\to Q$ are restrictions of $\ca F$-homomorphisms $\phi:PU\to QU$ 
such that $\phi$ induces the identity map on $U$. One observes that $N_{\ca F}(U)$ is a fusion system on 
$N_S(U)$ and that $C_{\ca F}(U)$ is a fusion system on $C_S(U)$. The following result is immediate from 1.3.

\proclaim {Lemma 1.5} Let $U,V\leq S$ be subgroups of $S$, and suppose that there exists an 
$\ca F$-isomorphism $\a:N_S(U)\to N_S(V)$ such that $U\a=V$. Then $\a$ is an isomorphism 
$N_{\ca F}(U)\to N_{\ca F}(V)$. Similarly, if $\b:C_S(U)U\to C_S(V)V$ is an $\ca F$-isomorphism 
such that $U\b=V$, then the restriction of $\b$ to $C_S(U)$ is an isomorphism 
$C_{\ca F}(U)\to C_{\ca F}(V)$. 
\qed 
\endproclaim

\definition {Definition 1.6} Let $T$ be a subgroup of $S$. Then $T$ is {\it weakly closed} in $\ca F$ if 
$T^{\ca F}=\{T\}$, {\it strongly closed} in $\ca F$ if $X^{\ca F}$ is a set of subgroups of $T$ for each 
subgroup $X$ of $T$, and {\it normal} in $\ca F$ if $\ca F=N_{\ca F}(T)$. 
\enddefinition

The following result is immediate from the definitions.

\proclaim {Lemma 1.7} If $U$ and $V$ are subgroups of $S$ which are normal in $\ca F$ then also $UV$ 
is normal in $\ca F$. Thus, there is a largest subgroup $O_p(\ca F)$ of $S$ which is normal in $\ca F$. 
\qed 
\endproclaim 

A set $\D$ of subgroups of $S$ is {\it $\ca F$-invariant} if $X\in\D\implies X^{\ca F}\sub\D$. An 
$\ca F$-invariant set $\D$ of subgroups of $S$ is {\it $\ca F$-closed} if $\D$ is non-empty and 
is closed with respect to overgroups in $S$ ($P\in\D$ and $P\leq Q\leq S$ $\implies$ $Q\in\D$).

\definition {Definition 1.8} Let $\ca F$ be a fusion system on $S$, and let $P\leq S$ be a subgroup of $S$. 
\roster 

\item "{(1)}" $P$ is {\it centric} in $\ca F$ (or $P$ is {\it $\ca F$-centric}) if $C_S(Q)\leq Q$ for all 
$Q\in P^{\ca F}$. 

\item "{(2)}" $P$ is {\it radical} in $\ca F$ (or $P$ is {\it $\ca F$-radical}) if there exists 
$Q\in P^{\ca F}$ such that $Q$ is fully normalized in $\ca F$ and such that $Q=O_p(N_{\ca F}(Q))$. 

\item "{(3)}" $P$ is {\it quasicentric} in $\ca F$ (or $P$ is {\it $\ca F$-quasicentric}) if there exists 
$Q\in P^{\ca F}$ such that $Q$ is fully centralized in $\ca F$ and such that $C_{\ca F}(Q)$ is the 
trivial fusion system on $C_S(Q)$. 

\item "{(4)}" $P$ is {\it subcentric} in $\ca F$ (or $P$ is {\it $\ca F$-subcentric}) if there exists 
$Q\in P^{\ca F}$ such that $Q$ is fully normalized in $\ca F$ and such that $O_p(N_{\ca F}(Q))$ is 
centric in $\ca F$. 

\endroster 
Write $\ca F^c$, $\ca F^q$, and $\ca F^s$, respectively, for the set of subgroups of $S$ which 
are $\ca F$-centric, $\ca F$-quasicentric, and $\ca F$-subcentric. Write $\ca F^{cr}$ for the set of 
subgroups of $S$ which are both centric and radical in $\ca F$. 
\enddefinition 

\definition {Remark} The above definition of $\ca F$-radical subgroup is different from the standard 
one (which is that $P$ is $\ca F$-radical if $Inn(P)=O_p(Aut_{\ca F}(P))$). But it will turn out to be 
equivalent to the standard definition in the case that $\ca F$ is the fusion system of a proper locality. 
\enddefinition 

\proclaim {Lemma 1.9} Let $\ca F$ be a fusion system on $S$. Then $\ca F^{cr}$ is $\ca F$-invariant, and 
$\ca F^c$ is $\ca F$-closed. 
\endproclaim 

\demo {Proof} Both $\ca F^c$ and $\ca F^{cr}$ are $\ca F$-invariant by definition.  
Let $P\in\ca F^c$, let $P\leq Q\leq S$, and let $\phi:Q\to S$ be an $\ca F$-homomorphism. 
Then $C_S(Q\phi)\leq C_S(P\phi)\leq P\phi\leq Q\phi$, and so $Q\in\ca F^c$. Thus $\ca F^c$ is closed 
with respect to overgroups in $S$. As $S\in\ca F^c$ it follows that $\ca F^c$ is $\ca F$-closed. 
\qed 
\enddemo

\proclaim {Lemma 1.10} Let $P\leq S$ be a subgroup of $S$ and let $Q\in P^{\ca F}$ such that $Q$ is 
fully centralized in $\ca F$. Then $P\in\ca F^c$ if and only if $C_S(Q)\leq Q$. 
\endproclaim 

\demo {Proof} Suppose that $C_S(Q)\leq Q$ and let $R\in Q^{\ca F}$. Then 
$$ 
|C_S(R)|\leq |C_S(Q)|=|Z(Q)|=|Z(R)|, 
$$ 
and so $C_S(R)=Z(R)$. That is, $C_S(R)\leq R$, and thus $Q$ is $\ca F$-centric. As $\ca F^c$ is 
$\ca F$-invariant by 1.9, $P$ is then $\ca F$-centric. That is: 
$$ 
C_S(Q)\leq Q\implies P\in\ca F^c. 
$$ 
The reverse implication is given by the definition of $\ca F^c$. 
\qed 
\enddemo 

Let $\Psi$ be a non-empty set of $\ca F$-isomorphisms. Then $\ca F$ is {\it generated} by $\Psi$ if 
every $\ca F$-isomorphism can be expressed as a composition of restrictions of members of $\Psi$, and  
write $\ca F=\<\Psi\>$ in that case. If $\G$ is a non-empty set of subgroups of $S$ and 
$$ 
\Psi=\<\bigcup\{Aut_{\ca F}(R)\mid R\in\G\>. \tag* 
$$ 
then we say that $\ca F$ is {\it $\G$-generated} if $\ca F=\<\Psi\>$. An important special case is that 
in which $\G$ is the set of all $P\ca F^{cr}$ such that $P$ is fully normalized in $\ca F$, where we 
we shall say that $\ca F$ is {\it $(cr)$-generated} if $\ca F$ is $\G$-generated.

\definition {Definition 1.11} Let $\ca F$ be a fusion system on $S$, and let $\G$ be an $\ca F$-closed 
set of subgroups of $S$. Then $\ca F$ is {\it $\G$-inductive} if: 
\roster 

\item "{(*)}" For each $U\in\G$, and each $V\in U^{\ca F}$ such that $V$ is fully normalized in $\ca F$,  
there exists an $\ca F$-homomorphism $\phi:N_S(U)\to N_S(V)$ with $U\phi=V$. 

\endroster 
If $\ca F$ is $\G$-inductive where $\G$ is the set of all subgroups of $S$, we shall simply say that $\ca F$  
is {\it inductive}. 
\enddefinition

\proclaim {Lemma 1.12} Let $G$ be a finite group, let $S$ be a Sylow $p$-subgroup of $G$, and set 
$\ca F=\ca F_S(G)$. Then $\ca F$ is inductive, and $\ca F$ is $(cr)$-generated.  
\endproclaim 

\demo {Proof} Let $V$ be fully normalized in $\ca F$. Equivalently: $N_S(V)\in Syl_p(N_G(V))$. Let 
$U\in V^{\ca F}$, and let $g\in\ca L$ with $U^g=V$. Then $N_S(U)^g\leq N_G(V)$, and there then exists 
$h\in N_G(V)$ with $N_S(U)^{gh}\leq N_S(V)$. Thus $\ca F$ is inductive. That $\ca F$ is $(cr)$-generated 
is a well-known consequence of the Alperin-Goldschidt fusion theorem [Gold].
\qed 
\enddemo

The next few results provide information about inductive fusion systems. All of these results 
can be re-stated (and proved) in an obvious way, so as to yield corresponding results about $\G$-inductive 
fusion systems.

\proclaim {Lemma 1.13} Assume that $\ca F$ is inductive, and let $V\leq S$ be fully normalized in $\ca F$. 
Then $V$ is fully centralized in $\ca F$. 
\endproclaim 

\demo {Proof} Let $U\in V^{\ca F}$, and let $\phi:N_S(U)\to N_S(V)$ be an $\ca F$-homomorphism such 
that $V=U\phi$. Then $\phi$ maps $C_S(U)$ into $C_S(V)$, and thus $|C_S(U)|\leq |C_S(V)|$. 
\qed 
\enddemo

\proclaim {Lemma 1.14} Assume that $\ca F$ is inductive, and let $U$ and $V$ be $\ca F$-conjugate 
subgroups of $S$. 
\roster 

\item "{(a)}" If $U$ and $V$ are fully normalized in $\ca F$ then $N_{\ca F}(U)\cong N_{\ca F}(V)$.  

\item "{(b)}" If $U$ and $V$ are fully centralized in $\ca F$ then $C_{\ca F}(U)\cong C_{\ca F}(V)$.  

\endroster  
\endproclaim 

\demo {Proof} Suppose that $U$ and $V$ are fully normalized in $\ca F$. Then (FL1) implies that there 
exisits an $\ca F$-isomorphism $\phi:N_S(U)\to N_S(V)$ with $U\phi=V$. By 1.5, $\phi$ is then an isomorphism 
$N_{\ca F}(U)\to N_{\ca F}(V)$. Thus (a) holds. Now suppose instead that $U$ and $V$ are fully centralized 
in $\ca F$, and let $X$ be a fully normalized $\ca F$-conjugate of $U$ (and hence also of $V$). There 
are then $\ca F$-homomorphisms $\r:N_S(U)\to N_S(X)$ and $\s:N_S(V)\to N_S(X)$ with $U\r=X=V\s$. The 
restriction of $\r$ to $C_S(U)U$ is then an isomorphism with $C_S(X)X$, and similarly $\s$ restricts to an 
isomorphism $C_S(V)V\to C_S(X)X$. A further application of 1.5  now yields (b). 
\qed 
\enddemo

\proclaim {Lemma 1.15} Assume that $\ca F$ is inductive, let $T$ be strongly closed in 
$\ca F$, let $Q\leq S$ be a subgroup of $S$, and set $V=Q\cap T$. Suppose that $V$ is fully normalized 
in $\ca F$ and that $Q$ is fully normalized in $N_{\ca F}(V)$. Then $Q$ is fully normalized in $\ca F$.  
\endproclaim 

\demo {Proof} Let $P\in Q^{\ca F}$ such that $P$ is fully normalized in $\ca F$, let $\phi:N_S(Q)\to N_S(P)$ 
be an $\ca F$-homomorphism with $Q\phi=P$, and set $U=V\phi$. Then $U=P\cap T$ as $T$ is strongly closed 
in $\ca F$, and then also $N_S(P)\leq N_S(U)$. By (FL1) there exists an $\ca F$-homomorphism 
$\psi:N_S(U)\to N_S(V)$ with $U\psi=V$. Set $Q'=P\psi$. Then $Q'=(Q\phi)\psi$ is an 
$N_{\ca F}(V)$-conjugate of $Q$. Since $N_S(Q)=N_{N_S(V)}(Q)$ (and similarly for $Q'$), and since $Q$ 
is fully normalized in $N_{\ca F}(V)$, we have $|N_S(Q)|\geq|N_S(Q')|$. The sequence 
$$ 
N_S(Q)@>\phi>>N_S(P)@>\psi>>N_S(Q')  
$$ 
of injective homomorphisms then shows that $|N_S(P)|=|N_S(Q)|$, and so $Q$ is fully normalized in $\ca F$. 
\qed 
\enddemo 

\proclaim {Lemma 1.16} Assume that $\ca F$ is inductive, let $V\leq S$ be a subgroup of $S$, and let 
$Q\leq N_S(V)$ be a subgroup of $N_S(V)$ containing $V$. Then each of the following conditions implies 
that $Q$ is fully centralized in $\ca F$. 
\roster 

\item "{(1)}" $V$ is fully normalized in $\ca F$ and $Q$ is fully centralized in $N_{\ca F}(V)$. 

\item "{(2)}" $V$ is fully centralized in $\ca F$, $Q=VY$ for some subgroup $Y$ of $C_S(V)$ such that 
$Y$ is fully centralized in $N_{\ca F}(V)$, 
and $N_{\ca F}(V')$ is inductive for each fully normalized $\ca F$-conjugate $V'$ of $V$.   

\endroster 
\endproclaim 

\demo {Proof} Assume first that (1) holds. 
As in the proof of 1.15 let $P\in Q^{\ca F}$ such that $P$ is fully normalized in $\ca F$,  
and let $\phi:N_S(Q)\to N_S(P)$ be an $\ca F$-homomorphism with $Q\phi=P$. Set $U=V\phi$ and let 
$\psi:N_S(U)\to N_S(V)$ be an $\ca F$-homomorphism with $U\psi=V$. Set $Q'=P\psi$. Then $Q'$ is an 
$N_{\ca F}(V)$-conjugate of $Q$. As $Q$ is fully centralized in $N_{\ca F}(V)$, and since 
$C_S(Q)=C_{C_S(V)}(Q)$ (and similarly for $Q'$), we then have $|C_S(Q)|\geq|C_S(Q')|$. The sequence 
$$ 
C_S(Q)@>\phi>>C_S(P)@>\psi>>C_S(Q')
$$ 
then shows that $|C_S(P)|=|C_S(Q)|$. On the other hand, as $P$ is fully normalized in $\ca F$, 
$P$ is also fully centralized in $\ca F$ by 1.13. Thus $|C_S(P)|\geq|C_S(Q)|$, and $Q$ is fully 
centralized in $\ca F$. 

Assume next that (2) holds. Let $V'\in V^{\ca F}$ with $V'$ fully normalized in $\ca F$, and let 
$\a:C_S(V)V\to C_S(V')V'$ be an $\ca F$-homomorphism with $V'=V\a$. Then $\a$ is an 
isomorphism, and induces an isomorphism $C_{\ca F}(V)\to C_{\ca F}(V')$. Thus $Y'$ is fully centralized 
in $C_{\ca F}(V')$. Here $\phi\i$ maps $C_S(V'Y')$ into $C_S(VY)$, so if 
$V'Y'$ is fully centralized in $\ca F$ then also $VY$ is fully centralized in $\ca F$. Thus we may 
assume from the beginning that $V$ is fully normalized in $\ca F$. By hypothesis, $N_{\ca F}(V)$ is 
then inductive.  

Let $P\in Q^{\ca F}$ be fully centralized in $\ca F$, and let $\phi:C_S(Q)Q\to C_S(P)P$ be an 
$\ca F$-homomorphism with $P=Q\phi$. Set $U=V\phi$ and $X=Y\phi$. As $V$ is fully centralized 
in $\ca F$ there exists an $\ca F$-homomorphism $\psi:C_S(P)P\to C_S(V)V$ with $U\psi=V$. 
Set $R=P\psi$. Then $R$ is an $N_{\ca F}(V)$-conjugate of $Q$ and $R$ is fully centralized in $\ca F$, 
whence also $R$ is fully centralized in $N_{\ca F}(V)$. We have $R=VZ$ where $Z=(Y)\phi\circ\psi$ is 
an $N_{\ca F}(V)$-conjugate of $Y$. As $Y$ is fully centralized in $N_{\ca F}(V)$ there exists an 
$N_{\ca F}(V)$-homomorphism $\r:C_S(VZ)VZ\to C_S(VY)VY$ with $Z\r=V$. Thus $C_S(P)P$ is mapped by 
$\psi\circ\r$ into $C_S(Q)Q$, and so $Q$ is fully centralized in $\ca F$. 
\qed 
\enddemo 

\proclaim {Lemma 1.17} Assume that $\ca F$ is inductive, and let $T\leq S$ be strongly closed in $\ca F$. 
\roster 

\item "{(a)}" For any $P\leq S$ there exists $Q\in P^{\ca F}$ such that both $Q$ and $Q\cap T$ are fully 
normalized in $\ca F$. 

\item "{(b)}" Let $\ca E$ be a fusion subsystem of $\ca F$ on $T$, and let $U\leq T$ be a subgroup of $T$ 
such that $U$ is fully normalized in $\ca F$. Then $U$ is fully normalized in $\ca E$. 

\endroster 
\endproclaim 

\demo {Proof} For the proof of (a), set $X=P\cap T$, let $Y\in X^{\ca F}$ be fully normalized in $\ca F$, 
and let $\phi:N_S(X)\to N_S(Y)$ be an $\ca F$-homomorphism with $X\phi=Y$. Without loss of generality we 
may assume that $P$ has itself been chosen so that $P$ is fully normalized in $\ca F$. Set $Q=P\phi$. 
We have $N_S(P)\leq N_S(X)$, so $\phi$ restricts to an isomorphism $N_S(P)\to N_S(Q)$. Thus $Q$ is 
fully normalized in $\ca F$, and (a) holds. 

For the proof of (b), 
let $V\in U^{\ca E}$ such that $V$ is fully normalized in $\ca E$. Then $V\in U^{\ca F}$. 
As $\ca F$ is inductive there exists an $\ca F$-homomorphism $\phi:N_S(V)\to N_S(U)$ with $V\phi=U$. 
Then $N_T(V)\phi\leq N_T(U)$ as $T$ is strongly closed in $\ca F$. As $V$ is fully normalized in $\ca E$ 
it follows that $N_T(V)\phi=N_T(U)$. Thus $U$ is fully normalized in $\ca E$, and (b) holds.  
\qed 
\enddemo

\proclaim {Lemma 1.18} Assume that $\ca F$ is inductive, and let $U\leq S$ be a subgroup of $S$. Then there 
exists $V\in U^{\ca F}$ such that both $V$ and $O_p(N_{\ca F}(V))$ are fully normalized in $\ca F$.  
\endproclaim 

\demo {Proof} Without loss of generality, $U$ is fully normalized in $\ca F$. Set $P=O_p(N_{\ca F}(U))$ 
and let $Q\in P^{\ca F}$ with $Q$ fully normalized in $\ca F$. Let $\phi:N_S(P)\to N_S(Q)$ be an 
$\ca F$-homomorphism which maps $P$ to $Q$, and set $V=U\phi$. As $N_S(U)\leq N_S(P)$ we have 
$N_S(U)\phi\leq N_S(V)$. But $|N_S(U)|\geq|N_S(V)|$ as $U$ is fully normalized, so $\phi$ induces an 
isomorphism $N_S(U)\to N_S(V)$. Thus $V$ is fully normalized in $\ca F$, and then $Q=O_p(N_{\ca F}(V))$ 
by 1.14(a). 
\qed 
\enddemo 

\proclaim {Lemma 1.19} Let $\ca F$ be an inductive fusion system on $S$, and let $V\leq S$ be a subgroup 
of $S$. Assume that $V$ is fully normalized in $\ca F$, and  
let $P$ be a subgroup of $N_S(V)$ containing $V$. Then $P$ is centric in $\ca F$ if and 
only if $P$ is centric in $N_{\ca F}(V)$. 
\endproclaim 

\demo {Proof} We are free to replace $P$ by any $N_{\ca F}(V)$-conjugate of $P$, 
so we may assume that $P$ is fully normalized in $N_{\ca F}(V)$. Then $P$ is  
fully centralized in $\ca F$ by 1.16, and so $P\in\ca F^c$ if and only if $C_S(P)\leq P$ by 1.10. As 
$V\leq P$ we have $C_S(P)=C_{N_S(V)}(P)$, and the lemma follows. 
\qed 
\enddemo 

THe next result concerns quotients of fusion systems. The hypotheses will be seen to be 
fulfilled in 2.12, in the contest of projections of localities.

\proclaim {Lemma 1.20} Let $\ca F$ and $\bar{\ca F}$ be fusion systems on $p$-groups $S$ and 
$\bar S$, and let $\l:S\to\bar S$ be a fusion-preserving homomomorphism. Denote also by $\l$ the 
corresponding homomorphism $\ca F\to\bar{\ca F}$ of fusion systems given by 1.3, and write $\bar X$ for 
the image under $\l$ of any subgroup or set of subgroups $X$ of $S$. Assume that $\l$ is 
surjective, and that each of the mappings 
$$ 
\l_{X,Y}:Hom_{\ca F}(X,Y)\to Hom_{\bar{\ca F}}(\bar X,\bar Y)\ \ \ (Ker(\l)\leq X\cap Y) 
$$ 
is surjective. Let $\bar P\leq\bar S$ and let $P$ be the $\l$-preimage of $\bar P$. Then the following 
hold. 
\roster 

\item "{(a)}" $P$ is fully normalized in $\ca F$ if and only if $\bar P$ is fully normalized in 
$\bar{\ca F}$. 

\item "{(b)}" $O_p(N_{\ca F}(P))\l\leq O_p(N_{\bar{\ca F}}(\bar P))$. 

\item "{(c)}" For each of the symbols $\star\in\{cr,c,q,s\}$ we have $\bar P\in\bar{\ca F}^*$ if 
$P\in\ca F^*$.  

\endroster 
\endproclaim 

\demo {Proof} For any subgroup or element $\bar U$ of $\bar S$ write $U$ for the 
preimage of $\bar U$ in $S$. For any such $U$ we have $\bar{N_S(U)}=N_{\bar S}(\bar U)$, and it is this 
observation that yields (a). 

Set $Q=O_p(N_{\ca F}(P))$ and set $\bar R=O_p(N_{\bar{\ca F}}(\bar P)$. Let $\bar\phi:\bar X\to\bar Y$ be a 
$N_{\bar{\ca F}}(\bar P)$-homomorphism. By hypothesis there exists an $N_{\ca F}(P)$-homomorphism 
$\phi:X\to Y$ with $\bar\phi=(\phi)\l$. Then $\phi$ extends to an $N_{\ca F}(P)$-homomorphism 
$\psi:QX\to QY$ which fixes $Q$, and then $(\psi)\l$ is an extension of $\bar\phi$ to an 
$N_{\bar{\ca F}}(\bar P)$-homomorphism $\bar Q\bar X\to\bar Q\bar Y$ which fixes $\bar Q$. This shows that 
$\bar Q\norm N_{\bar{\ca F}}(\bar P)$, and so $\bar Q\leq\bar R$. A similar argument 
shows that $R\leq Q$, and establishes (b). 

The statement in (c) concerning $\bar{\ca F}^c$ and $\ca F^c$ is immediate from the observation that 
$\bar{C_S(U)}\leq C_{\bar S}(\bar U)$ for any subgroup $U\leq S$. For any fusion system $\ca E$ on a 
$p$-group $T$ write $\ca E^r$ for the set of $\ca E$-radical subgroups of $T$ (cf. 1.8). 
Suppose that $\bar P\in\bar{\ca F}^r$. By definition 1.8 there exists  
an $N_{\bar{\ca F}}(\bar P)$-conjugate $\bar P_1$ of $\bar P$ such that $\bar P_1$ is fully normalized in 
$\bar{\ca F}$ and such that $\bar P_1=O_p(N_{\bar{\ca F}}(\bar P_1))$. Then $P_1$ is an $\ca F$-conjugate 
of $P$, $P_1$ is fully normalized in $\ca F$ by (a), and $P_1=O_p(N_{\ca F}(P_1))$ by (b). Thus 
$P\in\ca F^r$, and this establishes the statement in (c) concerning ``$cr$". 
The statement concerning ``$s$" follows from (a) and 
from the statement concerning ``$c$". Finally, if $C_{\bar{\ca F}}(\bar P)$ is 
the trivial fusion system on $C_{\bar S}(\bar P)$ then $C_{\ca F}(P)$ is the trivial fusion system on 
$C_S(P)$, and then (a) completes the proof of the statement in (c) 
concerning ``$q$". 
\qed 
\enddemo 

\proclaim {Lemma 1.21} Assume that $\ca F$ is inductive and that $N_{\ca F}(X)$ is inductive for all 
$X\leq S$ such that $X$ is fully normalized in $\ca F$. Let $A$ and $B$ be subgroups of $S$ 
such that $[A,B]\leq A\cap B$, and assume: 
\roster 

\item "{(*)}" $A$ is fully normalized in $\ca F$ and $B$ is fully normalized in $N_{\ca F}(A)$. 

\endroster 
Let $\bar B\in B^{\ca F}$ be fully normalized in $\ca F$, let $\phi:N_S(B)\to N_S(\bar B)$ be an 
$\ca F$-homomorphism such that $\bar B=B\phi$, and set $\bar A=A\phi$. Set $D=N_S(A)\cap N_S(B)$, 
set $\bar D=N_S(\bar A)\cap N_S(\bar B)$, and let $\psi$ be the homomorphism $D\to\bar D$ induced by $\phi$. 
Then $\psi$ is an isomorphism, $\bar A$ is fully normalized in $N_{\ca F}(\bar B)$, and $\psi$ is then  
an isomorphism 
$$ 
N_{N_{\ca F}(A)}(B)\cong N_{N_{\ca F}(\bar B)}(\bar A). 
$$  
\endproclaim 

\demo {Proof} Let $\bar A'$ be an $N_{\ca F}(\bar B)$-conjugate of $\bar A$ such that $\bar A'$ is fully 
normalized in $N_{\ca F}(\bar B)$. Set $\bar D'=N_S(\bar A)\cap N_S(\bar B)$, and 
let $\eta:D\phi\to\bar D'$ be an $N_{\ca F}(\bar B)$-homomorphism such that $\bar A\eta=\bar A'$. 
By (*) and the hypothesis concerning inductivity there exists an $\ca F$-homomorphism $\l:\bar D'\to D$ with 
$\bar A'\l=A$ and with $\bar B\l=B$. Thus $\psi$ and $\eta$ are isomorphisms, and so $\bar A$ is fully 
normalized in $N_{\ca F}(\bar B)$. In order to complete the proof it remains to show that 
$\psi$ is fusion-preserving. This is achieved by the observation that if $\a:U\to V$ is an 
$N_{N_{\ca F}(A)}(B)$-isomorphism with 
$AB\leq U$. Then $\psi\i\circ\a\circ\psi$ is an $N_{N_{\ca F}(\bar B)}(\bar A)$-isomorphism 
$U\psi\to V\psi$. 
\qed 
\enddemo

\vskip .2in 
\noindent 
{\bf Section 2: Fusion systems of localities} 
\vskip .1in 

Fix a locality $(\ca L,\D,S)$ throughout this section. For each $g\in\ca L$ there is then an isomorphism  
$$ 
c_g:S_g\to S_{g\i}
$$ 
given by conjugation by $g$. Define $\ca F_S(\ca L)$ to be the fusion system $\ca F$ on $S$ which is generated 
by these isomorphisms. We then say that $(\ca L,\D,S)$ is a locality {\it on} $\ca F$. 

It will be useful to have the following explicit description of the set of $\ca F$-isomorphisms. For 
$w=(g_1,cdots,g_n)\in\bold W(\ca L)$ and for $U\leq S_w$ define $V=U^w$ to be the subgroup of $S_{w\i}$ 
obtained by successively applying the conjugation maps $c_{g_i}$ to $U$: 
$$ 
U=U_0@>c_{g_1}>>U_1@>c_{g_2}>>\cdots@>c_{g_n}>>U_n=V.  
$$ 
The set of such mappings $c_w:U\to V$ ($U,V\leq S$) is then the set of $\ca F$-isomorphisms.

\proclaim {Lemma 2.1} Let $(\ca L,\D,S)$ be a locality on $\ca F$, and let $P\in\D$. Then 
$$ 
P^{\ca F}=\{P^g\mid g\in\ca L,\ P\leq S_g\}. \tag*
$$ 
Moreover, $P$ is fully normalized in $\ca F$ if and only if $N_S(P)\in Syl_p(N_{\ca L}(P))$, and $P$ is 
fully centralized in $\ca F$ if and only if $C_S(P)\in Syl_p(C_{\ca L}(P))$. 
\endproclaim 

\demo {Proof} Let $\phi:P\to Q$ be an $\ca F$-isomorphism. As noted above, $\phi=c_w$ for some 
$w\in\bold W(\ca L)$ with $P\leq S_w$. Then $w\in\bold D$, and $c_w=c_{\Pi(w)}$ by I.2.3(c). This yields 
(*). 

As $P$ is in $\D$, $N_{\ca L}(P)$ is a subgroup of $\ca L$, and $C_{\ca L}(P)$ is a normal 
subgroup of $N_{\ca L}(P)$. Let $X$ be a Sylow $p$-subgroup of $N_{\ca L}(P)$ containing $N_S(P)$. By 
I.2.11 there exists $g\in\ca L$ with $X^g\leq S$, and conjugation by $g$ induces an isomorphism 
$N_{\ca L}(P)\to N_{\ca L}(P^g)$ by I.2.3(b). If $P$ is fully normalized in $\ca F$ then 
$|N_S(P)|\geq |N_S(P^g)|\geq |X^g|=|X|$, so $N_S(P)=X$, and $N_S(P)\in Syl_p(N_{\ca L}(P)$ in that case. 
Conversely, suppose that $N_S(P)$ is a Sylow subgroup of $N_{\ca L}(P)$ and let $Q\in P^{\ca F}$. We have   
$P=Q^h$ for some $h\in\ca L$ by (*), so $N_S(Q)^h\leq N_{\ca L}(P)$. By Sylow's theorem there exists 
$f\in N_{\ca L}(P)$ with $(N_S(Q)^g)^f\leq N_S(P)$, and thus $|N_S(P)|\geq |N_S(Q)|$. This establishes 
the first of the two ``if and only ifs" of the lemma. The proof of the second ``if and only if" is obtained 
in similar fashion.  
\qed 
\enddemo

\proclaim {Proposition 2.2} Let $(\ca L,\D,S)$ be a locality on $\ca F$. Then $\ca F$ is $\D$-inductive. 
Moreover, for each $P\in\D$ such that $P$ is fully normalized in $\ca F$: 
\roster

\item "{(a)}" $N_{\ca F}(P)=\ca F_{N_S(P)}(N_{\ca L}(P))$ and $C_{\ca F}(P)=\ca F_{C_S(P)}(C_{\ca L}(P))$. 

\item "{(b)}" $N_{\ca F}(P)$ and $C_{\ca F}(P)$ are $(cr)$-generated. 

\endroster 
\endproclaim 

\demo {Proof} That $\ca F$ is $\D$-inductive is immediate from the preceding lemma and from I.2.10. Let 
$P\in\D$ with $P$ fully normalized in $\ca F$, and let $\phi:X\to Y$ be an $N_{\ca F}(P)$-isomorphism 
between two subgroups $X$ and $Y$ of $N_S(P)$ containing $P$. As in the proof of (*) in 2.1, we find that 
$\phi=c_g$ for some $g\in N_{\ca L}(P)$, and this shows that $N_{\ca F}(P)=\ca F_{N_S(P)}(N_{\ca L}(P))$. 
Then $N_{\ca F}(P)$ is $N_{\ca F}(P)^{cr}$-generated by 1.12. 

By 1.13 $P$ is fully centralized in $\ca F$. Let $\phi:X\to Y$ be a $C_{\ca F}(P)$-isomorphism 
between two subgroups $X$ and $Y$ of $C_S(P)$. By definition of $C_{\ca F}(P)$, $\phi$ extends to an 
$\ca F$-isomorphism $\psi:XP\to YP$ such that $\psi$ restricts to the identity map on $P$. Then 
$\psi=c_g$ for some $g\in C_{\ca L}(P)$, and thus $C_{\ca F}(P)=\ca F_{C_S(P)}(C_{\ca L}(P))$. We again 
appeal to 1.12, obtaining $(cr)$-generation for $C_{\ca F}(P)$. 
\qed 
\enddemo

A finite group $G$ is of {\it characteristic $p$} if $C_G(O_p(G))\leq O_p(G)$.

\proclaim {Lemma 2.3} Let $(\ca L,\D,S)$ be a locality on $\ca F$, and assume that $N_{\ca L}(P)$ is of 
characteristic $p$ for all $P\in\D$. Let $R\leq S$ be a subgroup of $S$. Then $R\norm\ca L$ if and only 
if $Z\norm\ca F$. In particular, we have $O_p(\ca F)=O_p(\ca L)$. 
\endproclaim 

\demo {Proof} As $\ca F$ is generated by the conjugation maps $c_g:S_g\to S$ for $g\in\ca L$, it is 
immediate that if $R\norm\ca L$ then $R\norm\ca F$. So assume that $R\norm\ca F$, and assume 
by way of contradiction that $R$ is not normal in $\ca L$. Among all 
elements of $\ca L$ not in $N_{\ca L}(R)$, choose $g$ so that $|S_g|$ is as large as possible. Then 
$R\nleq S_g$, since $c_g:S_g\to S$ is an $\ca F$-homomorphism. In particular, $S_g\neq S$, and $S_g$ is a 
proper subgroup of $N_S(S_g)$. 

Set $P=S_g$, $P'=P^g$, and let $Q\in P^{\ca F}$ be fully normalized in $\ca F$. As $\ca L$ is 
$\D$-inductive by 2.2, 
there exists $x\in\ca L$ such that $P^x=Q$ and such that $N_S(P)\leq S_x$. Then also $R\leq S_x$ by the 
maximality of $|P|$ in the choice of $g$. Since $Q\in(P')^{\ca F}$ there exists also $y\in\ca L$ 
such that $(P')^y=Q$ and such that $N_S(P')R\leq S_y$. 

Note that $(x\i,g,y)\in\bold D$ via $Q$, and that $f:=\Pi(x\i,g,y)\in N_{\ca L}(Q)$. Note also that 
$(x,x\i,g,y,y\i)\in\bold D$ via $P$, and that 
$$ 
\Pi(x,f,y\i)=\Pi(x,x\i,g,y,y\i)=\Pi(g) 
$$ 
by $\bold D$-associativity (I.1.4). If $R\leq S_f$ then $R\leq S_{(x,f,y\i)}$, and then $R\leq S_g$. 
Thus $R\nleq S_f$, 
and we may therefore replace $g$ with $f$, and $P$ with $Q$. That is, we may assume that  
$P$ is fully normalized in $\ca F$ and $g\in N_{\ca L}(P)$. 

Set $M=N_{\ca L}(P)$. Then $M$ is a subgroup of $\ca L$ and, by hypothesis, $M$ is of characteristic $p$. 
Set $D=N_R(P)$, and set $\ca E=\ca F_{N_S(P)}(M)$. Then $D\norm\ca E$, and so the conjugation map 
$c_g:P\to P$ extends to an $\ca E$-automorphism of $PD$. Thus, there exists $h\in N_M(PD)$ such 
that $c_h:PD\to PD$ restricts to $c_g$ on $P$. Then $gh\i\in C_M(P)\leq P$, and so $gh\i\in N_M(D)$. 
This yields $D^g=D$, so $D\leq P$, and then $R\leq P$. Then $R^g=R$. This result is 
contrary to the choice of $g$, and completes the proof. 
\qed 
\enddemo

\definition {Definition 2.4} Let $(\ca L,\D,S)$ be a locality on $\ca F$. Then $\ca L$ is {\it proper} if: 
\roster 

\item "{(PL1)}" $\ca F^{cr}\sub\D$, and 
\vskip .05in 

\item "{(PL2)}" $N_{\ca L}(P)$ is of characteristic $p$ for each $P\in\D$. 

\endroster 
\enddefinition

The next two results are well known, and are important for an understanding of the structure of finite 
groups of characteristic $p$.

\proclaim {Lemma 2.5} Let $G$ be a finite group, and let $A$ and $B$ be subgroups of $G$ such that 
$|A|$ is relatively prime to $|B|$. Suppose that $[A,B]\leq C_B(A)$. Then $[A,B]=1$. 
\endproclaim 

\demo {Proof} Let $a\in A$ and $b\in B$. Then $a\i a^b=[a,b]$ commutes with $a$, by hypothesis, and 
so $a^b$ commutes with $\<a\>$. As $|a|=|a^b|$ is relatively prime to $|B|$, the same is then true of 
$a\i a^b$. As $[a,b]\in B$ by hypothesis, we conclude that $[a,b]=1$, and thus $[A,B]=1$. 
\qed 
\enddemo

\proclaim {Lemma 2.6 (Thompson's $A\times P$ Lemma)} Let $G$ be a finite group and let $H\leq G$ be a 
subgroup of $G$ such that: 
\roster 

\item "{(i)}" $H=PA$, where $P$ is a normal $p$-subgroup of $H$, and where $A$ is a $p'$-subgroup of $H$. 

\item "{(ii)}" There exists a subgroup $B$ of $C_P(A)$ such that $[C_P(B),A]=1$. 

\endroster 
Then $[P,A]=1$ and $H$ is isomorphic to the direct product $A\times P$. 
\endproclaim 

\demo {Proof} We have $A\cap P=1$ since $|A|$ and $|P|$ are relatively prime, and so it suffices to 
show that $P=C_P(A)$. Suppose false, so that $C_P(A)$ is a proper subgroup of $P$. Set $Q=N_P(C_P(A))$. 
Thus $B\leq C_P(A)\norm Q$, and so $[Q,B]\leq C_P(A)$. One may express this in the standard way as 
$[Q,B,A]=1$. Also $[B,A,Q]=1$ since $B\leq C_P(A)$. The Three Subgroups Lemma (2.2.3 in [Gor]) 
then yields $[A,Q,B]=1$. That is, we have $[A,Q]\leq C_Q(B)$, and hence $[A,Q]\leq C_Q(A)$ by (ii). Then 
$Q\leq C_P(A)$ by 2.5. That is, $N_P(C_P(A))\leq C_P(A)$, and hence $P=C_P(A)$, as required. 
\qed 
\enddemo

\proclaim {Lemma 2.7} Let $G$ be a finite group of characteristic $p$. 
\roster 

\item "{(a)}" Every normal subgroup of $G$ is of characteristic $p$. 

\item "{(b)}" Every $p$-local subgroup of $G$ is of characteristic $p$. 

\item "{(c)}" Let $V$ be a normal $p$-subgroup of $G$, and let $X$ be the set of elements $x\in C_G(V)$ 
such that $[O_p(G),x]\leq V$. Then $X$ is a normal $p$-subgroup of $G$. 

\item "{(d)}" Let $N\norm G$ be a normal subgroup such that $C_G(N)\leq N$. Then $C_G(O_p(N))\leq O_p(G)$, 
and every overgroup of $N$ in $G$ is of characteristic $p$. 

\endroster 
\endproclaim 

\demo {Proof} Let $K\norm G$ be a normal subgroup of $G$, and set $R=O_p(K)$. Then 
$[O_p(G),K]\leq R$, and so $[O_p(G),C_K(R),C_K(R)]=1$. Then $O^p(C_K(R))\leq C_G(O_p(G))$ by 2.5, 
and thus $O^p(C_K(R))\leq O_p(G)$. Then $C_K(R)$ is a normal $p$-subgroup of $K$, and so 
$C_K(R)\leq R$. This establishes point (a). 

Next, let $U$ be a $p$-subgroup of $G$, and set $H=N_G(U)$, $P=O_p(H)$, and $Q=O_p(G)$. Then 
$N_Q(U)\leq P$. Let $A$ be a $p'$-subgroup of $C_H(P)$. Then $[N_Q(U),A]=1$, and so $[Q,A]=1$ by 2.6. 
Then $A\leq Q$, and thus $A=1$, proving (b). 

For the proof of (c), notice that $C_G(V)\norm G$, and that $X$ is the intersection of $C_G(V)$ with 
the preimage in $G$ of the normal subgroup $C_{G/V}(O_p(G)/V)$ of $G/V$. Thus $X\norm G$. 
Each $p'$-subgroup of $X$ centralizes $O_p(G)$ by 2.5, so $X$ is a $p$-group. 

Finally, let $N\norm G$ with $C_G(N)\leq N$, and set $Q=O_p(N)$.  
Then $[N,C_G(Q)]\leq C_N(Q)\leq Q$, and so $[N,C_G(Q),C_G(Q)]=1$. 
Set $Z=[C_G(Q),C_G(Q)]$. The Three Subgroups Lemma then yields $Z\leq C_G(N)$. As 
$C_G(N)\leq N$ by hypothesis, we conclude that $Z$ is a $p$-group and that $Z\leq Z(C_G(Q))$. Then 
$O^p(C_G(Q))$ is a $p'$-group. As $C_G(Q)\norm G$ and $G$ is of characteristic $p$, we conclude that 
$O^p(C_G(Q))=1$. Thus $C_G(Q)\leq O_p(G)$, as desired. Now let $H$ be a subgroup of $G$ containing $N$. 
Then $C_H(O_p(H))\leq C_H(Q)\leq O_p(G)\cap H\leq O_p(H)$, and thus $H$ is of characteristic $p$. 
This completes the proof of (d). 
\qed 
\enddemo 

The next result refers to the terminology and notation of 1.8. 

\proclaim {Lemma 2.8} Let $(\ca L,\D,S)$ be a proper locality on $\ca F$ and let $P\in\D$. Then 
$P$ is subcentric in $\ca F$, and the following hold.  
\roster 

\item "{(a)}" $P\in\ca F^{cr}$ if and only if $P=O_p(N_{\ca L}(P))$. 

\item "{(b)}" $P$ is centric in $\ca F$ if and only if $C_{\ca L}(P)=Z(P)$. 

\item "{(c)}" $P$ is quasicentric in $\ca F$ if and only if $C_{\ca L}(P)\leq O_p(N_{\ca L}(P))$. 

\endroster 
\endproclaim 

\demo {Proof} Set $M=N_{\ca L}(P)$, and let $Q\in P^{\ca F}$ with $Q$ fully normalized in $\ca F$. Then 
$Q=P^g$ for some $g\in\ca L$  by 2.1, and then the conjugation map $c_g:M\to M^g$ is an isomorphism by 
I.2.3(b). As $\ca F^{cr}$, $\ca F^c$, $\ca F^q$, and $\ca F^s$ are $\ca F$-invariant by 1.9, it follows that 
it suffices to establish the lemma under the assumption (which we now make) that $P$ itself is fully 
normalized in $\ca F$. Since $M$ may be regarded as a proper locality whose set of objects is the set 
of all subgroups of $N_S(P)$, it follows from 2.2 and 2.3 that $O_p(N_{\ca F}(P))=O_p(M)$. As 
$\ca F$ is $\D$-inductive by 2.3, and $M$ is of characteristic $p$, it follows from 1.10 and 1.13 that 
$O_p(M)$ is centric in $\ca F$. Thus $P\in\ca F^s$. 

Set $K=C_{\ca L}(P)$. As $\ca L$ is proper, $M$ is of characteristic $p$, and then $K$ is of characteristic 
$p$ by 2.7(a). Further, we have $C_{\ca F}(P)=\ca F_{C_S(P)}(K)$ by 2.2. If $P\in\ca F^c$ then 
$C_S(P)=Z(P)\leq Z(K)$, so $O^p(K)$ is a normal $p'$-subgroup of $K$. Then $O^p(K)=1$, and $K=Z(P)$. 
Conversely, if $K=Z(P)$ then $C_S(P)\leq P$. As $P$ is fully normalized in $\ca F$, $P$ is also fully 
centralized by 1.13, and we then conclude from 1.10 that $P\in\ca F^c$. This establishes (b). 

Suppose next that $P\in\ca F^q$, so that $C_{\ca F}(P)$ is the trivial fusion system on $C_S(P)$. Then 
$N_K(U)/C_K(U)$ is a $p$-group for every subgroup $U$ of $C_S(P)$. Take $U=O_p(K)$. Thus  
$K/Z(U)$ is a $p$-group, so $K$ is a $p$-group, and then $K=C_S(P)$ is a normal 
$p$-subgroup of $M$. Conversely, if $K\leq O_p(M)$ then $C_{\ca F}(P)$ is trivial, so we have (c). 

Set $R=O_p(M)$ and let $\G$ be the set of all overgroups of $R$ in $N_S(P)$. As $M$ is of characteristic $p$, 
$R$ is centric in $N_{\ca F}(P)$, and then $\G\sub N_{\ca F}(P)^c$. We may view $M$ as a locality 
$(M,\G,N_S(P))$ which happens to be a group, and this locality is proper by 2.7(b). Then 
$R=O_p(N_{\ca F}(P))$ by 2.3. This completes the proof of (a). 
\qed 
\enddemo

The next result shows that if $(\ca L,\D,S)$ is a locality on $\ca F$, such that the set $\D$ of objects is 
not ``too small" ($\ca F^{cr}\sub\D$) and not ``too large" ($\D\sub\ca F^q$), then $\ca L$ has a canonical 
homomorphic image which is a proper locality on $\ca F$.

\proclaim {Proposition 2.9} Let $(\ca L,\D,S)$ be a locality on $\ca F$, with $\ca F^{cr}\sub\D\sub\ca F^q$. 
Let $\Theta$ be the union of the groups $\Theta(P):=O_{p'}(C_{\ca L}(P))$ over all $P\in\D$. Then $\Theta$ 
is a partial normal subgroup of $\ca L$, and $\ca L/\Theta$ is a proper locality on $\ca F$.   
\endproclaim 

\demo {Proof} Let $P\in\D$ and let $Q\in P^{\ca F}$ be fully normalized in $\ca F$. Set $M=N_{\ca L}(Q)$ 
and $K=C_{\ca L}(Q)$. Then $C_{\ca F}(Q)=\ca F_{C_S(Q)}(K)$ by 2.2. As $Q$ is $\ca F$-quasicentric by 
hypothesis, $C_{\ca F}(Q)^{cr}$ is the trivial fusion system on $C_S(Q)$, and a classical theorem of 
Frobenius [Theorem 7.4.5 in Gor] then implies that $K$ has a normal $p$-complement. That is, 
we have $K=\Theta(Q)C_S(Q)$. 

Set $\bar M=M/\Theta(Q)$, set $\bar R=O_p(\bar M)$, and let $R$ be the preimage of $\bar R$ in $N_S(Q)$. 
Let $\bar X$ be a $p'$-subgroup of $C_{\bar M}(\bar R)$, and let $X$ be the preimage of $\bar X$ in $M$. 
Since $M=N_M(R)\Theta(Q)$ by the Frattini argument, we obtain $X=N_X(R)\Theta(Q)$. But 
$$ 
[R,N_X(R)]\leq R\cap[R,X]\leq R\cap\Theta(Q)=1, 
$$ 
so 
$$ 
N_X(R)\leq C_X(R)\leq O^p(C_K(Q))=\Theta(Q), 
$$ 
and thus $X=\Theta(Q)$. This shows that $C_{\bar M}(\bar R)$ is a $p$-group, and thus $\bar M$ is of 
characteristic $p$. In this way the hypothesis (*) of I.4.12 is fulfilled, and we conclude that 
$\Theta\norm\ca L$, $(\ca L/\Theta,\D,S)$ is a locality, $\ca F=\ca F_S(\ca L/\Theta)$, and 
$N_{\ca L/\Theta}(P)$ is of characteristic $p$ for all $P\in\D$. As $\ca F^{cr}\sub\D$, $\ca L/\Theta$ is 
a proper locality on $\ca F$. 
\qed 
\enddemo

There are proper localities $(\ca L,\D,S)$, with fusion system $\ca F$, such that $\D$ is strictly larger 
than $\ca F^q$. For example, if $G$ is a finite group of Lie type, defined over a field of characteristic 
$p$, and $\D$ is the set of all non-identity subgroups of a Sylow $p$-subgroup $S$ of $G$, then a 
theorem of Borel and Tits shows that $N_G(P)$ is of characteristic $p$ for all $P\in\D$. Then  
$(G\mid_\D,\D,S)$ is a proper locality on $\ca F_S(G)$, whereas in general there are non-identity 
subgroups of $S$ which are not quasicentric in $\ca F$.

\proclaim {Lemma 2.10} Let $(\ca L,\D,S)$ be a locality on $\ca F$, and suppose that $\ca F^{cr}\sub\D$. 
Then $\ca F$ is $(cr)$-generated. 
\endproclaim 

\demo {Proof} Let $\G$ be the set of all $R\in\ca F^{cr}$ and such that $R$ is fully normalized in $\ca F$, 
and let $\ca F_0$ be the fusion system on $S$ generated by the union of the groups $Aut_{\ca F}(R)$, for 
$R\in\G$. Assuming the lemma to be false, there  exists an $\ca F$-isomorphism $\phi:P\to P'$ such that 
$\phi$ is not an $\ca F_0$-isomorphism. 

By definition, $\ca F$ is 
generated by the conjugation maps $c_g:S_g\to S_{g\i}$, so we may take $\phi$ to be such a $c_g$, with 
$P=S_g$ and $P'=P^g$ (and where $P$ and $P'$ are in $\D$). Among all such obstructions $g$ to the lemma, 
choose $g$ so that $|P|$ is as 
large as possible. As $S\in\G$, we have $P\neq S$, so $P$ is a proper subgroup of $N_S(P)$. 

Let $Q\in P^{\ca F}$ with $Q$ fully normalized in $\ca F$. As $\ca F$ is $\D$-inductive by 2.2, there 
exist elements $x,y\in\ca L$ such that $N_S(P)^x\leq N_S(Q)\geq N_S(P')^y$, and such that 
$P^x=Q=(P')^y$. As $P$ is a proper subgroup of $N_S(P)$, $c_x$ and $c_y$ are $\ca F_0$-isomorphisms. 
If also $c_x\i\circ c_g\circ c_y$ is an $\ca F_0$-isomorphism, then so is $c_g$, which is contrary to the 
case. We may therefore replace $g$ with $x\i gy$, so that $g\in N_{\ca L}(Q)$. Points (a) and (b) of 2.2 
then imply that $N_{\ca F}(Q)$ is $(cr)$-generated and that $N_{\ca F}(Q)$ is the fusion system of 
$N_{\ca L}(Q)$ over $N_S(Q)$. In this way the maximality of $Q$ reduces the problem to the case where 
$Q\in N_{\ca F}(Q)^{cr}$. Then $Q=O_p(N_{\ca F}(Q))$ and $Q\in\ca F^c$, and so $Q\in\ca F^{cr}$. Thus 
$\phi$ is in $\ca F_0$, and this contradiction proves the lemma.   
\qed 
\enddemo

\proclaim {Lemma 2.11} Let $(\ca L,\D,S)$ be a locality on $\ca F$ and let $\ca N\norm\ca L$ be a partial 
subgroup. Set $T=S\cap\ca N$ and set $\ca H=N_{\ca L}(T)$. Assume that $x^h$ is defined for all $x\in\ca N$ 
and all $h\in\ca H$. Let $\phi:X\to Y$ be an $\ca F$-isomorphism. Then there exists $u\in\bold W(\ca N)$ 
with $X\leq S_v$, and there exists $v\in\bold W(\ca H)$ with $X^u\leq S_v$, such that: 
\roster 

\item "{(a)}" For each entry $h$ of $v$ and each $x\in\ca N$ we have $S_{(h,x)}=S_{(x^{h\i},h)}$, and 

\item "{(b)}" $\phi=c_u\circ c_v$. 

\endroster
\endproclaim 

\demo {Proof} Among all counterexamples choose $\phi:X\to Y$ so that the length $\ell(w)$ of a shortest 
word $w\in\bold W(\ca L)$ with $\phi=c_w$ is as small as possible. Then $\phi$ is not an identity map, 
and so $w$ is not the empty word. Write $w=w_0\circ(g)$ with $g\in\ca L$, and let 
$\phi_0:X\to X^{w_0}$ be the conjugation map induced by $w_0$. The minimality of $\ell(w)$ 
implies that there exist $u_0\in\bold W(\ca N)$ and $v_0\in\bold W(\ca H)$ such that (a) and (b) hold 
with respect to $\phi_0$, $u_0$, and $v_0$. The splitting lemma (I.3.12) 
implies that $g=xh$ with $x\in\ca N$ and $h\in\ca H$, and such that (a) holds. Write 
$v_0=(h_1,\cdots,h_n)$, set $x_0=x$, and recursively define $x_{k+1}=(x_k)^{h_{n-k}\i}$ for $0\leq k<n$.  
Set $y=x_n$, and set $u=u_0\circ(y)$ and $v=v_0\circ(h)$. Then $u$ and $v$ satisfy the conclusion of 
the lemma, and hence $\phi$ is not a counterexample. 
\qed 
\enddemo

\proclaim {Lemma 2.12} Let $(\ca L,\D,S)$ and $(\ca L',\D',S')$ be localities, with fusion systems 
$\ca F=\ca F_S(\ca L)$ and $\ca F'=\ca F_{S'}(\ca L')$. Let $\l:\ca L\to\ca L'$ be a homomorphism of 
partial groups such that $S\l\leq S'$. Then the restriction $\l\mid_S$ of $\l$ to $S$ is fusion-preserving 
with respect to $\ca F$ and $\ca F'$. That is, $\l\mid_S$ is a morphism $\ca F\to\ca F'$ in the category 
of fusion systems. Moreover, if $\l$ is a projection (I.4.4) then $\l\mid_S$ maps $Ob(\ca F)$ onto 
$Ob(\ca F')$, and maps $Mor(\ca F)$ onto $Mor(\ca F')$. 
\endproclaim 

\demo {Proof} For any subgroup $X\leq S$ and any $g\in\ca L$ write $\bar X$ for $X\l$ and $\bar g$ for 
$g\l$. Also, for $X\leq S$ write $\l_X$ for $\l\mid X$. Thus, $\l_S$ is fusion-preserving if for each 
$\ca F$-homomorphism $\a:X\to Y$ there exists an $\ca F'$-homomorphism $\b:\bar X\to\bar Y$ such that 
$$ 
\l_X\circ\b=\a\circ\l_Y. \tag*
$$ 
The uniqueness of $\b$ then follows from the observation that $\bar x\b=(x\a)\l_Y$ for all $x\in X$. 
In the case that $\a=c_g$ for some $g\in\ca L$ the formula (*) is immediate from the definition of 
homomorphism of partial groups, and $\b=c_{\bar g}$ in that case. Since $\ca F$ is generated by such 
homomorphisms $c_g$, (*) then follows from the further observation that if also $h\in\ca L$ and 
$Y^h\leq Z\leq S$, then 
$$ 
c_g\circ c_h\circ\l=\l\mid_X\circ c_{\bar g}\circ c_{\bar h}.  
$$ 
Thus $\l\mid_S$ is fusion-preserving.  

Suppose that $\l$ is a projection. Then $\l$ maps $S$ onto $S'$ by I.4.3(c), and so 
$\l\mid_S$ maps $Ob(\ca F)$ onto $Ob(\ca F')$. Let $\bar w\in\bold W(\ca L')$. By I.3.13(e) there exists 
$w\in\bold W(\ca L)$ such that that the $\l\mid_S$-preimage of $S'_{\bar w}$ in $S$ is a subgroup of 
$S_w$. The mapping $\l^*:\bold W(\ca L)\to\bold W(\ca L')$ of free monoids, induced by $\l$, sends 
$w$ to $\bar w$, and then $\l\mid_S$ sends $c_w$ to $c_{\bar w}$. Thus, $\l\mid_S$ maps 
$Mor(\ca F)$ onto $Mor(\ca F')$. 
\qed 
\enddemo 

\proclaim {Lemma 2.13} Let $(\ca L,\D,S)$ be a locality on $\ca F$, and let $D\leq S$ be a subgroup of 
$S$ such that every $\ca F$-isomorphism $X\to Y$ extends to an $\ca F$-isomorphism $XD\to YD$. Then 
$D\leq O_p(\ca F)$. 
\endproclaim 

\demo {Proof} The set of $\ca F$-isomorphisms is the set of all conjugation maps $c_w:X\to X^w$ with 
$w\in\bold W(\ca L)$ and with $X\leq S_w$. Set $E=\bigcap{S_w\mid w\in\bold W(\ca L)}$. Then 
$E^v=E$ for all $v\in\bold W(\ca L)$ since $\bold W(\ca L)$ is closed under concatenation of words. 
Thus $D\leq E\leq O_p(\ca L)$. 
\qed 
\enddemo

\vskip .2in 
\noindent 
{\bf  Section 3: Restrictions and expansions} 
\vskip .1in 

This section is based closely on [section 5 in Ch1]. 

\definition {Definition 3.1} Let $(\ca L,\D,S)$ be a locality and let $\ca H$ be a partial subgroup of 
$\ca L$. Set $R=S\cap\ca H$ and let $\G$ be a non-empty collection of subgroups of $R$, closed with 
respect to overgroups in $R$. Set 
$$ 
\ca H\mid_\G=\{g\in\ca H\mid S_g\cap R\in\G\}.  
$$ 
Then $(\ca H,\G)$ is a {\it localizable pair} in $\ca L$ if the following conditions hold.  
\roster 

\item "{(1)}" $w\in\bold W(\ca H)$ with $S_w\cap R\in\G\implies w\in\bold D(\ca L)$. In particular, 
$N_{\ca H}(P)$ is a subgroup of $\ca H$ for all $P\in\G$. 

\item "{(2)}" $X^g\in\G$ for all $(X,g)\in\G\times\ca H$ such that $X\leq S_g$. 

\item "{(3)}" $R$ is a maximal $p$-subgroup of $\ca H$. 

\endroster 
Define $\ca F_\G(\ca H)$ to be the fusion system $\ca E$ on $R$ which is generated by the conjugation 
maps $c_g:X\to Y$ where $g\in\ca H$ and where $X,Y\in\G$. If $\ca E^{cr}\sub\G$ and each of the 
groups $N_{\ca H}(P)$ for $P\in\G$ is of characteristic $p$ then 
we say that $(\ca H,\G)$ is {\it properly localizable}. 
\enddefinition

\proclaim {Lemma 3.2} Let $(\ca H,\G)$ be a localizable pair in the locality $\ca L$, and set  
$\ca E:=\ca F_\G(\ca H)$. Then $(\ca H\mid_\G,\G,S\cap\ca H)$ is a locality on $\ca E$, with 
$$
\bold D(\ca H\mid_\G)=\{w\in\bold W(\ca H)\mid S_w\cap R\in\G\},  
$$ 
and the inclusion map $\ca H\mid_\G\to\ca L$ is a homomorphism of partial groups. Moreover, if 
$(\ca H,\G)$ is properly localizable then $(\ca H\mid_\G,\G,S\cap\ca H)$ is a proper locality. 
\endproclaim 

\demo {Proof} By 3.1(1) and (2) $\ca H\mid_\G$ is an objective partial group. As $\G$ is a set of 
subgroups of $R$, 3.1(3) then implies that $(\ca H\mid_\G,\G,S\cap\ca H)$ is a locality, and plainly 
a locality on $\ca E$. Also, 3.1(1) implies that the inclusion $\ca H\mid_\G\to\ca L$ is a homomorphism of 
partial groups. It is immediate from definition 2.4 that $\ca H\mid_\G$ is proper if 
$(\ca H,\G)$ is properly localizable. 
\qed 
\enddemo 

\definition {Example 3.3} Take $\ca H=\ca L$ and let $\G$ be an $\ca F$-closed set of subgroups of $S$. 
The conditions (1) through (3) in definition 3.1 are fulfilled in an obvious way, so $(\ca L,\G)$ is a 
localizable pair. The locality $(\ca L\mid_\G,\G,S)$ will be called the {\it restriction} of $\ca L$ 
to $\G$. 
\enddefinition 

\definition {Example 3.4} Let $R\leq S$ be fully normalized in $\ca F$ and assume that $N_S(R)\in\D$.  
Set $\D_R=\{P\in\D\mid R\norm P\}$. By I.2.12(a) $N_{\ca L}(R)$ is a partial subgroup of $\ca L$.
Then $(N_{\ca L}(R),\D_R)$ is a localizable pair (conditions (1) and (2) of definition 3.1 are 
obvious, while condition (3) is given by I.2.11(b)). Set $\ca F_R=\ca F_\G(N_{\ca L}(R)$. 
We write $\ca L_R$ for the locality $N_{\ca L}(R)\mid_{\D_R}$. 
\enddefinition

\definition {Definition 3.5} Let $(\ca L,\D,S)$ be a locality. An {\it expansion} of $\ca L$ consists of 
a locality $(\ca L',\D',S')$ together with an injective homomorphism $\eta:\ca L\to\ca L'$ such that 
$\eta$ restricts to an isomorphism $S\to S'$, and such that $Im(\eta)$ is a restriction of $\ca L'$. 
Such an expansion is {\it elementary} if there exists $R\in\D'$ such that $\D'=\D\eta\bigcup R^{\ca F}$. 
\enddefinition 

Of course, one always has trivial expansions of $\ca L$, for which $\eta$ is an isomorphism. 
Our goal in the remainder of this section is to construct a non-trivial elementary expansion of 
$\ca L$, under conditions (see Hypothesis 3.7 below) which will be found to hold if $\ca L$ 
is proper, and provided that $\D$ is not already ``as large as possible". The following lemma prepares 
the way for such a construction.

\proclaim {Lemma 3.6} Let $(\ca L,\D,S)$ be a locality on $\ca F$, let $R\leq S$ with $R$ fully 
normalized in $\ca F$, and suppose that $\<U,V\>\in\D$ for every pair $U\neq V$ of distinct 
$\ca F$-conjugates of $R$. Set $\D_R=\{P\in\D\mid R\norm P\}$. Then the following hold. 
\roster 

\item "{(a)}" $\D_R=\{N_Q(R)\mid R\leq Q\in\D\}$, and $(N_{\ca L}(R),\D_R,N_S(R))$ is a locality. 
Moreover, if every strict overgroup of $R$ in $S$ is in $\D$, then  
$\ca F_{N_S(R)}(N_{\ca L}(R)=N_{\ca F}(R))$. 

\item "{(b)}" $R^{\ca F}=\{R^g\mid g\in\ca L,\ R\leq S_g\}$. 

\item "{(c)}" $\ca F$ is $\D\cup R^{\ca F}$-inductive. 

\endroster 
\endproclaim 

\demo {Proof} Let $Q\in\D$ with $R\leq Q$. If $R\norm Q$ then $N_Q(R)=Q\in\D$, while otherwise 
$N_Q(R)$ contains a pair of distinct $Q$-conjugates of $R$. Thus $N_Q(R)\in\D$ in either case, and  
thus $\D_R=\{N_Q(R)\mid R\leq Q\in\D\}$. 

For any $g\in N_{\ca L}(R)$ we have $N_{S_g}(R)\in\D$, so 
$(N_{\ca L}(R),\D_R,N_S(R))$ is a locality by 3.4. Set $\ca F_R=\ca F_{N_S(R)}(N_{\ca L}(R))$. 
Then $\ca F_R$ is a fusion subsystem of $N_{\ca F}(R)$, and it remains to show that $\ca F_R=N_{\ca F}(R)$ 
in order to complete the proof of (a). The proof of this last step in point (a) will require (c), which we 
now prove together with (b). 

Let $V\in R^{\ca F}$, and let $\bold Y_V$ be the set of elements $y\in\ca L$ such that $R^y=V$ and such that 
$N_S(V)\leq S_{y\i}$. By definition, each $\ca F$-isomorphism $\phi:V\to R$ can be factored as a composition 
of conjugation maps 
$$ 
V @>c_{x_1}>> V_1@>>> \cdots @>>>V_{k-1} @>c_{x_k}>>R \tag* 
$$
for some word $w=(x_1,\cdots,x_k)\in\bold W(\ca L)$. Write $V^w=R$ to indicate this. Assume now that 
there exists $V\in R^{\ca F}$ such that $\bold Y_V=\nset$. Among all such $V$, 
choose $V$ so that the minimum length $k$, taken over all words 
$w$ for which $V^w=R$, is as small as possible. Then, subject to this condition, choose $w$ so that 
$|N_{S_w}(V)|$ is as large as possible. We evidently have $k>0$, and $R\neq S$. 

Suppose that $k=1$. Thus $w=(x)$ for some $x\in\ca L$ with $V^x=R$. Set $P=N_{S_x}(V)$ and 
set $\w P=N_{N_S(V)}(P)$. Then $P\in\D$, and conjugation by $x$ then induces an isomorphism 
$N_{\ca L}(P)\to N_{\ca L}(P^x)$ by I.2.3(b). Thus $\w P^x$ is a $p$-subgroup of the locality $N_{\ca L}(R)$.  
By I.2.11 there exists $z\in N_{\ca L}(R)$ with 
$(\w P^x)^z\leq N_S(R)$, and then $(x,z)\in\bold D$ via $P$. Now $\w P^{xz}\leq N_S(R)$, and the 
maximality of $|N_{S_w}(V)|$ in the choice of $V$ and $w$ yields $P=\w P$. Then $P=N_S(V)$, so 
$x\i\in\bY_V$. This shows that $k>1$. 

The minimality of $k$ now yields $\bY_{V_1}\neq\nset$, where $V_1$ is defined by (*). Thus $V^{(c,d)}=R$, 
where $d$ may be chosen in $\bY_{V_1}$. Set $A=N_{S_c}(V)$. Then $A\in\D$ as $S_c\in\D$,  
and $A^{c}\leq N_S(V_1)$. As $b\in\bY_{V_1}$ we then have $w\in\bold D$ via $A$. Thus $k=1$, contrary to 
the result of the preceding paragraph. We conclude that $\bold Y_V=\nset$ for $V\in R^{\ca F}$, and 
this proves (b) and (c). 

We now return to the final step in the proof of (a). Set $\ca F_R=\ca F_{N_S(R)}(N_{\ca L}(R))$, and  
assume that $\D_R$ contains the set of all strict overgroups of $R$ in $S$. Assume also, by way of 
contradiction, that $\ca F_R$ is a proper subsystem of $N_{\ca F}(R)$. Then there exists a 
subgroup $Y$ of $N_S(R)$, containing $R$ and fully normalized in $N_{\ca F}(R)$, such that 
$Aut_{\ca F_R}(Y)$ is a proper subgroup of $Aut_{N_{\ca F}(R)}(Y)$. If $R$ is a proper subgroup of 
$Y$ then $Y\in\D$, and then every $N_{\ca F}(R)$-automorphism of $Y$ is induced by conjugation 
by an element of $N_{\ca L}(Y)$ which fixes $R$. That is, we have 
$Aut_{\ca F_R}(Y)=Aut_{N_{\ca F}(R)}(Y)$ as is contrary to the case. Thus $Y=R$. Now let 
$\a\in Aut_{\ca F}(R)$ and let $w\in\bold W(\ca L)$ with $c_w=\a$. Write $w=(g_1,\cdots, g_n)$, set 
$R_0=R$, and recursively define $R_k$ for $0<k\leq n$ by $R_k=R^{g_k}$. Thus $R_n=R$. 
Set $a_0=a_n=\1$. As a consequence of (b) we may choose and for each $k$ with $0<k<n$ an element 
$a_k\in\ca L$ such that $(R_k)^{a_k}=R$ and such that $N_S(R_k)\leq S_{a_k}$. Then the word 
$u_k=(a_{k-1},g_k,a_k)$ (with $0<k\leq n$) is in $\D$ via $N_{S_{a_{k-1}}}(R_{k-1})^{a_{k-1}}$. 
Setting $f_k=\Pi(u_k)$ and $w'=(f_1,\cdots,f_n)$, we obtain $f_k\in N_{\ca L}(R)$ and $\a=c_{w'}$. 
Thus $\a$ is an $\ca F_R$-automorphism of $R$, which is again contrary to the case. We conclude that 
$\ca F_R=N_{\ca F}(R)$, completing the proof of (a) and of the lemma.  
\qed 
\enddemo

We assume the following setup for the remainder of this section.

\definition {Hypothesis 3.7} $(\ca L,\D,S)$ is a locality on $\ca F$, and $R\leq S$ is 
a subgroup of $S$ such that: 
\roster 

\item "{(1)}" Each strict overgroup of $R$ in $S$ is in $\D$. 

\item "{(2)}" Both $R$ and $O_p(N_{\ca F}(R))$ are fully normalized in $\ca F$. 

\item "{(3)}" $N_{\ca L}(R)$ is a subgroup of $\ca L$, and $N_{\ca F}(R)=\ca F_{N_S(R)}(N_{\ca L}(R))$.  

\endroster 
\enddefinition 

For each $V\in R^{\ca F}$ let $\bY_V$ be the set of elements $y\in\ca L$ such that $R^y=V$ and such that 
$N_S(V)\leq S_{y\i}$. Define $\bY$ to be the union of the sets $\bY_V$, taken over all $V\in R^{\ca F}$. 
Notice that points (b) and (c) of 3.6  are equivalent to the condition that $\bY_V$ be non-empty   
for each $V\in R^\ca F$.

\proclaim {Theorem 3.8} Assume Hypothesis 3.7, and set $\D^+=\D\bigcup R^{\ca F}$. Then there exists a 
locality $(\ca L^+,\D^+,S)$ such that $\ca L$ is the restriction of $\ca L^+$ to $\D$, 
and such that $N_{\ca L^+}(R)=N_{\ca L}(R)$. Moreover, we then have $\ca F_S(\ca L^+)=\ca F$, and 
the following hold. 
\roster 

\item "{(a)}" For any locality $(\w{\ca L},\D^+,S)$ such that $\w{\ca L}\mid_\D=\ca L$, and 
such that $N_{\w{\ca L}}(R)=N_{\ca L}(R)$, we have $\ca F_S(\w{\ca L})=\ca F$, and there is a unique 
isomorphism $\b:\ca L^+\to\w{\ca L}$ which restricts to the identity map on $\ca L$. 

\item "{(b)}" Let $\ca L_0^+$ be the set of all $g\in\ca L^+$ such that $S_g$ contains an $\ca F$-conjugate 
of $R$. Then $\ca L_0^+$ is a partial subgroup of $\ca L$, and $\ca L^+$ is the pushout in the 
category of partial groups of the diagram 
$$ 
\ca L_0^+@<<<\ca L\cap\ca L_0@>>>\ca L 
$$ 
of inclusion homomorphisms. 

\endroster 
\endproclaim

Along the way to proving Theorem 3.8, we shall explicitly determine the partial group structure of $\ca L^+$  
in 3.14 through 3.18. These results will then play an important computational role in section 4. 
\vskip .1in 
Set $\bX=\{x\i\mid x\in\bY\}$ and set 
$$ 
\Phi=\bX\times N_{\ca L}(R)\times\bY. 
$$ 
For any $\phi=(x\i,h,y)\in\Phi$ set $U_\phi=R^x$ and $V_\phi=R^y$. Thus: 
$$ 
U_\phi@>x\i>>R@>h>>R@>y>>V_\phi, \quad\text{and}\quad N_S(U_\phi)@>x\i>>N_S(R)@<y\i<<N_S(V_\phi)
$$ 
are diagrams of conjugation maps, labelled by the conjugating elements. (In the first of these diagrams 
the conjugation maps are isomorphisms, and in the second they are homomorphisms into $N_S(R)$.)

\definition {(3.9)} Define a relation $\sim$ on $\Phi$ as follows. For $\phi=(x\i,h,y)$ and 
$\bar\phi=(\bar x\i,\bar h,\bar y)$ in $\Phi$, write $\phi\sim\bar\phi$ if 
\roster 

\item "{(i)}" $U_\phi=U_{\bar\phi}$, $V_\phi=V_{\bar\phi}$, and 

\item "{(ii)}" $(\bar xx\i)h=\bar h(\bar y y\i)$. 

\endroster 
\enddefinition 

The products in 3.9(ii) are well-defined. Namely, by 3.9(i), $(\bar x,x\i)\in\bold D$ via 
$N_S(U)^{\bar x\i}$, $\bar xx\i\in N_{\ca L}(R)$, and then $(\bar xx\i,h)\in\bold D$ since 
$N_{\ca L}(R)$ is a group (3.7(3)). The same considerations apply to $(\bar y,y)$ and 
$(\bar h,\bar y y\i)$. 

One may depict the relation $\sim$ by means of a commutative diagram, as follows. 
$$ 
\CD 
U   @>x\i>> R@>g>>R@>y>>V \\ 
@|  @A \bar x x\i AA  @AA\bar y y\i A  @| \\ 
U   @>>\bar x\i> R@>>\bar g> R@>>\bar y> V 
\endCD 
$$

\proclaim {Lemma 3.10} $\sim$ is an equivalence relation on $\Phi$.
\endproclaim

\demo {Proof} Evidently $\sim$ is reflexive and symmetric. Let $\phi_i=(x_i\i,g_i,y_i)\in\Phi$ 
($1\leq i\leq 3$) with $\phi_1\sim\phi_2\sim\phi_3$. Then $R^{x_1}=R^{x_3}$ and 
$R^{y_1}=R^{y_3}$. Notice that 
$$ 
\align 
(x_3,x_2\i,x_2,x_1\i)&\in\bold D\quad\text{via $N_S(U)^{x_3\i}$ and,} \\ 
(y_3,y_2\i,y_2,y_1\i)&\in\bold D\quad\text{via $N_S(V)^{y_3\i}$}. 
\endalign 
$$ 
Computation in the group $N_{\ca L}(R)$ then yields 
$$ 
\align 
(x_3 x_1\i)g_1&=(x_3 x_2\i)(x_2 x_1\i)g_1=(x_3 x_2\i)g_2(y_2 y_1\i) \\ 
&=g_3(y_3 y_2\i)(y_2 y_1\i)=g_3(y_3 y_1\i), 
\endalign 
$$
which completes the proof of transitivity. 
\qed 
\enddemo

\proclaim {Lemma 3.11} Let $\psi\in\Phi$ and set $U=U_\psi$ and $V=V_\phi$. Let $x\in\bY_U$ and 
$y\in\bY_V$. Then there exists a unique $h\in N_{\ca L}(V)$ such that $\psi\sim(x\i,h,y)$. 
\endproclaim 

\demo {Proof} Write $\psi=(\bar x\i,\bar h,\bar y)$. Then $(\bar x,x\i)\in\bold D$ via $N_S(U)^{\bar x\i}$, 
and $\bar xx\i\in N_{\ca L}(R)$. Similarly $(\bar y,y\i)\in\bold D$ and $\bar yy\i\in N_{\ca L}(R)$. 
As $N_{\ca L}(R)$ is a subgroup of $\ca L$ we may form the product $h:=(x\bar x\i)\bar h(\bar yy\i)$ 
and obtain $(\bar xx\i)h=\bar h(\bar yy\i)$. Setting $\phi=(x\i,h,y)$, we thus have $\phi\sim\psi$. 
If $h'\in N_{\ca L}(R)$ with also $(x\i,h',y)\sim\psi$ then $h'=(x\bar x\i)\bar h(\bar yy\i)=h$. 
\qed 
\enddemo 

For ease of reference we record the following observation, even though it is simply part of 
the definition of the relation $\sim$.

\proclaim {Lemma 3.12} Let $C$ be a $\sim$-class of $\Phi$, let $\phi=(x\i,g,y)\in C$, and 
set $U=R^x$ and $V=R^y$. Then the pair $(U,V)$ depends only on $C$, and not on the choice of  
representative $\phi$. 
\qed 
\endproclaim

\proclaim {Lemma 3.13} Let $\ca L_0$ be the set of all $g\in\ca L$ such that $S_g$ contains an 
$\ca F$-conjugate $U$ of $R$. For each $g\in\ca L_0$ set 
$$ 
\Phi_g=\{\phi\in\Phi\cap\bold D\mid \Pi(\phi)=g\},   
$$ 
and set 
$$ 
\ca U_g=\{U\in R^{\ca F}\mid U\leq S_g\}. 
$$ 
\roster 

\item "{(a)}" $\ca U_g$ is the set of all $U_\phi$ such that $\phi\in\Phi_g$. 

\item "{(b)}" $\Phi_g$ is a union of $\sim$-classes.  

\item "{(c)}" Let $U\in\ca U_g$, set $V=U^g$, and let $x\in\bX_U$ and $y\in\bX_V$. Then 
$(x,g,y\i)\in\bold D$, $h:=xgy\i\in N_{\ca L}(U)$, and $\phi:=(x\i,h,y)\in\Phi_g$. If also  
$(x\i,h',y)\in\Phi_g$ then $h=h'$. 

\item "{(d)}" Let $\phi$ and $\psi$ in $\Phi_g$. Then  
$$ 
\phi\sim\psi\iff U_\phi=U_\psi\iff V_\phi=V_\psi. 
$$ 

\endroster 
\endproclaim 

\demo {Proof} Let $g\in\ca L_0$, let $U\in\ca U_g$, and set $V=U^g$. Let $(x\i,y)\in\bX_U\times\bY_V$ and 
set $w=(x,g,y\i)$. Then $w\in\bold D$ via 
$$
(N_{S_g}(U)^{x\i},N_{S_g}(U),N_{S_{g\i}}(V),N_{S_{g\i}}(V)^{y\i}). 
$$ 
Set $h=\Pi(w)$. Then $h\in N_{\ca L}(R)$ since $R^x=U$ and $R^y=V$. Set $w'=(x\i,x,g,y\i,y)$. Then 
$w'\in\bold D$ via $N_{S_w}(U)$, and $g=\Pi(w')=\Pi(x\i,h,y)$ by $\bold D$-associativity (I.1.4(b)). If 
also $(x\i,h',y)\in\bold D$ with $\Pi(x\i,h',y)=g$ then $h=h'$ by the cancellation rule (I.1.4(e)). This 
establishes (c), and shows that $\ca U_g\sub\{U_\phi\mid\phi\in\Phi_g\}$. The opposite inclusion is 
immediate (cf. I.2.3(c)), so also (a) is established. 

Let $\phi\in\Phi_g$ and let $\bar\phi\in\Phi$ with $\phi\sim\bar\phi$. Write $\phi=(x\i,h,y)$ and 
$\bar\phi=(\bar x\i,\bar h,\bar y)$, and set 
$$ 
w'=(\bar x\i,\bar x,x\i,h,y,\bar y\i,\bar y). \tag*
$$ 
As $\phi\in\bold D$ we have $N_{S_\phi}(U)\in\D$ by 3.7(1), and then $w'\in\bold D$ via $N_{S_\phi}(U)$. 
Now $\Pi(w')=\Pi(\phi)=g$, while also 
$$ 
\Pi(w')=\Pi(\bar x\i,\Pi((\bar xx\i),h,(y\bar y\i)),\bar y)=\Pi(\bar x\i,\bar h,\bar y)=\Pi(\bar\phi), 
\tag**
$$ 
and thus $\bar\phi\in\Phi_g$. This proves (b). 

It remains to prove (d). So, let $\phi,\psi\in\Phi_g$. If $\phi\sim\psi$ then $U_\phi=U_\psi$ and 
$V_\phi=V_\psi$ by 3.12. On the other hand, assume that $U_\phi=U_\psi$ or that $V_\phi=V_\psi$. Then 
both equalities obtain, since $V_{\phi}=(U_{\phi})^g$ and $V_{\psi}=(U_{\psi})^g$. Write 
$\phi=(x\i,h,y)$ and $\psi=(\bar x\i,\bar h,\bar y)$ in the usual way, and define $w'$ as in (*). 
Then $w'\in\bold D$ via $N_{S_\phi}(V_\phi)$, and $\Pi(w')=\Pi(\phi)=g$. Set 
$$ 
h'=\Pi((\bar xx\i)h(y\bar y\i)). 
$$ 
Then $g=\Pi(w')=\Pi(\bar x\i,h',\bar y)$, so $(\bar x\i,h',\bar y)\in\Phi_g$. Then 3.11 yields 
$h'=\bar h$, and thus $(\bar x\i,h',\bar y)=\bar\phi$. This shows that $\phi\sim\bar\phi$, completing 
the proof of (d). 
\qed 
\enddemo 

\definition {(3.14)} We now have a partition of the disjoint union $\ca L\bigsqcup\Phi$ (and a 
corresponding equivalence relation $\approx$ on $\ca L\bigsqcup\Phi$) by means of three types of 
$\approx$-classes, as follows.   
\roster 

\item "{$\cdot$}" Singletons $\{f\}$, where $f\in\ca L$ and where $S_f$ contains no 
$\ca F$-conjugate of $R$ (classes whose intersection with $\Phi$ is empty). 

\item "{$\cdot$}" $\sim$-classes $[\phi]$ such that $[\phi]\cap\bold D=\nset$  
(classes whose intersection with $\ca L$ is empty).  

\item "{$\cdot$}" Classes $\Phi_g\cup\{g\}$, where $g\in\ca L$ and where $S_g$ contains an 
$\ca F$-conjugate of $R$ (classes having a non-empty intersection with both $\ca L$ and $\Phi$). 

\endroster 
For any element $E\in\ca L\cup\Phi$, write $[E]$ for the $\approx$-class of $E$. (Thus $[E]$ is also 
a $\sim$-class if and only if $[E]\cap\bold D=\nset$.) Let $\ca L^+$ be the set of all $\approx$-classes. 
Let $\ca L_0^+$ be the subset of $\ca L^+$, consisting of those $\approx$-classes whose  
intersection with $\Phi$ is non-empty. That is, the members of $\ca L_0^+$ are the $\approx$-classes 
of the form $[\phi]$ for some $\phi\in\Phi$. 
\enddefinition

Recall from Part I that there is an inversion map $w\maps w\i$ on $\bold W(\ca L)$, given by 
$(g_1,\cdots,g_n)\i=(g_n\i,\cdots,g_1\i)$. The following result is then a straightforward 
consequence of the definitions of $\Phi$, $\sim$, and $\approx$. 

\proclaim {Lemma 3.15} The inversion map on $\bold W(\ca L)$ preserves $\Phi$. Further,  
for each $E\in\ca L\cup\Phi$ we have $E\i\in\ca L\bigcup\Phi$, and $[E\i]$ is the 
set $[E]\i$ of inverses of members of $[E]$. 
\qed 
\endproclaim

\definition {Definition 3.16} For any $\g=(\phi_1,\cdots,\phi_n)\in\bold W(\Phi)$ let $w_\g$ be the 
word $\phi_1\circ\cdots\circ\phi_n$ in $\bold W(\ca L)$. Let $\G$ be the set of all $\g\in\bold W(\Phi)$ 
such that  $S_{w_\g}$ contains an $\ca F$-conjugate of $R$. Let $\bold D_0^+$ be the set of all sequences 
$w=([\phi_1],\cdots,[\phi_n])\in\bold W(\ca L_0^+)$ for which there exists a sequence 
$\g$ of representatives for $w$ with $\g\in\G$. We shall say that $\g$ is a {\it $\G$-form} of $w$. 
\enddefinition

The following lemma shows how to define a product $\Pi_0^+:\bold D_0^+\to\ca L_0^+$.

\proclaim {Lemma 3.17} Let $w=([\phi_1],\cdots,[\phi_n])\in\bold D_0^+$, and let $\g=(\phi_1,\cdots,\phi_n)$ 
be a $\G$-form of $w$. Write $\phi_i=(x_i\i,h_i,y_i)$. 
\roster 

\item "{(a)}" $(y_i,x_{i+1}\i)\in\bold D$ and $y_ix_{i+1}\i\in N_{\ca L}(R)$ for each $i$ with $1\leq i<n$. 

\item "{(b)}" Set 
$$ 
w_0=(h_1,y_1x_2\i,\cdots,y_{n-1}x_n\i,h_n). 
$$ 
Then $w_0\in\bold W(N_{\ca L}(R))$ and $(x_1\i,\Pi(w_0),y_n)\in\Phi$. Moreover: 

\item "{(c)}" The $\approx$-class $[x_1\i,\Pi(w_0),y_n]$ depends only on $w$, and not on 
the choice of $\G$-form of $w$. 

\endroster 
\endproclaim 

\demo {Proof} Let $U\in R^{\ca F}$ and let $x,y\in\bY_U$. Then $(y,x\i)\in\bold D$ via $N_S(U)^{y\i}$, 
and then $yx\i\in N_{\ca L}(R)$. This proves (a), and shows that $w_0\in\bold W(N_{\ca L}(R))$. As 
$N_{\ca L}(R))$ is a subgroup of $\ca L$, (b) follows. 

Let $\bar\g=(\bar\phi_1,\cdots,\bar\phi_n)$ be any $\G$-form of $w$, write 
$\bar\phi_i=(\bar x_i\i,\bar h_i,\bar y_i)$, and define $\bar w_0$ in analogy with $w_0$. Set 
$U_0=U_{\phi_1}$ and $\bar U_0=U_{\bar\phi_1}$. For each $i$ with $1\leq i\leq n$ set $U_i=V_{\phi_i}$ 
and $\bar V_i=U_{\bar\phi_i}$. Thus:  
$$ 
U_{i-1}@>x_i\i>>R@>h_i>>R@>y_i>>U_i\quad\text{and}\quad 
\bar U_{i-1}@>\bar x_i\i>>R@>\bar h_i>>R@>\bar y_i>>\bar U_i. 
$$ 

Suppose that there exists an index $j$ with $U_j\neq\bar U_j$. As $\phi_j\approx\bar\phi_j$ it follows 
from 3.10 and 3.11 that $\phi_j$ and $\bar\phi_j$ are in $\bold D$, and that there is an element 
$g_j\in\ca L_0$ such that $\phi_j$ and $\bar\phi_j$ are in $\Phi_g$. If $j<n$ then  
$$ 
U_{j+1}=(U_j)^g\neq(\bar U_j)^g=\bar U_{j+1},  
$$ 
and if $j>0$ one obtains $U_{j-1}\neq\bar U_{j-1}$ in similar fashion, by consideration of $w\i$ 
via 3.14. Thus, for each index $i$ 
we have $U_i\neq\bar U_i$. Then $\phi_i$ and $\bar\phi$ lie in distinct $\sim$-classes by 3.10. As  
$\phi_i\approx\bar\phi_i$ it follows that $\phi_i$ and $\bar\phi_i$ are members of $\Phi\cap\bold D$, 
and that there is a word $v=(g_1,\cdots,g_n)\in\bold W(\ca L)$ with $g_i=\Pi(\phi_i)=\Pi(\bar\phi_i)$. 
As $\<U_0,\bar U_0\>\leq S_v$, we have $v\in\bold D$. 

Set $P=S_v$ and set $P_0=N_P(U_0)$. Then $P_0\in\D$ by 3.7(1). Set $P_i=(P_{i-1})^{g_i}$ for $1\leq i\leq n$.
Then 
$$
P_{i-1}^{x_i\i}\leq N_S(R)\ \ \text{and}\ \ (P_i)^{y_i\i}\leq N_S(R).  
$$ 
Thus conjugation by $g_i$ maps $P_{i-1}^{x_i\i}$ to $(P_i)^{y_i\i}$, and this shows: 
\roster 

\item "{(*)}" Set $w_\g=\g_1\circ\cdots\circ\g_n$. Then $w_\g\in\bold D$ via $P_0$. 

\endroster 
Similarly, one obtains $w_{\bar\g}\in\bold D$ via $N_P(\bar U_0)$, where 
Now $\bold D$-associativity yields: 
$$ 
\Pi(x_1\i,\Pi(w_0),y_n)=\Pi(w_\g)=\Pi(v)=\Pi(w_{\bar \g})=\Pi(\bar x_1\i,\Pi(\bar w_0),\bar y_n).  
$$ 
Thus $\Pi(x_1\i,\Pi(w_0),y_n)$ and $(\bar x_1\i,\Pi(\bar w_0),\bar y_n)$ lie in the same fiber of 
$\Pi:\Phi\cap\bold D\to\ca L_0$, and thus 
$(x_1\i,\Pi(w_0),y_n)\approx(\bar x_1\i,\Pi(\bar w_0),\bar y_n)$. This reduces (c) to the following claim. 
\roster 

\item "{(**)}" Suppose that $U_i=\bar U_i$ for all $i$. Then 
$(x_1\i,\Pi(w_0),y_n)\sim(\bar x_1\i,\Pi(\bar w_0),\bar y_n)$. 

\endroster 
Among all counter-examples to (**), let $w$ be chosen so that $n$ is as small as possible. Then 
$n\neq 0$ (i.e. $w$ and $\bar w$ are non-empty words), as there is otherwise 
nothing to verify. If $n=1$ (so that $w=(\phi_1)$ and $\bar w=(\bar\phi_1)$) 
then $w_0=(h)$, $\bar w_0=(\bar h)$, and (**) follows since $\phi\sim\bar\phi$. Thus, $n\geq 2$. 

Set $\psi=(x_1\i,h_1(y_1x_2\i)h_2,y_2)$, and similarly define $\bar\psi$. Then $\psi$ and $\bar\psi$ 
are in $\Phi$. As $\phi_i\sim\bar\phi_i$ for $i=1,2$, one has the commutative diagram 
$$ 
\CD 
U_0@>x_1\i>>R@>h_1>>R@>y_1>>U_1@>x_2\i>>R@>h_2>>R@>y_2>>U_2  \\  
@|  @A{\bar x_1x_1\i}AA @AA{\bar y_1y_1\i}AA @| @A{\bar x_2x_2\i}AA @AA{\bar y_2y_2\i}A @|   \\  
U_0@>\bar x_1\i>>R@>\bar h_1>>R@>\bar y_1>>U_1@>\bar x_2\i>>R@>\bar h_2>>R@>\bar y_2>>U_2
\endCD 
$$ 
of conjugation maps. This diagram then collapses to the commutative diagram 
$$ 
\CD 
U_0 @>x_1\i>> R @>h_1(y_1x_2\i)h_2>> R @>y_2>>U_2  \\  
@|  @A{\bar x_1x_1\i}AA  @AA{\bar y_2y_2\i}A  @|   \\
U_0 @>\bar x_1\i>> R @>\bar h_1(\bar y_1\bar x_2\i)\bar h_2>> R @>\bar y_2>>U_2 
\endCD 
$$
which shows that $\psi\sim\bar\psi$. Here $(\psi,\phi_3,\cdots,\phi_n)$ and 
$(\bar\psi,\bar\phi_3,\cdots,\bar\phi_n)$ are $\G$-forms of the word 
$u=([\psi],[\phi_3],\cdots,[\phi_n])$. Set 
$$ 
u_0=(h_1(y_1x_2\i)h_2,y_2x_3\i,\cdots,y_{n-1}x_n\i,h_n), 
$$ 
and similarly define $\bar u_0$. Then $\Pi(u_0)=\Pi(w_0)$ and $\Pi(\bar u_0)=\Pi(\bar w_0)$. 
The minimality of $n$ then yields 
$$ 
(x\i,\Pi(w_0),y_n)\approx(\bar x_1\i,\Pi(\bar w_0),\bar y_n). 
$$ 
This proves (**), and thereby completes the proof of (c). 
\qed 
\enddemo

\proclaim {Proposition 3.18} There is a mapping $\Pi^+:\bold D_0^+\to\ca L_0^+$, given by 
$$ 
\Pi_0^+(\nset)=[\1,\1,\1],    
$$ 
and by 
$$
\Pi^+_0(w)=[x_1\i,\Pi(w_0),y_n], \tag* 
$$ 
on non-empty words $w\in\bold D_0^+$; where $w_0$ is given by a $\G$-form of $w$ as in 3.15(b).  
Further, there is an involutory bijection on $\ca L_0^+$ given by 
$$ 
[x\i,g,y]\i=[y\i,g\i,x]. \tag** 
$$
With these structures, $\ca L_0^+$ is a partial group. 
\endproclaim 

\demo {Proof} The reader may refer to I.1.1 for the conditions (1) through (4) defining the notion of 
partial group. Condition (1) requires that $\bold D_0^+$ contain all words of length 1 in the 
alphabet $\ca L_0^+$, and that $\bold D_0^+$ be closed with respect to decomposition. (That is, 
if $u$ and $v$ are two words in the free monoid $\bold W(\ca L_0^+)$, and the 
concatenation $u\circ v$ is in $\bold D_0^+$, then $u$ and $v$ are in $\bold D_0^+$.) Both of these   
conditions are immediate consequences of the definition of $\bold D_0^+$. 

That $\Pi_0^+$ is a well-defined mapping is given by 3.15. 
The proof that $\Pi_0^+$ satisfies the conditions I.1.1(2) 
($\Pi_0^+$ restricts to the identity map on words of length 1) and I.1.1(3): 
$$ 
\bu\circ\bv\circ\w\in\bold D_0^+\implies\Pi_0^+(\bu\circ\bv\circ\bv)=\Pi_0^+(\bu\circ\Pi_0^+(\bv)\circ\bw)
$$ 
are then straightforward, and may safely be omitted. 

The inversion map $[x\i,h,y]\maps [y\i,h\i,x]$ is well-defined by 3.13. Evidently this mapping is an 
involutory bijection, and it extends to an involutory bijection 
$$ 
(C_1,\cdots,C_n)\i=(C_n\i,\cdots,C_1\i) 
$$ 
on $\bold W(\ca L_0^+)$. It thus remains to show I.1.1(4). That is, we must check that  
$$ 
w\in\bold D_0^+\implies w\i\circ w\in\bold D_0^+\ \ \text{and}\ \ 
\Pi_0^+(w\i\circ w)=[\1,\1,\1]. 
$$ 
In detail: take $w=(C_1,\cdots,C_n)$ and let $\g=(\phi_1,\cdots,\phi_n)$ be a $\G$-form of $w$, where  
$\phi_i$ is written as $[x_i\i,h_i,y_i)$. 
One easily verifies that $\g\i\circ\g\in\G$, and hence $w\i\circ w\in\bold D_0^+$. Now 
$$ 
\Pi_0+(w\i\circ w)=[y_n\i,\Pi(u_0),y_n], 
$$ 
where 
$$  
u_0=(g_n\i,x_ny_{n-1}\i,\cdots,x_2y_1\i,g_1\i,x_1x_1\i,g_1,y_1x_2\i),\cdots,y_{n-1}x_n\i,g_n). 
$$ 
One observes that $\Pi(u_0)=\1$, and so $\Pi_0+(w\i\circ w)=[y_n\i,\1,y_n]$. Now observe that 
$(y_n\i,\1,y_n)\equiv(\1,\1,\1)$, since $(y_n\i,\1,y_n)$ and $(\1,\1,\1)$ are in $\bold D$ and  
since $\Pi(y_n\i,\1,y_n)=\1=\Pi(\1,\1,\1)$. Thus $\Pi_0+(w\i\circ w)=[\1,\1,\1]=\Pi_0^+(\nset)$. Thus 
I.1.1(4) holds, and the proof is complete. 
\qed 
\enddemo

\proclaim {Lemma 3.19} Let $\bold D_0$ be the set of all $w\in\bold D$ such that $S_w$ contains an 
$\ca F$-conjugate of $R$, and let $\ca L_0$ be the set of words of length 1 in $\bold D_0$, regarded as 
a subset of $\ca L$. Let $\Pi_0:\bold D_0\to\ca L_0$ be the restriction of $\Pi$ to $\bold D_0$. Then: 
\roster 

\item "{(a)}" $\ca L_0$, with $\Pi_0$ and the restriction of the inversion map on $\ca L$ to $\ca L_0$, 
is a partial group. 

\item "{(b)}" Let $\iota_0:\ca L_0\to\ca L$ be the inclusion map, and let $\l_0:\ca L_0\to\ca L_0^+$ 
be the mapping $g\maps\Phi_g\cup\{g\}$. Then $\iota_0$ and $\l_0$ are homomorphisms of partial groups. 

\endroster 
\endproclaim 

\demo {Proof} The verification of (a) is straightforward, and is left to the reader (see I.1.1). 
Moreover, since the product in $\ca L_0$ is inherited from $\ca L$, it is immediate that $\iota_0$ is a 
homomorphism of partial groups. 

Let $v=(g_1,\cdots,g_n)\in\bold D_0$, set $w=([g_1],\cdots,[g_n])$, and let $U\in R^{\ca F}$ be 
chosen so that $U\leq S_v$. By 3.17 there exists a word $\g=(\phi_1,\cdots,\phi_n)\in\bold W(\Phi)$ 
such that $\phi_i\in[g_i]$ and such that $U\leq S_{w_\g}$, where $w_\g$ is the word 
$\phi_1\circ\cdots\circ\phi_n\in\bold W(\ca L)$. Then $\g$ is a $\G$-form of $w$, and so 
$w\in\bold D_0^+$. Set $P=S_w$. The proof of the intermediary result (*) in the proof of 3.15 
shows that $w_\g\in\bold D$ via $N_P(U)$. 
Write $\phi_i=(x_i\i,h_i,y_i)$, and form the word $w_0\in\bold W(N_{\ca L}(R))$ as in 3.15(b). 
Set $g=\Pi(v)$. Then $g=\Pi(w_\g)=\Pi(x_1\i,w_0,y_n)$ by $\bold D$-associativity, and thus 
$(x_1\i,w_0,y_n)\in\Phi_g$. That is, we have $[g]=[x_1\i,w_0,y_n]$. This shows that 
$\l_0$ is a homomorphism of partial groups, completing the proof of (b). 
\qed
\enddemo 

We remark that it is easily verified that $Im(\l_0)$ is in fact a partial subgroup of $\ca L_0^+$. But 
there is no reason to suppose that $Im(\iota_0)$ is a partial subgroup of $\ca L$, as it may be the 
case that $\bold W(\ca L_0)\cap\bold D$ is not contained in $\bold D_0$.  

\vskip .1in 
For any $w=(g_1,\cdots,g_n)\in\bold W(\ca L)$ set $[w]=([g_1],\cdots,[g_n])$. For any subset $W$ of 
$\bold W(\ca L)$ set $[X]=\{[w]\mid w\in W\}$. Since each $\approx$-class $[g_i]$ intersects $\ca L$ in 
$\{g_i\}$ it follows that the product $\Pi:\bold D\to\ca L$ may be regarded as a mapping 
$[\bold D]\to\ca L$. 

Recall from Theorem I.1.17 that the category of partial groups has all colimits (and all limits). 
In particular, pushouts are available.

\proclaim {Proposition 3.20} Set $\bold D^+=[\bold D]\cup\bold D_0^+$. 
\roster 

\item "{(a)}" The products $\Pi:\bold D\to\ca L$ and $\Pi_0^+:\bold D_0^+\to\ca L_0^+$ agree on 
$[\bold D]\cap\bold D_0^+$, and $\Pi\cup\Pi_0^+$ is a mapping 
$$ 
\Pi^+:\bold D^+\to\ca L^+.  
$$ 

\item "{(b)}" $\ca L^+$ is a partial group via the product $\Pi^+$ and the involutory bijection 
given by 3.13. 

\item "{(c)}" Let $\l:\ca L\to\ca L^+$ be the mapping $g\maps[g]$, and let $\iota:\ca L_0^+\to\ca L^+$ 
be the inclusion map. Then $\l$ and $\iota$ are injective homomorphisms of partial groups, and 
$$ 
\CD 
\ca L_0^+  @>\iota>>  \ca L^+  \\
@A{\l_0}AA            @AA{\l}A   \\
\ca L_0   @>>\iota_0>  \ca L
\endCD 
$$ 
is a pushout diagram in the category of partial groups. 

\endroster
\endproclaim 

\demo {Proof} That $\Pi$ and $\Pi_0^+$ agree on $[\bold D]\cap\bold D_0^+$ is one way of interpreting 
the fact (3.17(b)) that $\l_0$ and $\iota_0$ are homomorphisms of partial groups. Thus (a) holds, and 
point (b) is then a straightforward exercise with definition I.1.1. 

Let $\ca L^*$ (with the appropriate diagram of homomorphisms) be a pushout for 
$$ 
\ca L_0^+@<\l_0<<\ca L_0@>\iota_0>>\ca L. \tag*
$$ 
By I.1.17 we may in fact take the underlying set of $\ca L^*$ to be the standard pushout of (*) as 
a diagram of mappings of sets. That is, we may take $\ca L^*$ to be the disjoint union 
$\ca L_0^+\bigsqcup\ca L$ modulo the relation $\equiv$ which identifies $g\in\ca L_0$ with $g\l_0$. Here 
$g\l_0=\Phi_g\cup\{g\}$, and the  elements of $\ca L$ which are not in $\ca L_0$ are by definition the 
singletons $\{f\}$ such that $S_f$ contains no $\ca F$-conjugate of $R$. By identifying such a singleton 
$\{f\}$ with its unique element we thereby obtain $\ca L^*=\ca L^+$ as sets. 

The domain $\bold D^*$ of the product in $\ca L^*$ is obtained by from the disjoint union 
$\bold D_0^+\bigsqcup\bold D$ by imposing the $\equiv$-relation componentwise. That is, we have 
$\bold D^*=\bold D^+$. The product $\Pi^*:\bold D^*\to\ca L^*$ is then the union of the products 
$\Pi_0^*$ and $\Pi$; which is to say that $\Pi^*=\Pi^+$. Similarly, the inversion maps on $\ca L^*$ 
and on $\ca L^+$ coincide, and so $\ca L^*=\ca L^+$ as partial groups. That $\l$ and $\iota$ are 
the homomorphisms which give the required pushout diagram in then immediate. 

Now let $f,g\in\ca L$ with $f\l=g\l$. Then 3.12 yields 
$$ 
\{f\}=[f]\cap\ca L=[g]\cap\ca L=\{g\} 
$$ 
and so $\l$ is injective. The inclusion map $\iota$ is of course injective, so the proof is complete. 
\qed 
\enddemo

Let $\D^+$ be the union of $\D$ with the set of all subgroups $P\leq S$ such that $P$ contains an 
$\ca F$-conjugate of $R$. The following lemma prepares the way for showing that $(\ca L^+,\D^+)$ is an 
objective partial group.

\proclaim {Lemma 3.21} Let $\l:\ca L\to\ca L^+$ be the homomorphism of partial groups given by $g\maps[g]$. 
Let $[\phi]\in\ca L_0^+$, and let $S_{[\phi]}$ be the set of all $a\in S$ such that $[a]^{[\phi]}$ is 
defined in $\ca L^+$, and such that $[a]^{[\phi]}\in [S]$. Then $S_{[\phi]}=S_\phi$. 
\endproclaim 

\demo {Proof} Let $a\in S$ such that $[a]^{[\phi]}=[b]$ for some $b\in S$. 
Suppose first that $[\phi]\cap\ca L$ is non-empty, and let $g$ be the unique element of 
$[\phi]\cap\ca L$.  The equality $[a]^{[\phi]}=[b]$ then simply means that 
$(a^g)\l=b\l$, and the injectivity of $\l$ (3.20(c)) yields $a^g=b$. Thus $a\in S_g$. Conversely, 
for any $x\in S_g$ we have $x\l=[x]\in S_{[\phi]}=g\l$, and thus the lemma holds in this case. 

Assume now that $[\phi]\cap\ca L=\nset$. Then $[\phi]$ is a $\sim$-class by 3.14, and 3.12 shows that the 
pair $(U,V):=(U_\phi,V_\phi)$ is constant over all $\phi\in[\phi]$. Set 
$w=([\phi]\i,[a],[\phi])$. Then $w\in\bold D_0^+$ by 
hypothesis, so there is a $\G$-form $\g=(\phi\i,\psi,\bar\phi)$ of $w$. This means that, upon 
setting $w\g=\phi\i\circ\psi\circ\bar\phi$, we have $V\leq S_{w_\g}$. The uniqueness of $(U,V)$ for 
$[\phi]$ (and of $(V,U)$ for $[\phi]\i$) yields $U=U_\psi=V_\psi$, and then $a\in N_S(U)$ 
since $a=\Pi(\psi)$. Since $[b]^{[\phi]\i}=[a]$ we similarly obtain $b\in N_S(V)$. Notice also 
that since $(U,V)$ is independant of the choice of representative for $[\phi]$ we may take $\bar\phi$ to 
be $\phi$. As $x\in\bY_U$ and $y\in\bY_V$ we obtain:  
\roster 

\item "{(*)}" $a^{x\i},b^{y\i}\in N_S(R)$. 

\endroster 

Write $\phi=(x\i,h,y)$, and set $\theta=(a\i x\i,\1,x a^2)$. 
As $a\in N_S(U)$ we get $xa^2\in\bold Y_U$ and $a\i x\i\in\bold X_U$. Thus 
$\theta\in\Phi$. Moreover, we have $\theta\in\bold D$ via $N_S(U)$, and $\Pi(\theta)=a$. Thus 
$\theta\in\Phi_a$, and $(\phi\i,\theta,\phi)$ is a $\G$-form of $w$. Then 
$$
\Pi^+(w)=[y\i,h\i(x(a\i x\i))((xa^2)x\i)h,y] 
$$
by the definition of $\Pi^+$ in 3.20. Observing now that 
$$ 
\text{$(x,a\i,x\i,x,a^2,x\i)\in\bold D$}, 
$$ 
via $N_S(U)^{x\i}$, we obtain 
$$
[b]=\Pi^+(w)=[y\i,h\i(xax\i)h,y]=[y\i,(a^{x\i})^h,y]. 
$$ 

We now claim that $(y\i,b^{y\i},y)\in [b]$. 
In order to see this, one observes first of all that $(y\i,b^{y\i},y)\in\Phi$. Also, since 
$b^{y\i}$ normalizes $N_S(V)^{y\i}$ we have $(y\i,y,b,y\i,y)\in\bold D$ via $N_S(V)$. Then 
$\Pi(y\i,b^{y\i},y)=b$, and the claim is proved. Thus: 
$$ 
[y\i,(a^{x\i})^h,y]=[y\i,b^{y\i},y]. \tag**
$$ 
Application of $\Pi$ to both sides of (**) then yields 
$$ 
((a^{x\i})^h)^y=(b^{y\i})^y,  
$$ 
and then $(a^{x\i})^h=b^{y\i}$ by the cancellation rule in $\ca L$. As $a^{x\i},b^{y\i}\in S$ by (*), we 
conclude that $a\in S_\phi$. This completes the proof of (a), and thereby completes the proof. 
\qed 
\enddemo

At this point it will be convenient (and need cause no confusion) to view $\l_0$ and $\l$ as inclusion maps. 
Then $\ca L^+=\ca L\cup\ca L_1^+$, where 
$$ 
\ca L_1^+=\{[\phi]\mid \phi\in\Phi,\ \phi\notin\bold D\}. 
$$

\proclaim {Proposition 3.22} $(\ca L^+,\D^+,S)$ is a locality. Moreover: 
\roster 

\item "{(a)}" $\ca L$ is the restriction $\ca L^+\mid_\D$ of $\ca L^+$ to $\D$. 

\item "{(b)}" $N_{\ca L^+}(R)=N_{\ca L}(R)$. 

\item "{(c)}" $\ca F_S(\ca L^+)=\ca F$. 

\endroster 
\endproclaim 

\demo {Proof} We have first to show that $(\ca L^+,\D^+)$ satisfies the conditions (O1) and (O2) in 
the definition I.2.1 of objectivity. Condition (O2) is the requirement that $\D^+$ be $\ca F$-closed 
(i.e. that $\D^+$ be preserved by $\ca F$-homomorphisms). Since $\D$ is $\ca F$-closed, and $\D^+$ is 
given by attaching to $\D$ an $\ca F$-conjugacy class $R^{\ca F}$ and all overgroups in $S$ of members 
of $R^{\ca F}$, (O2) holds for $(\ca L^+,\D^+)$. 

Condition (O1) requires that $\bold D^+$ be equal to $\bold D_{\D^+}$. This means: 
\roster 

\item "{(*)}" The word $w=(C_1,\cdots,C_n)\in\bold W(\ca L^+)$ is in $\bold D^+$ $\iff$ there exists 
a sequence $(X_0,\cdots,X_n)\in\bold W(\D^+)$ such that $X_{i-1}^{C_i}=X_i$ for all $i$, $(1\leq i\leq n)$. 

\endroster 
Here we need only be concerned with the case $w\in\bold W(\ca L_0^+)$ since 
$\bold D^+=\bold D_0^+\cup[\bold D]$, and since $\bold D=\bold D_\D$. The implication $\implies$ in (*) is 
then given by the definition of $\bold D_0^+$, with $X_0\in R^{\ca F}$. The reverse implication is given 
by 3.21, and thus $(\ca L^+,\D^+)$ is objective. We note that $\ca L^+$ is finite since $\ca L$ and $\Phi$ 
are finite. 

Let $\w S$ be a $p$-subgroup of $\ca L^+$ containing $S$, and let $a\in N_{\w S}(S)$. As $S\in\D$ we get 
$a\notin\ca L_1^+$, so $a\in N_{\ca L}(S)$, and then $a\in S$ since $S$ is a maximal $p$-subgroup of $\ca L$. 
Thus $\w S=S$, $S$ is a maximal $p$-subgroup of $\ca L^+$, and $(\ca L^+,\D^+,S)$ is a locality. The 
restriction of $\ca L^+$ to $\D$ is by definition the partial group whose product is the restriction 
of $\Pi^+$ to $\bold D_\D$, whose underlying set is the image 
of $\Pi^+\mid\bold D_\D$, and whose inversion map is inherited from $\ca L^+$. That is, (a) holds. 

Let $[\phi]\in\ca L_1^+$, let $\phi=(x\i,h,y)\in[\phi]$, and set $U=U\phi$. Then $U\norm S_{[\phi]}=S_\phi$ 
by 4.17(a), and the conjugation map $c_{[\phi]}:S_\phi\to S$ is then the composite 
$c_x\i\circ c_h\circ c_y$ applied to $S_\phi$. Thus $c_{[\phi]}$ is an $\ca F$-homomorphism, and this 
yields (c). Suppose now that $[\phi]\in N_{\ca L^+}(R)$. Then $x,h$, and $y$ are in $N_{\ca L}(R)$, 
and then $\phi\in\bold D$ as $N_{\ca L}(R)$ is a subgroup of $\ca L$. Then $[\phi]\notin\ca L_1^+$, 
and so (b) holds. 
\qed 
\enddemo  

\proclaim {Proposition 3.23} Let $(\w{\ca L},\D^+,S)$ be a locality having the same set $\D^+$ of 
objects as $\ca L^+$. Assume that $\w{\ca L}\mid_\D=\ca L$ and that $N_{\w{\ca L}}(R)=N_{\ca L}(R)$. 
Then the identity map on $\ca L$ extends in a unique way to an isomorphism $\ca L^+\to\w{\ca L}$. 
\endproclaim 

\demo {Proof} Write $\w\Pi:\w{\bold D}\to\w{\ca L}$ for the product in $\w{\ca L}$. (It isn't necessary 
to distinguish the inversion map on $\w\Pi$ in any way, since by I.1.13 it restricts to the inversion map 
on $\ca L$.) Let $\phi=(x\i,h,y)$ 
and $\bar\phi=(\bar x\i,\bar h,\bar y)$ be members of $\Phi$ such that $\phi\sim\bar\phi$. Set 
$U=U_\phi$. Then $U=U_{\bar\phi}$, and $(\bar x,x\i)\in\bold D$ via $N_S(U)^{\bar x\i}$. Then 
$\w\Pi(\bar x,x\i)=\Pi(\bar x,x\i)$ as $\ca L$ is the restriction of $\w{\ca L}$ to $\D$. Similarly, 
we obtain $\w\Pi(\bar y,y\i)=\Pi(\bar y,y\i)$. Then: 
$$ 
\w\Pi(\bar x,x\i,h)=\w\Pi(\bar xx\i,h)=\Pi(\bar xx\i,h)=\Pi(\bar h,\bar yy\i)=\w\Pi(\bar h,\bar y,y\i).
\tag1  
$$ 
Observe that both $(\bar x\i,\bar x,x\i,h,y)$ and $(\bar x\i,\bar h,\bar y,y\i,y)$ are in 
$\bold D^+\cap\w{\bold D}$ (via the obvious conjugates of $U$). Then (1) yields  
$$ 
\w\Pi(\phi)=\w\Pi(\bar x\i,\bar x,x\i,h,y)=\w\Pi(\bar x\i,\bar h,\bar y,y\i,y)=\w\Pi(\bar\phi).  
$$ 
If $\phi\in\bold D$ then $\w\Pi(\phi)=\Pi(\phi)$, so we have shown that there is a well-defined 
mapping 
$$ 
\b:\ca L^+\to\w{\ca L} 
$$ 
such that $\b$ restricts to the identity map on $\ca L$ and sends the element $[\phi]\in\ca L_0^+$ to 
$\w\Pi(\phi)$.  

We now show that $\b$ is a homomorphism of partial groups. As $\ca L^+$ is a pushout, in the manner 
described in 3.20(c), and since the inclusion $\ca L\to\w{\ca L}$ is a homomorphism, it suffices to 
show that the restriction $\b_0$ of $\b$ to $\ca L_0^+$ is a homomorphism. 

So, let $w=([\phi_1],\cdots,[\phi_n])\in\bold D^+$, let $\g=(\phi_1,\cdots,\phi_n)$ be 
a $\G$-form of $w$, and set $w_\g=(\phi_1,\cdots,\phi_n)$. Write $\phi_i=(x_i\i,h_i,y_i)$. Then 
$\Pi^+(w)=[x_1\i,\Pi(w_0),y_n]$, where $w_0\in\bold W(N_{\ca L}(R))$ is the word 
$$ 
w_0=(h_1,y_1x_2\i,\cdots,y_{n-1}x_n\i,h_n). 
$$
given by 3.17(b). Let $\b*$ be the induced mapping $\bold W(\ca L^+)\to\bold W(\w{\ca L})$ of free monoids. 
Then 
$$ 
\w\Pi(w\b^*)=\w\Pi([\phi_1]\b,\cdots,[\phi_n]\b)=\w\Pi(\w\Pi(\phi_1),\cdots,\w\Pi(\phi_n))=\w\Pi(w_\g) 
$$ 
by $\w{\bold D}$-associativity, and  
$$ 
(\Pi^+(w))\b=[x_1\i,\Pi(w_0),y_n]\b=\w\Pi(x_1\i,\Pi(w_0),y_n).       
$$ 
Set $w_0'=(h_1,y_1,x_2\i,\cdots,y_{n-1},x_n\i,h_n)$. Then $w_0'\in\w{\bold D}$ via $U_{\phi_1}$, and 
$\w\Pi(w_0')=\w\Pi(w_0)=\Pi(w_0)$. Then  
$$
(\Pi^+(w))\b=\w\Pi((x_1\i)\circ\w\Pi(w_0')\circ(y_n))=\w\Pi(w_\g)=\w\Pi(w\b^*),  
$$ 
and so $\b_0$ is a homomorphism. As already mentioned, this result implies that $\b$ is a homomorphism. 

Let $f\in\w{\ca L}$ with $f\notin\ca L$. Then $S_f\notin\D$ (as $\ca L=\w{\ca L}\mid_\D$), and 
so $S_f$ contains a unique $U\in R^{\ca F}$. Set $V=U^f$, and recall that the inversion map on $\ca L$ 
is induced from the inversion map on $\w{\ca L}$. If $V\in\D$ then $f\i\in\w{\ca L}$, and then $f\in\ca L$. 
Thus $V\notin\D$, and hence $V\in R^{\ca F}$. By 3.6 there exist elements $x,y\in\ca L$ such that $U^x=R=V^y$ 
and such that $N_S(U)^{x\i}\leq N_S(R)\geq N_S(V)^{y\i}$. We find that 
$(x,f,y\i)\in\w{\bold D}$ via $R$ and that $h:=\w\Pi(x,f,y\i)\in N_{\w{\ca L}}(R)$. Then $h\in N_{\ca L}(R)$ 
by hypothesis, and $f=\w\Pi(x\i,h,y)$. In particular, we have shown: 
\roster 

\item "{(*)}" Every element of $\w{\ca L}$ is a product of elements of $\ca L$. 

\endroster 
We may now show that $\b^*$ maps $\bold D_0^+$ onto $\w{\bold D}$. Thus, let 
$\w w=(f_1,\cdots,f_n)\in\w{\bold D}$ and set $X=S_{\w w}$. If $X\in\D$ then $\w w\in\bold D$ and 
$\w w=\w w\b^*$. So assume that $X\notin\bold D$. Then $X$ contains a unique $\ca F$-conjugate $U_0$ 
of $R$, and there is a sequence $(U_1,\cdots,U_n)\in\bold W(\D^+)$ given by $U_i=(U_{i-1})^f_i$. As 
seen in the preceding paragraph, we then have $U_i\in R^{\ca F}$ for all $i$, and there exists a 
sequence $\g=(\phi_1,\cdots,\phi_n)\in\bold W(\Phi)$ such that $\w\Pi(\phi_i)=f_i$. Set 
$w_\g=[\phi_1]\circ\cdots\circ[\phi_n]$. Then $w_\g\in\bold D^+$ and $(w_\g)\b^*=\w w$. Thus
$\b^*:\bold D_0^+\to\w{\bold D}$ is surjective. That is, $\b$ is a projection, as defined in I.4.5. 

Let $g\in Ker(\b)$. If $g\in\ca L$ then $g=\1$ since $\b\mid_{\ca L}$ is the inclusion map. So assume 
that $g\notin\ca L$. Then $g=[\phi]$ for some $\phi=(x\i,h,y)\in\Phi$, and $[\phi]$ is a $\sim$-class. 
Set $U=U_\phi$ and $V=V_\phi$. As $\1=g\b=\w\Pi(\phi)$ it follows that $U=V=R$. But then $\phi\in\bold D$ 
as $N_{\ca L}(R)$ is a subgroup of $\ca L$, and so $[\phi]$ is not a $\sim$-class. Thus $Ker(\b)=\1$. 
As $\b$ is a projection, $\b$ is then an isomorphism by I.4.3(d). It follows from (*) that $\b$ is the 
unique isomorphism $\ca L^+\to\w{\ca L}$ which restricts to the identity on $\ca L$. 
\qed 
\enddemo 

Notice that Propositions 3.22 and 3.23 complete the proof of Theorem 3.8.

\vskip .2in 
\noindent 
{\bf Section 4: Elementary expansions and partial normal subgroups} 
\vskip .1in

We continue to assume Hypothesis 3.7 throughout this section. Our aim is to prove the following result.

\proclaim {Theorem 4.1} Assume Hypothesis 3.7, and let $(\ca L^+,\D^+,S)$ be the expansion of ($\ca L,\D,S)$ 
given by Theorem 3.8. Let $\ca N\norm\ca L$ be a partial normal subgroup of $\ca L$, and let 
$\Omega:=\ca N^{\ca L^+}$ be the set of all $\ca L^+$-conjugates of elements of $\ca N$ (the set of all 
$f^g$ such that $f\in\ca N$, $g\in\ca L^+$ and $(g\i,f,g)\in\bold D^+$). Let 
$\ca N^+$ be partial subgroup of $\ca L^+$ generated by $\Omega$ (cf. I.1.9). Then:   
\roster 

\item "{(a)}" $\ca N^+\norm\ca L^+$.  

\item "{(b)}" $\ca N^+\cap\ca L=\ca N$. 

\item "{(c)}" If $\ca M\norm\ca L^+$ is a partial normal subgroup of $\ca L^+$ such that 
$\ca M\cap\ca L=\ca N$ then $\ca M=\ca N^+$. 

\endroster 
\endproclaim

We will employ all of the notation from section 3. Thus, the reader will need to have in mind the 
meanings of $\ca L_0$, $\ca L_0^+$, $\ca L^+$, $\Phi$, $\sim$, and $\approx$. Set $\Omega=\ca N^{\ca L^+}$.

\proclaim {Lemma 4.2} Let $w\in\bold D^+$, and suppose that $w\notin\bold W(\ca L)$. 
Then $S_w$ is an $\ca F$-conjugate of $R$, and $S_w=S_{\Pi^+(w)}$ if $\Pi^+(w)\notin\ca L$. 
\endproclaim 

\demo {Proof} As $w\notin\ca L$ it follows that $S_w\notin\D$. Then, as $w\in\ca L^+$, it follows from 
3.5(1) that $S_w$ is an $\ca F$-conjugate of $R$. Write $w=(g_1,\cdots,g_n)$ with $g_i\in\ca L^+$, 
set $g=\Pi^+(w)$, and suppose that $g\notin\ca L$. Then $S_g$ is an $\ca F$-conjugate, and since 
$S_w\leq S_g$ we obtain $S_w=S_g$, as required.  
\qed 
\enddemo

\proclaim {Lemma 4.3} Assume that 4.1(b) holds. Then: 
\roster 

\item "{(a)}" Theorem 4.1 holds in its entirety. 

\item "{(b)}" Let $\r:\ca L\to\ca L/\ca N$ and $\r^+:\ca L^+/\ca N^+$ be the canonical projections, and 
let $\iota:\ca L\to\ca L^+$ be the inclusion map. Then there is a unique homomorphism 
$\eta:\ca L/\ca N\to\ca L^+/\ca N^+$ such that $\iota\circ\r^+=\r\circ\eta$. Moreover, $\eta$ is injective, 
and $\ca L^+/\ca N^+$ is an elementary expansion of $\ca L/\ca N$ via $\eta$ (cf. 3.1). 
Further, if $S\cap\ca N\nleq R$ 
then $\eta$ is an isomorphism, while if $S\cap\ca N\leq R$ then 
$N_{\ca L/\ca N}(R\eta)=N_{\ca L^+/\ca N^+}(R\eta)$. 

\endroster 
\endproclaim 

\demo {Proof} We assume that $\ca L\cap\ca N^+=\ca N$ (4.1(b)). Set $\Omega_0=\Omega$, 
and recursively define $\Omega_n$ for $n>0$ by 
$$ 
\Omega_n=\{\Pi^+(w)\mid w\in\bold W(\Omega_{n-1})\cap\bold D^+\}. 
$$ 
Then $\ca N^+=\<\Omega\>$ is the union of the sets $\Omega_n$, by I.1.9. In order to show that 
$\ca N^+\norm\ca L^+$ it then suffices to show that each $\Omega_n$ is closed with respect to 
conjugation (whenever defined) in $\ca L^+$. 

Let $g\in\Omega_0$. Then there exists $f\in\ca N$ and $a\in\ca L^+$ such that $(a\i,f,a)\in\bold D^+$ 
and with $g=f^a$. Now let $b\in\ca L^+$ such that $(b\i,g,b)\in\bold D^+$. If $g\in\ca L$ then 
$g\in\ca N$ by 4.1(b), and then $g^b\in\Omega_0$. On the other hand, assume $g\notin\ca L$. Then 
$S_g\in R^{\ca F}$, and 4.2 yields $S_g=S_{(a\i,f,a)}$. Set $u=(b\i,g,b)$ and $v=(b\i,a\i,f,a,b)$. 
Then $S_u=S_v=(S_g)^b$, so $v\in\bold D^+$, and $g^b=\Pi^+(u)=\Pi^+(v)=f^{ab}$. Thus $g^b\in\Omega_0$ in 
any case, and so $\Omega_0$ is closed with respect to conjugation in $\ca L^+$. 

Assume now that there exists a largest integer $k$ such that $\Omega_k$ is closed with respect to 
conjugation in $\ca L^+$, and let now $f\in\Omega_{k+1}$ and $b\in\ca L^+$ such that $f^b$ is defined 
in $\ca L^+$ and is not in $\Omega_{k+1}$. By definition there exists 
$w=(f_1,\cdots,f_n)\in\bold W(\Omega_k)\cap\bold D^+$ such that $f=\Pi^+(w)$. Suppose $f\in\ca L$. Then 
$f\in\ca N$ by 4.1(b), and then $f^b\in\Omega$. As $\Omega\sub\Omega_m$ for all $m$, we have 
$f^b\in\Omega_{k+1}$ in this case. So assume that $f\notin\ca L$. Then $S_f=S_w\in R^{\ca F}$ by 
4.2. Set $w'=(b\i,f_1,b,\cdots,b\i,f_n,b)$. Then $w'\in\bold D^+$ via $(S_w)^b$, and we obtain 
$f^b=\Pi^+(w')=\Pi^+(f_1^b,\cdots,f_n^b)$. As $\Omega_k$ is closed with respect to conjugation in $\ca L^+$ 
we conclude that $f^b\in\Omega_{k+1}$, contrary to the maximality of $k$; and this contradiction 
completes the proof that $\ca N^+\norm\ca L^+$. That is, 4.1(a) holds. 

Next: let $\ca M\norm\ca L^+$ be any partial normal subgroup of $\ca L^+$ such that $\ca M\cap\ca L=\ca N$. 
Let $\phi$ be the composite of the inclusion map $\ca L\to\ca L^+$ followed by the quotient map 
$\r:\ca L^+\to\ca L^+/\ca M$. Then $Ker(\r)=\ca M$ (cf. I.4.3), and so $Ker(\phi)=\ca M\cap\ca L=\ca N$. 
Let $\r_S$ be the restriction of $\r$ (and hence also of $\phi$) to $S$. 

Set $\bar\D=\{P\phi\mid P\in\D\}$, and write $\bar P$ for $P\phi$. Let $\bar w\in\bold W(\ca L^+/\ca M)$, 
and set $\bar X=\bar S_{\bar w}$. That is (and in the usual way) $\bar X$ is the largest subgroup of 
$\bar S$ which is conjugated successively into $\bar S$ by the entries of $\bar w$. Write 
$\bar w=(\bar g_1,\cdots,\bar g_n)$, and for each $i$ let $g_i\in\ca L^+$ be a representative of $\bar g_i$ 
such that $g_i$ is $\up$-maximal with respect to $\ca M$. Set $w=(g_1,\cdots,g_n)$, and let  
$X$ be the preimage of $\bar X$ in $S$ via $\r_S$. Then $X\leq S_w$ by I.3.14(e).  
If $\bar X\in\bar\D$ then $X\in\D$ and $w\in\bold D(\ca L^+\mid_\D)$. As $\ca L^+\mid\D=\ca L$, we have 
$w\in\bold D(\ca L)$ in this case, and then $\bar w\in\bold D((\ca L^+/\ca M)\mid_{\bar D})$. 
On the other hand, if $\bar w\in\bold D((\ca L^+/\ca M)\mid_{\bar D})$ then $\bar X\in\bar\D$. Thus: 
$$ 
Im(\phi)=(\ca L^+/\ca M)\mid_{\bar D}.  
$$ 
In particular, $Im(\phi)$ is a partial group (and indeed a locality) in this way, and the mapping 
$$ 
\phi_0:\ca L\to Im(\phi) 
$$ 
given by $\phi$ is then a homomorphism of partial groups. One observes that the induced mapping 
$$ 
\phi_0^*:\bold W(\ca L)\to\bold W(Im(\phi)) 
$$ 
of free monoids carries $\bold D(\ca L)$ onto $\bold W(Im(\phi))$. That is, $\phi_0$ is a projection.  
Then, by I.4.6 there is an isomorphism $\bar\phi:\ca L/\ca N\to Im(\phi)$ induced by $\phi$. 

Let $\eta:\ca L/\ca N\to\ca L^+/\ca M$ be the composition of $\bar\phi$ followed by the inclusion of 
$Im(\phi)$ into $\ca L^+/\ca M$. Our aim now is to show that $\ca L^+/\ca M$ is an expansion of $\ca L/\ca N$ 
via $\eta$. For this, what is required (cf. definition 3.1) 
is: (A) $\eta$ is injective (obvious, since 
$\eta$ is a composite of injective mappings); (B) $Im(\eta)$ is the restriction of $\ca L^+/\ca M$ to $\D\phi$ 
(which has been verified above);  and (C) $\D\phi$ is closed in the fusion system 
$\ca F'=\ca F_{\bar S}(\ca L^+/\ca M)$. In order to verify (C) we note first of all that 
$\ca F=\ca F_S(\ca L^+)$ by 3.4(a). Since $\r$ and $\phi$ are projections, (C) then follows from 1.20.  
Thus $\eta$ is an expansion, as claimed. Moreover, $\eta$ is an elementary expansion, since 1.20 
shows that the set $\{\bar P\mid P\in\D^+\}$ of objects of $\ca L^+/\ca M$ is the union of $\bar D$ with 
the set of $\ca F_{\bar S}(\ca L^+/\ca M)$-conjugates of $\bar R$. 

We have now to verify that Hypothesis 3.3 holds with respect to $\ca L/\ca M$ and $\bar R$, in order to 
apply Theorem 3.4 to the expansion $\ca L^+/\ca M$ of $Im(\phi)$. Set $T=S\cap\ca N$, and suppose first that 
$T\nleq R$. Then $TR\in\D$, so $\bar R\in\bar\D$, and then 3.4 holds trivially. On the other hand, suppose 
that $T\leq R$. Then every subgroup $\bar Q$ of $\bar S$ which properly contains an 
$\ca F_{\bar S}(\ca L^+/\ca M)$-conjugate of $\bar R$ is the image under $\phi$ of a subgroup $Q$ of $S$ 
which properly contains an $\ca F$-conjugate of $R$. That is, 3.3(1) obtains. Now let 
$\bar g\in N_{\ca L^+/\ca M}(\bar R)$, and let $g\in\ca L^+$ be a $\r$-preimage of $g$ such that $g$ is 
$\up$-maximal with respect to $\ca M$. Then $g\in N_{\ca L^+}(R)$ by I.3.14(e). Here 
$N_{\ca L^+}(R)=N_{\ca L}(R)$ by 3.4 applied to $\ca L^+$, so we have shown that $\r$ maps 
$N_{\ca L}(R)$ onto $N_{\ca L^+/\ca M}(\bar R)$. This verifies 3.3(2). 

The preceding analysis applies to $\ca N^+$ in the role of $\ca M$. We now observe that 
$\ca N^+\leq\ca M$ since $\ca N\sub\ca M$, and since $\ca N^+=\<\ca N^{\ca L^+}\>$. 
Let $\eta_0:\ca L/\ca N\to\ca L^+/\ca N^+$ be the homomorphism 
which defines the expansion of $\ca L/\ca N$ to $\ca L^+/\ca N^+$. By Theorem 3.4(b) there is then a  
unique homomorphism $\b:\ca L/\ca N^+\to\ca L/\ca M$ such that $\eta_0\circ\b=\eta$, and $\b$ is in fact 
an isomorphism. On the other hand, let $\g:\ca L/\ca N^+\to\ca L/\ca M$ be the mapping which sends a 
maximal coset $[g]_{\ca N^+}$ of $\ca N^+$ in $\ca L^+$ to the unique maximal coset $[g]_{\ca M}$ containing 
it. One easily verifies that $\g$ is a homomorphism and that $\eta_0\circ\g=\eta$, and thus $\g$ is an 
isomorphism. For $g\in\ca M$ we then have $[g]_{\ca N^+}=\ca M$, and this shows that $\ca M=\ca N^+$. 
That is, 4.1(c) holds, and we have obtained point (a) of the lemma. Point (b) was obtained in the course 
of proving (a).  
\qed 
\enddemo

\definition {(4.4) Notation} Set $M=N_{\ca L}(R)$, and $K=M\cap\ca N$. As in the preceding section,  
for each $g\in\ca L^+$ set 
$$ 
\ca U_g=\{U\in R^{\ca F}\mid U\leq S_g\},  
$$ 
and for each $U\in R^{\ca F}$ set 
$$ 
\bY_U=\{y\in\ca L\mid R^y=U,\ N_S(U)^{y\i}\leq N_S(R)\}. 
$$ 
\enddefinition

The notation (4.5) will remain fixed until the proof of Theorem 4.1 is complete. Note that since $M$ is a 
subgroup of $\ca L$ by 3.3, $K$ is a normal subgroup of $M$ by I.1.8(c). 
\vskip .1in

\proclaim {Lemma 4.6} Suppose $T\leq R$. Then $\ca N^+=\Omega$, and $\Omega\cap\ca L=\ca N$.  
\endproclaim 

\demo {Proof} Let $f\in\ca N$, and suppose that $S_f$ contains an $\ca F$-conjugate of $R$. Then 
$T\leq S_f$ since, by I.3.1(a) $T$ is weakly closed in $\ca F$. Then I.3.1(b) yields the following result. 
\roster 

\item "{(1)}" Let $f\in\ca N$ such that $S_f$ contains an $\ca F$-conjugate $U$ of $R$. Then 
$P^f=P$ for each subgroup $P$ of $S_f$ containing $U$. In particular, $U^f=U$. 

\endroster 

Let $f'\in\Omega$, and let $f\in\ca N$ and $g\in\ca L^+$ such that $f'=f^g$. That is, assume that 
$v:=(g\i,f,g)\in\bold D^+$ and that $f'=\Pi^+(v)$. If $v\in\bold D$ then $f'=\Pi(v)\in\ca N$. On the other 
hand, suppose that $v\notin\bold D$. Then $S_v$ is an $\ca F$-conjugate of $R$ by 4.1(1). Set 
$U=(S_v)^{g\i}$. Then $U=U^f$ by (1). Now choose $a\in\bY_U$, and set $v'=(g\i,a\i,a,f,a\i,a,g)$. 
Then $v'\in\bold D^+$ via $U^g$, and $\Pi^+(v)=\Pi^+(v')$. Notice that (1) implies that 
$(a,f,a\i)\in\bold D$ via $P:=N_{S_f}(U)$, and that 
$$ 
T@>a>>U@>f>>U@>a\i>>T, 
$$ 
so that $f^{a\i}\in M$. Thus $f^{a\i}\in K$, and $\Pi^+(v)=\Pi^+(g\i a\i,f^{a\i},ag)$. This shows:   
\roster 

\item "{(2)}" $\Omega$ is the union of $\ca N$ with the set of all $\Pi^+(g\i,f,g)$ such that $f\in K$ 
and such that $(g\i,f,g)\in\bold D^+$. Moreover, for any such $v=(g\i,f,g)$, $\Pi^+(v)$ normalizes 
each $V\in R^{\ca F}$ such that $V\leq S_v$. 

\endroster 

Assume now that we have $v=(g\i,f,g)$ as in (2) (so that $f\in K$), and let $A$ be an $\ca F$-conjugate 
of $R$ contained in $S_v$. In order to analyze these things further, we shall need to be able to compute 
products in $\ca L^+$ in the manner described in 4.13 and 4.14. To that end, note first of all that since 
$A^g=T$ there exists a unique $h\in M$ and $y\in\bY_A$ such that $g\approx(\1,h,y)\in\Phi$, by 4.7. 
Set $\phi=(\1,h,y)$ and set $\psi=(\1,f,\1)$. Then $(\phi\i,\psi,\phi)$ is a $\G$-form of $(g\i,f,g)$, 
as defined in 4.12. We then compute via 4.13 that 
$$ 
f'=\Pi^+(g\i,f,g)=[y\i,h\i,\1][\1,f,\1][\1,h,y]=[y\i,f^h,y].  
$$ 
\roster 

\item "{(3)}" Let $f'\in\Omega$, such that $S_{f'}$ contains an $\ca F$-conjugate of $R$. Then $f'$ is 
an $\approx$-class $[y\i,k,y]$ with $k\in K$. 

\endroster 
If $f'\in\ca L$ then $(y\i,k,y)\in\bold D$ by 4.10, and so $f'\in\ca N$. Thus: 
\roster 

\item "{(4)}" $\Omega\cap\ca L=\ca N$. 

\endroster

Now let $w=(f_1',\cdots,f_n')\in\bold W(\Omega)\cap\bold D^+$, and set $B=S_w$. 
Suppose that $\Pi^+(w)\notin\Omega$. Then $\Pi^+(w)\notin\ca N$, so $\Pi^+(w)\notin\ca L$ by (4). 
Thus $w\notin\bold D$, so $B\in R^{\ca F}$, and then (2) shows that each $f_i'$ normalizes $B$. 
Fix $b\in\bold Y_B$. Then (3) implies that there exist elements $k_i\in K$ such that 
$f_i'=[b\i,k_i,b]$. One observes that the sequence of elements $(b\i,k_i,b)$ of $\Phi$ is a 
$\G$-form for $w$, and then 4.13 yields $\Pi^+(w)=[b\i,k,b]$ where $k=\Pi(k_1,\cdots,k_n)\in K$. 
This simply means that $\Pi^+(w)=k^b$, since $b\i=[b\i,\1,\1]$, $k=[\1,k,\1]$, and $b=[\1,\1,b]$. 
Thus $\Pi^+(w)\in\Omega$. Since $\Omega$ is closed under inversion, we have thus shown that 
$\Omega$ is a partial subgroup of $\ca L$. 

Finally, let $c\in\ca L^+$ be given so that $(c\i,f',c)\in\bold D^+$ (and where $f'\in\Omega$). Suppose  
that $\Pi^+(c\i,f',c)\notin\Omega$. Then $f'\notin\ca N$, so $f'\notin\ca L$ by (4), and then 
$f'=\Pi^+(g\i,f,g)$ for some $f\in\ca N$ and some $g\in\ca L^+$, and where $S_{f'}=S_{(g\i,f,g)}$. 
Set $u=(c\i,g\i,f,g,c)$. Then $u\in\bold D^+$ via $(S_{f'})^c$, and $\Pi^+(u)=\Pi^+(c\i g\i,f,gc)$. 
Thus $\Pi^+(c\i,f',c)\in\Omega$ after all, and $\Omega\norm\ca L^+$. 
\qed 
\enddemo

Let $\bar L$ be the quotient locality $\ca L/\ca N$ (cf. 4.4), let $\r:\ca L\to\bar{\ca L}$ be the 
quotient map ,and let $\r^*$ be the induced homomorphism $\bold W(\ca L)\to\bold W(\bar{\ca L})$ of 
free monoids.  For any subset or element $X$ of $\ca L$, $\bar X$ shall denote the image of $X$ 
under $\r$. We extend this convention to subsets and elements of $\bold W(\bar{\ca L})$ in the obvious 
way. Set $\bar\D=\{\bar P\mid P\in\D\}$. 

\proclaim {Lemma 4.7} Assume that $T\nleq R$, and set $\bar{\ca F}=\ca F_{\bar S}(\bar{\ca L})$. Then 
$\bar R^{\bar{\ca F}}\sub\bar\D$. 
\endproclaim 

\demo {Proof} Let $U\in R^{\ca F}$. As $T\nleq R$ we then have $T\nleq U$ by I.3.1(a). Then 
$U$ is a proper subgroup of $UT$, so $UT\in\D$ by 3.3(1). Then $\bar U=\bar{UT}\in\bar\D$ by 4.3.  
\qed 
\enddemo

\proclaim {Lemma 4.8} Assume that $T\nleq R$. There is then a homomorphism $\s:\ca L^+\to\bar{\ca L}$ 
such that the restriction of $\s$ to $\ca L$ is the quotient map $\r$. 
\endproclaim 

\demo {Proof} Set $\bar{\bold D}=\bold D(\bar{\ca L})$ and let $\bar\Pi:\bar{\bold D}\to\bar{\ca L}$ be 
the product in $\bar{\ca L}$. As $\bar R^{\bar{\ca F}}\sub\bar\D$ by 4.7, we have 
$\bar\Phi\sub\bar{\bold D}$, and so there is a mapping $\l:\Phi\to\bar{\ca L}$ given by 
$\phi\l=\bar\Pi(\phi\r^*)$, where $\r^*:\bold W(\ca L)\to\bold W(\bar{\ca L})$ is the homomorphism 
of free monoids induced by $\r$. That is: $\phi\l=\bar\Pi(\bar\phi)$. 

Let $\phi_1,\phi_2\in\Phi$ with $\phi_1\sim\phi_2$, and write 
$\phi_i=(x_i\i,h_i,y_i)$. Then $(x_2 x_1\i)h_1=h_2(y_2 y_1\i)$ (in $M$), and so 
$$ 
(\bar x_2\bar x_1\i)\bar h_1=\bar h_2(\bar y_2\bar y_1\i)\tag 1
$$ 
in $\bar M$. As $(x_2,x_1\i,h_1)$ and $(h_2,y_2,y_1\i)$ are in $\bold D^+$ via the appropriate conjugates 
of $R$, we have $(\bar x_2,\bar x_1\i,\bar h_1)$ and $(\bar h_2,\bar y_2,\bar y_1\i)$ in $\bar{\bold D}$, 
and then 
$$ 
\bar\Pi(\bar x_2,\bar x_1\i,\bar h_1)=\bar\Pi(\bar h_2,\bar y_2,\bar y_1\i) 
$$ 
by (1) and $\bar{\bold D}$-associativity. A standard cancellation argument (cf. I.2.4(a)) then yields 
$\bar\Pi(\bar\phi_1)=\bar\Pi(\bar\phi_2)$, and thus $\l$ is constant on $\sim$-classes. 

Now suppose that $\phi\in\Phi_g$ for some $g\in\ca L$. That is, suppose that $\phi\in\Phi\cap\bold D$ 
and that $g=\Pi(\phi)$. Then $\phi\l=\bar g$ since $\r$ is a homomorphism of partial groups, and 
so $\l$ is constant on $\approx$-classes. This shows that there is a (well-defined) mapping 
$\s:\ca L^+\to\bar{\ca L}$ given by $\r$ on $\ca L$ and by $\l$ on $\ca L_0^+$. It remains only to 
check that $\s$ is a homomorphism. 

Let $w=(f_1,\cdots,f_n)\in\bold D^+$. If $w\in\bold D$ then $w\s^*=w\r^*\in\bar{\bold D}$ and 
$\bar\Pi(w\s^*)=(\Pi(w))\s$. On the other hand, suppose that $w\notin\bold D$. Then $f_i=[\phi_i]$ 
for some $\phi_i\in\Phi$, and there is a $\G$-form $\g=(\phi_1,\cdots,\phi_n)$ of $w$. Thus the word 
$w_\g=\phi_1\circ\cdots\circ\phi_n$ has the property that $S_{w_\g}$ is an $\ca F$-conjugate of $R$, 
and hence $\bar{w_\g}\in\bar{\bold D}$. Then 
$$ 
\bar\Pi((w_\g)\s^*)=\bar\Pi(\bar{w_\g}). 
$$ 
Write $\phi_i=(x_i\i,h_i,y_i)$. Then 
$$ 
(\Pi^+(w_\g))\s=[x_1\i,\Pi(w_0),y_n]\s 
$$ 
where $w_0\in\bold D(M)$ is given by the formula in 4.13(b). One then obtains 
$\bar\Pi((w_\g)\s^*)=(\Pi^+(w_\g))\s$ via $\bar{\bold D}$-associativity, and so $\s$ is a homomorphism, 
as desired. 
\qed 
\enddemo

\proclaim {Lemma 4.9} $\ca N^+\cap\ca L=\ca N$. 
\endproclaim 

\demo {Proof} Let $\s:\ca L^+\to\bar{\ca L}$ be a homomorphism, as in 4.8, whose restriction to 
$\ca L$ is the quotient map $\r:\ca L\to\bar{\ca L}$. Then $Ker(\s)\cap\ca L=Ker(\r)$, where 
$Ker(\r)=\ca N$, by 4.4. As $Ker(\s)\norm\ca L^+$, by 1.14, the lemma follows. 
\qed 
\enddemo

Theorem 4.1 follows at once from the union of lemmas 4.3, 4.6, and 4.9.

\vskip .2in 
\noindent 
{\bf Section 5: Theorem A} 
\vskip .1in

Recall from 1.8 that for any fusion system $\ca F$ on $S$, $\ca F^s$ denotes the set of 
$\ca F$-subcentric subgroups of $S$. These are the subgroups $U\leq S$ 
such that there exists an $\ca F$-conjugate $V$ of $U$ with $V$ fully normalized in $\ca F$ and with 
$O_p(N_{\ca F}(V))\in\ca F^c$.

\proclaim {Lemma 5.1} Let $(\ca L,\D,S)$ be a proper locality on $\ca F$, and let $R$ be a subgroup of $S$ 
which satisfies 3.3(1) and 3.3(2). Assume that $R\in\ca F^s$. Then $R$ satisfies 3.3(3). That is, 
$N_{\ca L}(R)$ is a subgroup of $\ca L$ of characteristic $p$, and 
$N_{\ca F}(R)=\ca F_{N_S(R)}(N_{\ca L}(R))$. 
\endproclaim 

\demo {Proof} Set $\ca L_R=N_{\ca L}(R)$ and set $\D_R=\{P\in\D\mid R\norm P\}$. By 3.2 
$(\ca L_R,\D_R,N_S(R))$ is a locality, and we set $\ca F_R=\ca F_{N_S(R)}(\ca L_R)$. Evidently 
$\ca F_R$ is a fusion subsystem of $N_{\ca F}(R)$. 

If $R\in\D$ then $\ca L_R$ is a subgroup of $\ca L$, of characteristic $p$ since $\ca L$ is proper; 
and $N_{\ca F}(R)=\ca F_{N_S(R)}(\ca L_R)$ by 2.2. Thus there is nothing to show in this case, and so 
we may assume that $R\notin\D$. As $R\in\ca F^s$ there exists an $\ca F$-conjugate $R'$ of $R$ such 
that $O_p(N_{\ca F}(R'))$ is centric in $\ca F$. As $\ca F$ is $\D\cup R^{\ca F}$-inductive by 3.2(c), 
it follows from 1.14(a) that $Q:=O_p(N_{\ca F}(R))$ is centric in $\ca F$. If $R=Q$ then 
$R\in\ca F^{cr}$, contrary to $R\notin\D$. Thus $R\neq Q$ and so $Q\in\D$ by 3.3(1). As 
$\ca F_R$ is a fusion subsystem of $N_{\ca F}(R)$, on $N_S(R)$, we have $Q\norm\ca F_R$, and then 
$Q$ is contained in every member of $(\ca F_R)^{cr}$. Thus: 
\roster 

\item "{(*)}" $(\ca F_R)^{cr}\sub\D_R$. 

\endroster 
For $P\in\D_R$ write $\ca L_P$ for the subgroup $N_{\ca L}(P)$ of $\ca L$. Then $R\leq\ca L_P$, and  
and $N_{\ca L_P}(R)=N_{\ca L_R}(P)$. As $\ca L$ is proper, $\ca L_P$ is of characteristic $p$, and  
$N_{\ca L_R}(P)$ is then of characteristic $p$ by II.2.7(b). This result, along with (*), shows that 
$\ca L_R$ is a proper locality on $\ca F_R$. Then $Q\norm\ca L_R$ by 2.3, and thus 
$\ca L_R=N_{\ca L_Q}(R)$ is a subgroup of $\ca L$ of characteristic $p$. 

It now remains only to show that $\ca F_R=N_{\ca F}(R)$, in order to complete the proof. Observe that 
$N_{\ca F}(R)$ is a fusion subsystem of $N_{\ca F}(Q)$, and that $N_{\ca F}(Q)$ is the fusion system 
$\ca F_{N_S(Q)}(\ca L_Q)$ of a finite group, by 2.2(a). Then each $N_{\ca F}(R)$-isomorphism is a 
conjugation map by an element of $N_{\ca L_Q}(R)$, and so $N_{\ca F}(R)$ is a fusion subsystem of 
$\ca F_R$. That $\ca F_R$ is a subsystem of $N_{\ca F}(R)$ has already been noted, so the required 
equality of fusion systems obtains. 
\qed 
\enddemo

\proclaim {5.2 (Theorem A1)} Let $(\ca L,\D,S)$ be a proper locality on $\ca F$ and let $\D^+$ be an 
$\ca F$-closed collection of subgroups of $S$ such that $\D\sub\D^+\sub\ca F^s$. 
\roster 

\item "{(a)}" There exists a proper locality $(\ca L^+,\D^+,S)$ on $\ca F$ such that $\ca L$ is the 
restriction $\ca L^+\mid_\D$ of $\ca L^+$ to $\D$. Moreover, $\ca L^+$ is generated by $\ca L$ as a 
partial group. 

\item "{(b)}" For any proper locality $(\w{\ca L},\D^+,S)$ on $\ca F$ whose restriction to $\D$ is $\ca L$, 
there is a unique isomorphism $\ca L^+\to\w{\ca L}$ which restricts to the identity map on $\ca L$. 

\endroster 
\endproclaim 

\demo {Proof} Suppose false, and among all counterexamples choose $(\ca L,\D,S)$ and $\D^+$ so that 
the set $\ca U=\D^+-\D$ has the smallest possible cardinality. 
We may choose $R\in\ca U$ so that the set $\D_1=\D\bigcup R^{\ca F}$ is $\ca F$-closed, and by 
1.15 $R$ may be chosen so that both $R$ and $O_p(N_{\ca F}(R))$ are fully normalized in 
$\ca F$. Thus the conditions (1) and (2) in Hypothesis 3.3 hold, with $\D_1$ in the role of $\D^+$. 
Then 3.3 holds in its entirety, by 5.1.  
Theorem 3.4 then yields a proper locality $(\ca L_1,\D_1,S)$ on $\ca F$, such 
that $\ca L_1\mid_\D=\ca L$, and such that $N_{\ca L_1}(R)=N_{\ca L}(R)$. The minimality of $|\ca U|$ 
yields the existence of a proper locality $(\ca L^+,\D^+,S)$ such that $\ca L^+\mid_{\D_1}=\ca L_1$ and 
with $\ca F_S(\ca L^+)=\ca F$. Then $\ca L^+\mid_\D=\ca L$. The explicit construction of $\ca L^+$ in 
section 3 shows that every element of $\ca L^+$ is a $\Pi^+$-product of elements of $\ca L$, so we have (a). 
Point (b) is then immediate from the corresponding uniqueness result for $\ca L_1$ (given by Theorem 3.4) 
and for $\ca L^+$ with respect to $\ca L_1$. 
\qed 
\enddemo

\proclaim {5.3 (Theorem A2)} Let the hypothesis and notation be as in Theorem A1. Let $\frak N$ be the 
set of partial normal subgroups of $\ca L$, and let $\frak N^+$ be the set of partial normal subgroups 
of $\ca L^+$. For each $\ca N\in\frak N$ let $\ca N^{\ca L^+}$ be the set of all elements of $\ca L^+$ of 
the form $\Pi^+(g\i,f,g)$, such that $f\in\ca N$ and such that $(g\i,f,g)\in\bold D^+$. 
Then there is an inclusion-preserving bijection 
$$ 
\frak N^+\to\frak N \ \ (\ca N^+\maps\ca N^+\cap\ca L).  
$$ 
In particular, $\ca N^+$ is the unique partial normal subgroup of $\ca L^+$ whose intersection with 
$\ca L$ is equal to $\ca N$, and $S\cap\ca N^+=S\cap\ca N$ for each $\ca N\in\frak N$. 
\endproclaim 

\demo {Proof} As in the proof of Theorem A1, assume that $\ca L$ is a counter-example with 
$|\D^+-\D|$ as small as possible. Define $R$, $\D_1$, and $\ca L_1$ as in the preceding proof. Let 
$\ca N\norm\ca L$ be a partial normal subgroup, let $\ca N_1$ be the partial subgroup of $\ca L_1$ 
generated by the set $\ca N^{\ca L_1}$ of $\ca L_1$-conjugates of elements of $\ca N$. Theorem 4.1 
then shows $\ca N_1\norm\ca L_1$ and $\ca N=\ca L\cap\ca N_1$. The minimality of $|\D^+-\D|$ yields 
the existence of a unique partial normal subgroup $\ca N^+\norm\ca L^+$ with $\ca L_1\cap\ca N^+=\ca N_1$, 
and then $\ca L\cap\ca N^+=\ca L\cap\ca N_1=\ca N$. Thus, the map 
$$ 
\mu:\frak N^+\to\frak N \ \ (\ca N^+\maps\ca N^+\cap\ca L)  
$$
is surjective. 

Now let $\w{\ca N}\norm\ca L^+$ with $\ca L\cap\w{\ca N}=\ca N$. Set $\w{\ca N}_1=\ca L_1\cap\w{\ca N}$. 
Then $\w{\ca N}_1=\ca N_1$ by 4.1, and then $\w{\ca N}=\ca N^+$ by the minimality of $|\D^+-\D|$. 
This shows that $\mu$ is injective, and so it only remains to show that $\mu$ is inclusion-preserving. 

Let $\ca M\norm\ca L$ with $\ca M\leq\ca N$, and set $\ca M_1=\<\ca M^{\ca L_1}\>$. Then 
$\ca M_1\leq\<\ca N^{\ca L_1}\>=\ca N_1$. The minimality of $|\D^+-\D|$ then gives 
$\ca M^+\leq\ca N^+$, completing the proof. 
\qed 
\enddemo

Here is an application. Note that by [He1] the product of partial normal subgroups of a locality is 
again a partial normal subgroup.

\proclaim {Corollary 5.4} Let $\ca M$ and $\ca N$ be partial normal subgroups of the proper locality 
$(\ca L,\D,S)$, and let $(\ca L^+,\D^+,S)$ be a proper expansion of $\ca L$. Then 
$(\ca M\ca N)^+=\ca M^+\ca N^+$, and if $\ca M\cap\ca N\leq S$ then 
$(\ca M\cap\ca N)^+=\ca M\cap\ca N$. 
\endproclaim 
 
\demo {Proof} Both $\ca M^+\ca N^+$ and $(\ca MN)^+$ are partial normal subgroups of $\ca L^+$, and 
$(\ca M\ca N)^+$ contains both $\ca M^+$ and $\ca N^+$, so $\ca M^+\ca N^+\leq (\ca M\ca N)^+$. Then 
$$ 
\ca L\cap\ca M^+\ca N^+\leq\ca L\cap(\ca M\ca N)^+=\ca M\ca N.  
$$ 
As $\ca M\leq\ca M^+$ and $\ca N\leq\ca N^+$ we obtain 
$$ 
\ca L\cap(\ca M\ca N)^+=\ca L\cap \ca M^+\ca N^+ , 
$$ 
and hence $(\ca M\ca N)^+=\ca M^+\ca N^+$. 

Set $R=\ca M\cap\ca N$, and assume that $R\leq S$. As $R\norm\ca L$ we then have $R\norm\ca L^+$ by 
2.3. We have 
$$ 
(\ca M^+\cap\ca N^+)\cap\ca L=(\ca M^+\cap\ca L)\cap(\ca N^+\cap\ca L)=\ca M\cap\ca N=R, 
$$ 
and also $R\cap\ca L=R$, so $R=\ca M^+\cap\ca N^+$ by Theorem A2.  
\qed 
\enddemo

We wish also to obtain a version of Theorem A for localities which are homomorphic images of proper 
localities. Thus, for the remainder of this section $(\ca L,\D,S)$ is a locality (not necessarily proper) 
on $\ca F$, and $\ca N\norm\ca L$ is a partial normal subgroup of $\ca L$. Set $T=S\cap\ca N$. 

For any $g\in\ca L$, $\ca Ng$ denotes the set of all products $xg$ such that $x\in\ca N$ and 
$(x,g)\in\bold D$, and we say that $\ca Ng$ is a right coset of $\ca N$ in $\ca L$. The analogous notion 
of left coset is obvious. The set of all cosets (left or right) of $\ca N$ in $\ca L$ is partially ordered 
by inclusion, and one thus has the notion of the maximal cosets of $\ca N$ in $\ca L$.  

By I.3.14 the maximal left cosets of $\ca N$ are the maximal right cosets, and these maximal cosets form a 
partition $\bar{\ca L}=\ca L/\ca N$ of $\ca L$. Let $\r:\ca L\to\bar{\ca L}$ be the map which sends 
$g\in\ca L$ to the unique maximal coset of $\ca N$ 
containing $g$. By I.3.16 there is a unique partial group structure on $\bar{\ca L}$ which makes 
$\r$ into a homomorphism of partial groups; and then the induced map 
$$ 
\r^*:\bold D(\ca L)\to\bold D(\bar{\ca L})
$$ 
is surjective. A homomorphism from a locality to a partial group is by definition a projection 
(cf. I.4.4) if the induced map between the domains of their products is surjective. Thus, $\r$ is a 
projection, and it is shown in I.4.3 that $\bar{\ca L}$ thereby has the structure of a locality  
$$
(\bar{\ca L},\bar\D,\bar S), 
$$ 
where $\bar D=\{P\r\mid P\in\D\}$, and where $\bar S=S\r\cong S/T$. An important (and elementary) property 
of any homomorphism of partial groups (I.1.15) is that it sends subgroups to subgroups. 

\vskip .1in 
The following lemma clarifies the relationship between the fusion system $\ca F$ of $\ca L$ and the 
fusion system of a homomorphic image of $\ca L$. 

\proclaim {Lemma 5.5} Let $(\ca L,\D,S)$ be a locality on $\ca F$, and let $\ca N\norm\ca L$ be a partial 
normal subgroup.  Let $(\bar{\ca L},\bar\D,\bar S)$ be the corresponding quotient 
locality, and let $\r:\ca L\to\bar{\ca L}$ be the canonical projection. Set $T=S\cap\ca N$, and
set $\bar{\ca F}:=\ca F_{\bar S}(\bar{\ca L})$.
\roster 

\item "{(a)}" $\r$ restricts to a surjective homomorphism $\s:S\to\bar S$, and $\s$ is fusion-preserving 
relative to $\ca F$ and $\bar{\ca F}$. That is, $\s$ is a homomorphism $\ca F\to\bar{\ca F}$ of fusion 
systems. Moreover, if $X$ and $Y$ are subgroups of $S$ containing $T$ then the induced map 
$$ 
\s_{X,Y}:Hom_{\ca F}(X,Y)\to Hom_{\bar{\ca F}}(\bar X,\bar Y)
$$ 
is surjective. 

\item "{(b)}" Let $X\leq S$ be a subgroup of $S$ containing $T$. Then $X$ is fully normalized in 
$\ca F$ if and only if $\bar X$ is fully normalized in $\bar{\ca F}$. 

\item "{(c)}" Let $\bar X\in\bar{\ca F}^c$, and let $X$ be the preimage of $\bar X$ in $S$. Then 
$X\in\ca F^c$. 

\item "{(d)}" Let $\bar X\in\bar{\ca F}^{cr}$, and let $X$ be the preimage of $\bar X$ in $S$. Then 
$X\in\ca F^{cr}$. 

\item "{(e)}" If $\ca F^{cr}\sub\D$ then $\bar{\ca F}^{cr}\sub\bar\D$. 

\endroster 
\endproclaim 

\demo {Proof} Point (a) is igiven by 1.20. For any subgroup $H$ of $\ca L$ write $\bar H$ for the 
image of $H$ under $\r$. 
Point (b) is given by the observation that if $X$ is a subgroup of $S$ containing $T$ then 
$\bar{N_S(X)}=N_{\bar S}(\bar X)$.  

Now let $\bar X\in\bar{\ca F}^c$, let $X$ be the preimage of 
$\bar X$ in $S$, and let $Y$ be an $\ca F$-conjugate of $X$. Then $\bar Y$ is an $\bar{\ca F}$-conjugate 
of $\bar X$, so $C_{\bar S}(\bar Y)\leq\bar Y$, and hence $C_S(Y)\leq Y$. This proves (c). Now 
assume further that $\bar X\in\bar{\ca F}^{cr}$, and let $Y\in X^{\ca F}$ with $Y$ fully normalized in 
$\ca F$. Then $\bar Y$ is fully normalized in $\bar{\ca F}$ by (b). Set $Q=O_p(N_{\ca F}(Y))$. Then 
$\bar Q\leq O_p(N_{\bar{\ca F}}(\bar Y))$ by (a), so $\bar Q=\bar Y$, and hence $Q=Y$. This shows that 
$X\in\ca F^{cr}$ (as defined in 1.1), and proves (d). Point (e) is immediate from (d). 
\qed 
\enddemo

If $(\ca L,\D,S)$ is a proper locality on $\ca F$, then the quotient locality $\ca L/\ca N$ (see I.4.5) 
need not be proper. For example, let $(\ca N,\G,T)$ be a direct product of pair-wise isomorphic 
proper localities $(\ca N_i,\G_i,T_i)$ with $1\leq i\leq k$. This means (cf. I.1.17) that the 
underlying set of $\ca N$ is the direct product of the sets $\ca N_i$, $\bold D(\ca N)$ is the direct 
product of the domains $\bold D(\ca N_i)$ (with the obvious product and inversion), $\G$ is the direct 
product of the collections $\G_i$, and $T$ is the direct product of the groups $T_i$. Then, let 
$\ca L$ be the partial group obtained as the natural semi-direct product $\ca N\rtimes H$, where $H$ is 
the symmetric group on $k$ letters. (The reader should have no difficulty, at this stage, in defining 
the partial group $\ca L$.) Let $S$ be a maximal $p$-subgroup of $\ca L$ containing $T$, and let 
$\D$ be the set of all $P\leq S$ such that $P\cap T\in\G$. Then $(\ca L,\D,S)$ is a proper locality, 
while $\ca L/\ca N$ is isomorphic to the group $H$, which is a proper locality if and only if $k=1$, or 
$(k,p)$ is $(2,2)$, $(3,3)$, or $(4,2)$. Thus, Theorems A1 and A2 are not directly applicable to 
homomorphic images of proper localities.

\proclaim {Theorem 5.6} Let $(\ca L,\D,S)$ be a proper locality on $\ca F$, let $\ca N\norm\ca L$, let 
$\bar{\ca L}=\ca L/\ca N$ be the quotient locality, and let $\r:\ca L\to\bar{\ca L}$ be the canonical 
projection. Let $\D^+$ be an $\ca F$-closed set of subgroups of $S$ 
such that $\D\sub\D^+\sub\ca F^s$. Let $(\ca L^+,\D^+,S)$ be the expansion of $\ca L$ to $\D^+$, 
set $\bar{\D}^+=\{P\r\mid P\in\D^+\}$, and set $\bar{\ca F}=\ca F_{\bar S}(\bar{\ca L})$.   
\roster 

\item "{(a)}" There is a locality $(\bar{\ca L}^+,\bar\D^+,\bar S)$ on $\bar{\ca F}$ whose restriction to 
$\bar\D$ is $\bar{\ca L}$, and $\bar{\ca L}^+$ is unique up to a unique isomorphism which restricts to 
the identity map on $\bar{\ca L}$. 

\item "{(b)}" There is a unique projection $\r^+:\ca L^+\to\bar{\ca L}^+$ whose restriction to $\ca L$ 
is $\r$. Moreover, $Ker(\r^+)=\ca N^+$, where $\ca N^+$ is the partial normal subgroup of $\ca L^+$ 
which intersects $\ca L$ in $\ca N$. 

\endroster 
\endproclaim 

\demo {Proof} Note first of all that since $\ca F^{cr}\sub\D$ we have $\bar{\ca F}^{cr}\sub\bar\D$ by 
5.5(e). Notice also that if $\D=\D^+$ then there is nothing to prove, so we may assume that $\D$ is 
properly contained in $\D^+$.  

Among all $\ca F$-closed sets $\D_1$ with $\D\sub\D_1\sub\D^+$, let $\D_1$ be maximal subject to the 
condition that (a) and (b) hold with $\D_1$ in the role of $\D^+$. Then $\D_1$ is a proper subset of 
$\D^+$, and by replacing $\D$ with $\D_1$ we reduce (as in the proof of Theorem A1) to the case 
where $\D^+=\D\cup R^{\ca F}$ for some $R\leq S$. 

Take $R$ to be fully normalized in $\ca F$, and suppose that $T\nleq R$. Then $RT\in\D$, and then 
$\D^+\r=\D\r=\bar\D$. Then $(\bar{\ca L},\bar\D,\bar S)$ is the unique locality on $\bar{\ca F}$ whose 
restriction to $\bar\D$ (namely, itself) is $\bar{\ca L}$, and thus (a) holds in this case. In order 
to verify (b) in this case recall that, by Theorem A1, $\ca L^+$ is the ``free amalgamated product" in the 
category of partial groups of the ``amalgam" given by the inclusion maps $\ca L_0\to\ca L$ and 
$\ca L_0\to\ca L_0^+$. Point (b) will follow if: 
\roster 

\item "{(1)}" There is a unique homomorphism $\r_0^+:\ca L_0^+\to\bar{\ca L}$ such that $\r_0^+$  
agrees with $\r$ on $\ca L_0$. 

\item "{(2)}" $\r_0+$ induces a surjection $\bold D(\ca L_0^+)\to\bold D(\bar{\ca L}_0)$, where 
$\bar{\ca L}_0$ is the set of all $\bar f\in\bar{\ca L}$ such that $\bar S_{\bar f}$ contains a  
$\bar{\ca F}$-conjugate of $\bar R$. 

\item "{(3)}" $Ker(\r^+)=\ca N^+$. 

\endroster 
Versions of the same three points will need to be verified in the case where $T\leq R$, and our 
approach will be to merely sketch the proofs in each case, leaving some of the entirely mechanical details 
to the reader. Thus, let $\Phi$ be the set of triples $(x\i,g,y)$ defined following 3.4, 
let $\approx$ be the relation on $\Phi$ defined in 3.10, and let 
$\r^*:\bold W(\ca L)\to\bold W(\bar{\ca L})$ be the $\r$-induced homomorphism of free monoids. Then 
$\r^*$ maps $\Phi$ into $\bold D(\bar{\ca L})$, and it may be routinely verified that the composition 
$$ 
\Phi@>\r^*>>\bold D(\bar{\ca L})@>\bar\Pi>>\bar{\ca L} 
$$
is constant on $\approx$-classes. Thus, there is a mapping $\r_0^+:\ca L_0^+\to\bar{\ca L}$ given by 
$[x\i,g,y]\maps\bar\Pi(\bar x\i,\bar g,\bar y)$. The verification that $\r_0^+$ is a homomorphism is also 
routine, and yields (1). Since $\r^*$ maps $\bold D(\ca L)$ onto $\bold D(\bar{\ca L})$, (2) is then 
immediate. As $Ker(\r^+)\cap\ca L=\ca N$, (3) is immediate from Theorem A2. Thus, the theorem holds if 
$T\nleq R$.  

Assume henceforth that $T\leq R$. Then $N_{\bar{\ca L}}(\bar R)=\bar{N_{\ca L}(R))}$, since each element of 
$\bar{\ca L}$ has a $\r$-preimage in $N_{\ca L}(T)$ by I.3.11. Then $N_{\bar{\ca L}}(\bar R)$ is a 
subgroup of $\bar{\ca L}$ by I.1.15. Thus Hypothesis 3.3 is satisfied, with $\bar{\ca L}$ and $\bar R$ in 
the roles of $\ca L$ and $R$. Theorem 3.4 then yields (a). For the same reasons as in the preceding 
case (and because the verification of (3) did not in fact make use of the hypothesis that 
$T\nleq R$) it now suffices to verify: 
\roster 

\item "{(4)}" There is a unique homomorphism $\r_0^+:\ca L_0^+\to\bar{\ca L}^+$ such that $\r_0^+$  
agrees with $\r$ on $\ca L_0$. 

\item "{(5)}" $\r_0+$ induces a surjection $\bold D(\ca L_0^+)\to\bold D(\bar{\ca L}_0^+)$, where 
$\bar{\ca L}_0^+$ is the set of all $\bar f\in\bar{\ca L}^+$ such that $\bar S_{\bar f}$ contains a  
$\bar{\ca F}$-conjugate of $\bar R$.  

\endroster 
Let $\Phi$ and $\r^*$ be as above, and define $\bar\Phi$ to be the set of triples 
$\bar\phi=(\bar x\i,\bar g,\bar y)$ such that one has 
$$ 
\bar U@>\bar x\i>>\bar R@>\bar g>>\bar R@>\bar y>>\bar V 
$$ 
(a sequence of conjugation isomorphisms between subgroups of $\bar S$, labeled by the conjugating elements), 
with $N_{\bar S}(\bar U)\leq\bar S_{\bar x\i}$, and with $N_{\bar S}(\bar V)\leq\bar S_{\bar y\i}$. Then 
$\r^*$ maps $\Phi$ onto $\bar\Phi$ by 5.5(a). It need cause no confusion to denote also by $\sim$ and 
$\approx$ the two equivalence relations on $\bar\Phi$ given by direct analogy with 3.5 and 3.10. Again by 
means of 5.5(a), one verifies that the restriction of $\r^*$ to $\Phi$ preserves these equivalence relations, 
and hence induces a surjective mapping $\r_0^+\ca L_0^+\to\bar{\ca L}_0^+$. Here $\ca L_0$ is by definition 
the set of elements $f\in\ca L$ such that $S_f$ contains an $\ca F$-conjugate of $R$, and any such $f$ 
is identified with its $\approx$-class $[f]$, consisting of all those $\phi\in\Phi\cap\bold D(\ca L)$ 
such that $\Pi(\phi)=f$. The analogous definition of $\bar{\ca L}_0$ leads to the conclusion that 
the restriction of $\r_0^+$ to $\ca L_0$ is the (surjective) homorphism $\r_0:\ca L_0\to\bar{\ca L}_0$ 
induced by $\r$. The verification that $\r_0^+$ is a homomorphism, and hence that (4) holds, is then 
straightforward. 

Set $\bold D_0^+=\bold D(\ca L_0^+)$ and $\bar{\bold D}_0^+=\bold D(\bar{\ca L_0}^+)$, let  
$\bar w\in\bar{\bold D}_0^+$, and write $\bar w=([\bar\phi_1],\cdots,[\bar\phi_n])$. Here 
$[\bar\phi_i]$ is the $\approx$-class of an element $\bar\phi_i=(\bar x_i\i,\bar g_i,\bar y_i)$ in 
$\bar\Phi$, and the representatives $\bar\phi_i$ may be chosen so that the sequence 
$\bar\g=(\bar\phi_1,\cdots,\bar\phi_n)$ is a $\bar\G$-form of $\bar w$, as defined in 3.12. That is, 
the word 
$$ 
\bar w_{\bar g}=\bar\phi_1\circ\cdots\circ\bar\phi_n 
$$ 
has the property that $\bar S_{\bar w_{\bar\g}}$ contains a $\bar{\ca F}$-conjugate of $\bar R$. Let 
$\phi_i$ be a $\r^*$-preimage of $\bar\phi_i$ in $\Phi$, set $\g=(\phi_1,\cdots,\phi_n)$, let 
$[\phi_i]$ be the $\approx$-class of $\phi_i$ (in $\ca L_0^+$), and set  
$w=([\phi_1],\cdots,[\phi_n])$. One verifies that $\g$ is a $\G$-form of $w$, and hence that 
$w\in\bold D_0^+$. Since $\r^*$ maps $w_\g$ to $\bar w_{\bar\g}$, it follows that $\r_0^+$ induces a 
surjection as required in (5). Thus (5) holds, and the proof is complete. 
\qed 
\enddemo

The preceding result is essentially a generalization of Theorem A1 to homomorphic images of proper 
localities. Here is the corresponding version of Theorem A2.

\proclaim {Theorem 5.7} Let the hypotheses and the setup be as in the preceding theorem, let 
$\frak N$ be the set of partial normal subgroups of $\ca L$ containing $Ker(\r)$, and let 
$\frak N^+$ be the set of partial normal subgroups of $\ca L^+$ containing $Ker(\r^+)$. Also, let
$\bar{\frak N}$ be the set of partial normal subgroups of $\bar{\ca L}$, and let $\bar{\frak N}^+$ be the 
set of partial normal subgroups of $\bar{\ca L}^+$. Then there is a commutative diagram of 
inclusion-preserving bijections:  
$$ 
\CD 
\frak N       @>\eta>>      \frak N^+ \\ 
@V{\b}VV                      @VV{\b^+}V  \\ 
\bar{\frak N} @>>\bar\eta> \bar{\frak N}^+ 
\endCD 
$$ 
where $\eta$ is the map $\ca N\maps\ca N^+$, $\b$ is the map $\ca N\maps\ca N\r$, and $\b^+$ is the 
map $\ca N^+\maps\ca N^+\r^+$. Moreover:  
\roster 

\item "{(a)}" $Ker(\r^+)=Ker(\r)\eta$, and 

\item "{(b)}" for each $\ca N\in\frak N$, $\bar{\ca N}^+$ is the unique 
partial normal subgroup of $\bar{\ca L}^+$ whose intersection with $\bar{\ca L}$ is $\bar{\ca N}$. 

\endroster 
\endproclaim 

\demo {Proof} By I.4.7 there are inclusion-preserving bijections 
$$ 
\b:\frak N\to\bar{\frak N}\ \ \text{and}\ \ \b^+:\frak N^+\to\bar{\frak N}^+, 
$$ 
by which a partial normal subgroup is sent to its image in $\bar{\ca L}$ (or in $\bar{\ca L}^+$) under 
the canonical projection $\r$ (or $\r^+$). Theorem A2 yields the inclusion-preserving bijection 
$\eta:\frak N\to\frak N^+$ which sends $\ca N\in\frak N$ to the unique $\ca N^+\norm\ca L^+$ 
such that $\ca N^+$ intersects $\ca L$ in $\ca N$. We obtain the above commutative square by 
taking $\bar\eta$ to be $\b\i\circ\eta\circ\b^+$, so it remains to prove (a) and (b), 

We have $Ker(\r)=\ca L\cap Ker(\r^+)$ since $\r$ is the restriction of $\r^+$ to $\ca L$. 
Thus $Ker(\r^+)$ is a partial normal subgroup of $\ca L^+$ which intersects $\ca L$ in $Ker(\r)$. As 
this condition characterized $Ker(\r)\eta$, we have (a). 
 
Fix $\ca N\in\frak N$, set $\bar{\ca N}=\ca N\r$, and set $\bar{\ca N}^+=\ca N^+\r^+$. Thus 
$\bar{\ca N}^+=\bar{\ca N}\bar\eta$. Identify $\bar{\ca L}^+$ with $\ca L^+/Ker(\r^+)$, let 
$g\in\ca N^+$, and let $[g]$ be the maximal coset of $Ker(\r^+)$ in $\ca L^+$ containing $g$. 
Then $[g]=Ker(\r^+)h$ for some $h\in[g]$ such that $h$ is $\up$-maximalin $\ca L^+$ with respect to 
$Ker(\r^+)$. The Splitting Lemma (I.3.12) yields $x\in Ker(\r^+)$ such that $(x,h)\in\bold D(\ca L^+)$, 
$xh=g$, and $S_{(x,h)}=S_g$. Then $(x\i,x,h)\in\bold D(\ca L^+)$, and $\bold D$-associativity yields 
$h=x\i g$. As $Ker(\r)\leq\ca N$, (a) implies that $Ker(\r^+)\leq\ca N^+$. Thus $h\in\ca N^+$, and 
then $[g]\sub\ca N^+$. Now suppose that also $[g]\in\bar{\ca L}$. As $\bar{\ca L}$ is the image of 
$\ca L$ under $\r^+$, we may then assume that $g$ was chosen so that $g\in\ca L$. That is, we may 
take $g\in\ca L\cap N^+$. Then $g\in\ca N$, and this shows that $\bar{\ca N}^+$ intersects 
$\bar{\ca L}$ in $\bar{\ca N}$. The uniqueness of $\bar{\ca N}^+$ for this property is given by the 
commutativity of the displayed diagram, and so (b) holds. 
\qed 
\enddemo

\proclaim {Proposition 5.8} Let $(\ca L,\D,S)$ be a proper locality on $\ca F$, let $(\ca L^+,\D^+,S)$ be a 
proper expansion of $\ca L$ on $\ca F$, and let $\L$ be a subgroup of $Aut(\ca L)$ such that both $\D$ and 
$\D^+$ are $\L$-invariant. Then each $\l\in\L$ extends in a unique way to an automorphism of $\ca L^+$. 
\endproclaim 

\demo {Proof} It suffices to consider the case where $\D^+=\D\cup\ca R$, where $\ca R$ is a set of 
subgroups of $S$ which is minimal subject to being both $\ca F$-invariant and $\L$-invariant. Thus 
$\ca R$ is the disjoint union of sets $(R_i)^{\ca F}$ $(1\leq i\leq m)$, where $\{R_1,\cdots,R_m\}$ is 
a $\L$-orbit of subgroups of $S$. By 3.20(c) $\ca L^+$ is a colimit in the category of 
partial groups (cf. I.1.17):
$$ 
\ca L^+=\ca L_1^+\underset{\ca L}\to\coprod\cdots\underset{\ca L}\to\coprod\ca L_m^+, 
$$ 
where each $\ca L_i^+$ is an elementary expansion of $\ca L$ whose set $\D_i^+$ of objects is 
$\D\cup(R_i)^{\ca F}$. 

We may assume that $R_1$ has been chosen so that both $R_1$ and $O_p(N_{\ca F}(R_1))$ are fully normalized 
in $\ca F$. Plainly $\L$ induces a group of automorphisms of $\ca F$, so the preceding condition applies 
also to each $R_i$. We may then form the set $\Phi_i$ if all triples $(x\i,g,y)$ as in 3.9, with 
respect to $R_i$ in the role of $R$. Then $\ca L_i^+$ is itself a colimit of the form 
$$ 
\ca L_i^+=\ca L\underset{\ca M}\to\coprod\ca M^+, 
$$ 
where $\ca M$ is the partial subgroup of $\ca L$ consisting of those elements $h\in\ca L$ that can be 
written as a product $\Pi(x\i,g,y)$ for some $(x\i,g,y)\in\Phi_i$, and where $\ca M^+$ is the partial 
group consisting of all of the equivalence classes $[x\i,g,y]$ obtained from $\Phi_i$. 

Let $\l\in\L$ and let an index $i$ be given. Then $R_i\l=R_j$ for some $j$. In order 
to extend $\l$ to a mapping $\ca L^+\to\ca L^+$ it now suffices to show that $\l$ induces a mapping 
$\ca M_i^+\to\ca M_j^+$ via 
$$ 
[x\i,g,y]\maps [(x\i)\l,g\l,y\l]. \tag*   
$$ 
For this it suffices to verify that $\l$ is equivariant with respect to the equivalence relations  
$\approx_i$ on $\Phi_i$ and $\approx_j$ on $\Phi_j$ given by 3.10 (with the obvious adjustment 
needed to accomodate the indices $i$ and $j$). This verification - that (*) is well defined - is 
entirely straightforward, and is left to the reader. That the mapping (*) is a homomorphism 
$\ca M_i^+\to\ca M_j^+$ is then immediate from the definition of the product in an expansion (see 3.13 
and 3.14). One then observes that $\l\i$ induces a homomorphism $\ca M_j^+\to\ca M_i^+$ which is the 
inverse of the homomorphism induced by $\l$. The functoriality of the colimit of partial groups then 
yields the desired extension of $\l$ to an automorphism $\l^+$ of $\ca L^+$. The uniqueness of $\l^+$ 
follows from the observation that $[x\i,g,y]$ can be written as the product $x\i gy$ in $\ca M_i^+$. 
Thus any extension of $\l$ to an endomorphism of $\ca L^+$ must necessarily satisfy (*) and will thereby be 
equal to $\l^+$. 
\qed 
\enddemo

\vskip .2in 
\noindent 
{\bf Section 6: Saturated fusion systems} 
\vskip .1in 

Throughout this section $\ca F$ is a fusion system on a finite $p$-group $S$. 
\vskip .1in 

We have been able to proceed, so far, with only very naive notions pertaining to fusion systems 
in general. The 
most useful of these has been the notion of inductivity (and $\G$-inductivity) introduced in 
1.11. This notion cannot be sufficient for our purposes, for two reasons. First, we have no criteria 
for deciding if inductivity follows from $\G$-inductivity for suitable $\ca F$-closed sets $\G$; and 
second, inductivity of $\ca F$ is not necessarily 
inherited by subsystems $N_{\ca F}(X)$, for $X$ fully normalized in $\ca F$. What is needed is a 
condition on $\ca F$ which implies inductivity and which resolves these two issues. 
This is provided by {\it saturation}, whose definition we take from [AKO], and which we now review. 

First: a subgroup $Q$ of $S$ is defined to be {\it fully automized} in $\ca F$ if $Aut_S(Q)$ is a 
Sylow $p$-subgroup of $Aut_{\ca F}(Q)$. Second: a subgroup $Q$ of $S$ is defined to be {\it receptive} 
in $\ca F$ if every $\ca F$-isomorphism $\phi:P\to Q$ extends to an $\ca F$-homomorphism 
$\bar\phi:N_\phi\to S$, where 
$$ 
N_\phi=\{x\in N_S(P)\mid \phi\i\circ c_x\circ\phi\in Aut_S(Q)\}. 
$$

\definition {Definition 6.1} Let $\G$ be an $\ca F$-closed set of subgroups of $S$. Then $\ca F$ is 
{\it $\G$-saturated} if for every $P\in\G$ there exists an $\ca F$-conjugate $Q$ of $P$ such that 
$Q$ is fully automized and receptive in $\ca F$. We say $\ca F$ is {\it saturated} if $\ca F$ is 
$\G$-saturated for the set $\G$ of all subgroups of $S$. 
\enddefinition 

\proclaim {Lemma 6.2} Let $\ca F$ be a fusion system on the finite $p$-group $S$, let $\G$ be an 
$\ca F$-closed set of subgroups of $S$, and assume that $\ca F$ is $\G$-saturated. 
\roster 

\item "{(a)}" A subgroup $Q\in\G$ is fully normalized in $\ca F$ if and only if $Q$ is fully automized 
and receptive in $\ca F$. 

\item "{(b)}" $\ca F$ is $\G$-inductive. 

\endroster 
\endproclaim 

\demo {Proof} Let $P\in\G$, let $Q$ be an $\ca F$-conjugate of $P$ such that $Q$ is fully automized and 
receptive in $\ca F$, and let $\phi:P\to Q$ be an $\ca F$-isomorphism. 
Then conjugation by $\phi$ is an isomorphism 
$$ 
c_\phi:Aut_{\ca F}(P)\to Aut_{\ca F}(Q). 
$$ 
Set $X=Aut_S(P)$ and $Y=Aut_S(Q)$. Then $Y\in Syl_p(Aut_{\ca F}(Q))$, so there exists $\a\in Aut_{\ca F}(Q)$ 
such that the inner automorphism $c_\b$ of $Aut_{\ca F}(Q)$ carries $(X)c_{\phi}$ into $Y$. Set 
$\psi=\phi\circ\b$. Then $c_\psi$ carries $X$ into $Y$, and thus $N_S(P)\leq N_\psi$. 
There is then an extension of $\psi$ to an $\ca F$-homomorphism $\bar\psi:N_S(P)\to N_S(Q)$. In 
particular, this shows that (a) implies (b), and so we need only prove (a). 

Suppose now that, in the preceding argument, $P$ is fully normalized in $\ca F$. Then $\bar\psi$ is an 
isomorphism, so $\bar\psi$ restricts to an isomorphism $C_S(P)\to C_S(Q)$, and thus $c_{\bar\psi}$ is an 
isomorphism $Aut_S(P)\to Aut_S(Q)$. Thus $P$ is fully automized (and $Q$ is fully normalized) in $\ca F$. 

Now let $R\in P^{\ca F}$, let $\g:R\to P$ be an $\ca F$-isomorphism, and set $\eta=\g\circ\psi$. Then 
$N_{\g}\leq N_\eta$, and $\eta$ extends to an $\ca F$-homomorphism $\w\eta:N_{\g}\to N_S(Q)$. Set  
$D=(N_{\g})\w\eta$, and let $\b$ be the restriction of ${\bar\psi}\i$ to $D$. Then $\w\eta\circ\b$ is 
an extension of $\g$ to $N_{\g}$, and thus $P$ is receptive. Thus (a) holds. 
\qed 
\enddemo

Recall from the discussion preceding 1.11 that if $\G$ is a non-empty set of subgroups of $S$ then $\ca F$ is 
defined to be $\G$-generated if every $\ca F$-isomorphism is a composite of restrictions of 
$\ca F$-isomorphisms $\phi:X\to Y$ with $X,Y\in\G$. 

\proclaim {Lemma 6.3} Let $(\ca L,\D,S)$ be a locality on $\ca F$. Then $\ca F$ is $\D$-saturated and 
$\D$-generated. 
\endproclaim 

\demo {Proof} By definition, $\ca F$ is generated by the set of conjugation maps $c_g:S_g\to S_{g\i}$. 
Since $S_g\in\D$ for all $g\in\ca L$, $\ca F$ is thus $\D$-generated. 

Let $Q\in\D$ with $Q$ fully normalized in $\ca F$. Then $N_S(Q)\in Syl_p(N_{\ca L}(Q))$ by I.2.10. Since 
$Aut_{\ca F}(Q)$ is a homorphic image of $N_{\ca L}(Q)$, $Q$ is then fully automized in $\ca F$. Now let 
$P\in Q^{\ca F}$ and let $\phi:P\to Q$ be an $\ca F$-isomorphism. Then $\phi=c_g$ for some $g\in\ca L$, 
and $c_g$ is an isomorphism $N_{\ca L}(P)\to N_{\ca L}(Q)$ by I.2.3(b). Here $N_\phi$ is simply 
$N_{S_g}(P)$, so $Q$ is receptive, and thus $\ca F$ is $\D$-saturated.  
\qed 
\enddemo

The following result, due to Broto, Castellana, Grodal, Levi, and Oliver [Theorem 2.2 in 5A1], 
is an essential tool for determining whether a fusion system is saturated. Since our definition 
of ``centric radical" is not the standard one (cf. 1.8) we provide a few lines of proof to address that 
point.

\proclaim {Theorem 6.4} Let $\ca F$ be a fusion system on $S$, and define $\ca P(\ca F)$ to be the set of all 
subgroups $P\leq S$ such that $P$ is $\ca F$-centric and such that for each $Q\in P^{\ca F}$ with $Q$ 
fully normalized in $\ca F$ we have: 
$$ 
O_p(Aut_{\ca F}(Q))\cap Aut_S(Q)=Inn(Q). \tag*
$$ 
Then $\ca P(\ca F)\sub\ca F^{cr}$. Further, let $\G$ be an $\ca F$-closed set of subgroups of $S$,  
and assume: 
\roster 

\item "{(1)}" $\G$ is $\G$-saturated and $\G$-generated. 

\item "{(2)}" $\ca P(\ca F)\sub\G$. 

\endroster 
Then $\ca F$ is saturated, and $\ca P(\ca F)=\ca F^{cr}$.  
\endproclaim 

\demo {Proof} As mentioned, the only issue here is that of showing $\ca P(\ca F)=\ca F^{cr}$. 
Let $P\in\ca P(\ca F)$ be fully normalized in $\ca F$, and set $\w P=O_p(N_{\ca F}(P))$. By definition, 
every $N_{\ca F}(P)$-isomorphism extends to an $N_{\ca F}(P)$-isomorphism whose domain contains $\w P$, so 
$Aut_{\w P}(P)\leq O_p(Aut_{\ca F}(P))\cap Aut_S(P)$. The condition (*) then yields 
$Aut_{\w P}(P)\leq Inn(P)$. As $P$ is centric in $\ca F$ we then have $\w P=P$, and so $P\in\ca F^{cr}$ 
by definition 1.8. Thus $\ca P(\ca F)\sub\ca F^{cr}$, and it remains to show the 
reverse inclusion under the assumption (which we now make) that $\ca F$ is saturated.  

Let $Q\in\ca F^c$, with $Q$ fully normalized in $\ca F$. Then definition 6.1 implies that every 
$\ca F$-automorphism of $Q$ extends to an $\ca F$-automorphism of the preimage $\w Q$ in $N_S(Q)$ of 
$O_p(Aut_{\ca F}(Q))$. Thus $\w Q\leq O_p(N_{\ca F}(Q))$, and so $Q=\w Q$ if $Q\in\ca F^{cr}$. 
Thus $\ca P(\ca F)=\ca F^{cr}$ if $\ca F$ is saturated. 
\qed 
\enddemo

\proclaim {Corollary 6.5} Let $(\ca L,\D,S)$ be a locality on $\ca F$, and suppose that  
$\ca P(\ca F)\sub\D$ or that $\ca F^{cr}\sub\D$. Then $\ca F$ is saturated, and $\ca P(\ca F)=\ca F^{cr}$. 
In particular, the fusion system of a proper locality is saturated.  
\endproclaim 

\demo {Proof} As $\ca P(\ca F)\sub\ca F^{cr}$ by 6.3, we have $\ca P(\ca F)\sub\D$. Let $\G$ be the 
overgroup-closure of $\ca P(\ca F)$ in $S$. Then $\G$ is $\ca F$-closed, since $\ca P(\ca F)$ is 
$\ca F$-invariant. By 6.2, $\ca F$ is $\G$-saturated and $\G$-generated. Thus the hypotheses (1) and (2) 
of 6.4 are fulfilled, hence $\ca F$ is saturated and $\ca P(\ca F)=\ca F^{cr}$.   
\qed 
\enddemo 

The following result is an immediate consequence of the Alperin-Goldschmidt theorem for fusion 
systems (cf. I.1.14 in [AKO]). Recall from the discussion preceding 1.11 that $\ca F$ is defined to be 
$(cr)$-generated if every $\ca F$-isomorphism is a composite of restrictions of $\ca F$-automorphisms 
of members of $\ca F^{cr}$ which are fully normalized in $\ca F$.

\proclaim {Lemma 6.6} Let $\ca F$ be a saturated fusion system on $S$. Then $\ca F$ is $(cr)$-generated.  
Thus, if $\ca F_0$ is a subsystem of 
$\ca F$ and $Aut_{\ca F}(Q)=Aut_{\ca F_0}(Q)$ for each $Q\in\ca F^{cr}$ such that 
$Q$ is fully normalized in $\ca F$, then $\ca F_0=\ca F$. 
\qed 
\endproclaim  

Recall from 3.1 the notion of a localizable pair. 

\proclaim {Lemma 6.7} Let $(\ca L,\D,S)$ be a locality and let $(\ca H,\G)$ be a localizable pair in 
$\ca L$. Set $R=S\cap\ca H$, set $\ca E=\ca F_R(\ca H)$, and assume that $\ca E^{cr}\sub\G$. Then 
$\ca E$ is $(cr)$-generated and $\ca E=\ca F_R(\ca H_\G)$. Moreover, if $\ca L$ is proper then 
$(\ca H,\G,R)$ is proper if either: 
\roster 

\item "{(1)}" $\ca H\norm\ca L$, or 

\item "{(2)}" $\ca H=N_{\ca L}(X)$ for some subgroup $X\leq S$. 

\endroster 
\endproclaim 

\demo {Proof} Set $\ca E_0=\ca F_R(\ca H_\G)$. Then $\ca E_0$ is a fusion subsystem of $\ca E$. But 
also $\ca E=\<Aut_{\ca E}(P)\mid P\in\G\>$ is a subsystem of $\ca E_0$since $\ca E$ is $(cr)$-generated 
and since $\ca E^{cr}\sub\G$. Thus $\ca E=\ca E_0$, and $\ca H_\G$ is a locality on $\ca E$. Now assume 
that $\ca L$ is proper and let $P\in\G$. Then $N_{\ca H}(P)$ is a subgroup of $\ca L$, so there 
exists $Q\in\D$ with $N_{\ca H}(P)\leq N_{\ca L}(Q)$ by I.2.11(a). Set $G=N_{\ca L}(Q)$. 
We may assume that $P\leq Q$, and then $N_{\ca H}(P)=N_{G\cap\ca H}(P)$. If $\ca H\norm\ca L$ then 
$G\cap\ca H\norm G$ and then $N_{\ca H}(P)$ is of characteristic $p$ by 2.7(a). On the other hand, if 
$\ca H=N_{\ca L}(X)$ with $X\leq S$ we may also take $X\leq Q$, and then 
$N_{\ca H}(P)=N_{N_{G\cap\ca H}(X)}(P))$. In this case $N_{\ca H}(P)$ is of characteristic $p$ by 
2.7(b). Thus $\ca H_\G$ is a proper locality if $\ca L$ is proper and if (1) or (2) obtains. 
\qed 
\enddemo

The following result (and references to several proofs) may be found in [AKO, Theorem I.11].

\proclaim {Proposition 6.8} Let $\ca F$ be a saturated fusion system on the finite $p$-group $S$, 
let $X$ be fully normalized in $\ca F$, and let $Y$ be fully centralized in $\ca F$. Then 
$N_{\ca F}(X)$ and $C_{\ca F}(Y)$ are saturated. 
\qed 
\endproclaim

The next few results concern the set $\ca F^s$ of $\ca F$-subcentric subgroups of the fusion 
system $\ca F$ of a proper locality, and they are due to Ellen Henke [He2]. We provide 
proofs for the convenience of the reader.  
The key insight into $\ca F^s$ is given by [Lemma 3.1 in He2], which is as follows.

\proclaim {Lemma 6.9} Let $\ca F$ be the fusion system on $S$ of a proper locality, and let $V\leq S$ be 
fully centralized in $\ca F$. Then $V\in\ca F^s$ if and only if $O_p(C_{\ca F}(V))$ is centric in 
$C_{\ca F}(V)$. 
\endproclaim 

\demo {Proof} If $U\in V^{\ca F}$ is fully normalized in $\ca F$ then 1.14 yields an isomorphism  
$C_{\ca F}(U)\cong C_{\ca F}(V)$ of fusion systems.  We may therefore assume to begin with that $V$ is 
fully normalized in $\ca F$. Set $\ca F_V=N_{\ca F}(V)$ and $R=O_p(\ca F_V)$. Also, set 
$\ca C_V=C_{\ca F}(V)$ and $Q=O_p(\ca C_V)$. 

Suppose first that the lemma holds with $\ca F_V$ in place of $\ca F$. That is, suppose that $R$ is 
centric in $\ca F_V$ if and only if $Q$ is centric in $\ca C_V$. Since $R\in(\ca F_V)^c$ if and only if 
$R\in\ca F^c$ by 1.19, we then have the lemma in general. Thus, we are reduced to the case where 
$V\norm\ca F$, and $R=O_p(\ca F)$. 

Fix a proper locality $(\ca L,\D,S)$ on $\ca F$. We may assume that $\ca F^c\sub\D$ by Theorem A1.  
Suppose that $V\in\ca F^s$. Then $R\in\ca F^c$, so $R\in\D$, and $\ca L$ is the group $N_{\ca L}(R)$. 
Notice that $[R,C_S(VQ)]\leq C_R(VQ)\leq Q$, and that both 
$R$ and $Q$ are normal subgroups of $\ca L$. Then $C_S(VQ)\leq R$ by 2.7(c). Then also 
$$ 
C_{C_S(V)}(Q)=C_S(VQ)\leq C_R(V)\leq Q, 
$$ 
and thus $Q$ is centric in $C_{\ca F}(V)$. On the other hand, assuming that $Q$ is centric in 
$C_{\ca F}(V)$, we obtain 
$$ 
C_S(R)\leq C_S(VQ)\leq Q\leq R, 
$$ 
and so $R\in\ca F^c$, as required. 
\qed 
\enddemo 

\proclaim {Corollary 6.10} Let $\ca F$ be a fusion system of a proper locality. Then 
$$ 
\ca F^{cr}\sub\ca F^c\sub\ca F^q\sub\ca F^s.   
$$
\endproclaim 

\demo {Proof} We have $\ca F^{cr}\sub\ca F^c$ by definition. Let $(\ca L,\D,S)$ be a proper locality on 
$\ca F$, and $P\leq S$ be fully normalized in $\ca F$. 
If $P\in\ca F^c$ then $C_{\ca F}(P)$ is the trivial fusion system on $Z(P)$, and so $P\in\ca F^q$. Now 
suppose instead that we are given $P\in\ca F^q$. Then $O_p(C_{\ca F}(Q))=C_S(Q)$, and then $P\in\ca F^s$ 
by 6.9. 
\qed 
\enddemo

\proclaim {Corollary 6.11} Let $\ca F$ be the fusion system on $S$ of a proper locality, and let  
$V\leq S$ with $V$ fully normalized in $\ca F$. Then 
$$
\{Q\in N_{\ca F}(V)^s\mid V\leq Q\}\sub\ca F^s. 
$$
\endproclaim 

\demo {Proof} One need only observe that $C_{\ca F}(V)=C_{C_{\ca F}(V)}(V)$, in order to obtain the 
desired result from 6.9. 
\qed 
\enddemo

The next result is [lemma 3.2 in He].

\proclaim {Theorem 6.12} Let $\ca F$ be the fusion system of a proper locality. Then $\ca F^s$ is 
$\ca F$-closed. 
\endproclaim 

\demo {Proof} Fix a proper locality $(\ca L,\D,S)$ on $\ca F$, with $\D=\ca F^c$. 
Let $U,V$ be subgroups of $S$. By 1.14, 
$N_{\ca F}(U)\cong N_{\ca F}(V)$ if $U$ and $V$ are $\ca F$-conjugate subgroups of $S$, so $\ca F^s$ is 
$\ca F$-invariant. Clearly $S\in\ca F^s$. Thus, we are reduced to 
showing that $\ca F^s$ is closed with respect to overgroups in $S$. 

Among all 
$V\in\ca F^s$ such that some overgroup of $V$ in $S$ is not subcentric in $\ca F$, choose $V$ so that 
$|V|$ is as large as possible. Then there exists an overgroup $P$ of $V$ in $S$ such that $P\notin\ca F^s$ 
and such that $V$ has index $p$ in $P$. Then $P\leq N_S(V)$. Let $U\in V^{\ca F}$ such that $U$ is fully 
normalized in $\ca F$. As $\ca F$ is inductive there is an $\ca F$-homomorphism 
$\phi:N_S(V)\to N_S(U)$ such that 
$V\phi=U$. Then $P\phi\notin\ca F^s$ as $\ca F^s$ is $\ca F$-invariant, and so we may 
assume to begin with that $V$ is fully normalized in $\ca F$. Set $\ca F_V=N_{\ca F}(V)$. Replacing $P$ 
with a suitable $\ca F_V$-conjugate of $P$, we may assume that 
$P$ is fully normalized in $\ca F_V$. 

Suppose that $P\in(\ca F_V)^s$. As $C_{\ca F}(P)=C_{\ca F_V}(P)$ we then obtain $P\in\ca F^s$ from 
6.9. Thus $P\notin(\ca F_V)^s$, so $\ca F_V$ is a counterexample, and we may therefore assume that 
$\ca F_V=\ca F$. Then $O_p(\ca F)\in\ca F^c$, so $P\in\D$, and $\ca L$ is a group of characteristic $p$. 
Then $N_{\ca L}(P)$ is of characteristic $p$ by 2.7(b), and since $P$ is fully normalized in $\ca F$ 
we have also $N_{\ca F}(P)=\ca F_{N_S(P)}(N_{\ca L}(P))$. 
Then $O_p(N_{\ca F}(P))\in N_{\ca F}(P)^c$, and hence $O_p(N_{\ca F}(P))\in\ca F^c$ by 1.19. 
Thus $P\in\ca F^s$, as required.  
\qed 
\enddemo 

The following result is immediate from 6.10 and Theorem A1. 

\proclaim {Corollary 6.13} Let $(\ca L,\D,S)$ be a proper locality on $\ca F$. Then there exists an 
expansion of $\ca L$ to a proper locality $(\ca L^s,\ca F^s,S)$ on $\ca F$. 
\qed 
\endproclaim 

The locality $(\ca L^s,\ca F^s,S)$ (unique up to isomorphism, by Theorem A1) will be referred to as the 
{\it subcentric closure} of $(\ca L,\D,S)$, or (if the role of $S$ is understood) of $\ca L$. 

\vskip .1in 
The following result is [lemma 3.3 in He]. 

\proclaim {Lemma 6.14} Let $\ca F$ be the fusion system on $S$ of a proper locality, and let $P\leq S$ 
with $O_p(\ca F)P\in\ca F^s$. Then $P\in\ca F^s$. 
\endproclaim 

\demo {Proof} Among all counter-examples, choose $P$ so that $|P|$ is as large as possible. Set 
$Q=O_p(\ca F)$, and let $\phi:P\to P'$ be an $\ca F$-isomorphism such that $P'$ is fully normalized in 
$\ca F$. Then $\phi$ extends to an $\ca F$-isomorphism $PQ\to P'Q$, so 
$P'Q\in\ca F^s$. Thus we may assume that $P$ is fully normalized in $\ca F$. 

Set $\ca E=N_{\ca F}(P)$ and set $D=N_Q(P)$. Then $D\leq O_p(\ca E)$. If $D\leq P$ then $P=PQ$, and $P$ 
is not a counter-example. Thus $D\nleq P$, and $DP\in\ca F^s$ by the maximality of $|P|$. Since 
$N_S(P)\leq N_S(PD)$ we may assume that $P$ was chosen so that also $PD$ is fully normalized in $\ca F$. 
Set $R=N_S(PD)$ and set $G=N_{\ca L}(DP)$. We may take $\D$ to be $\ca F^s$ by 6.13, so $G$ is a subgroup 
of $\ca L$, and is of characteristic $p$. Here $R\in Syl_p(G)$, and 
$\ca F_{R}(G)=N_{\ca F}(DP)$. Since each $\ca E$-homomorphism extends to 
an $N_{\ca F}(DP)$-homomorphism, $\ca E$ is then the fusion system of $N_G(P)$ at $N_S(P)$. As 
$G$ is of characteristic $p$, so is $N_G(P)$ by 2.7(b). Thus $O_p(\ca E)$ is centric in $\ca E$, 
and hence centric in $\ca F$ by 1.19. Thus $P\in\ca F^s$. 
\qed 
\enddemo

We end the section with a few further results which will be needed in Part III.  

\proclaim {Lemma 6.15} Let $\ca L$ be a proper locality on $\ca F$, let $P\in\ca F^{cr}$, let $T\leq S$ be 
strongly closed in $\ca F$, and let $\ca E$ be an inductive fusion system on $T$ such that $\ca E$ is a 
subsystem of $\ca F$. Then there exists $Q\in P^{\ca F}$ with $Q\cap T\in\ca E^c$. 
\endproclaim 

\demo {Proof} Set $U=P\cap T$ and set $A=N_{C_T(U)}(P)$. Then $[P,A]\leq U$, and thus $A$ centralizes the 
chain $(P\geq U\geq 1)$ of normal subgroups of the group $N_{\ca L}(P)$. Then $A\leq O_p(N_{\ca L}(P))$ 
by 2.7(c), and then $A\leq P$ by an application of 2.3 to the fusion systen $N_{\ca F}(P)$. Thus 
$C_T(U)\leq P$, and so $C_T(U)\leq U$. 

Let $V\in U^{\ca F}$ be fully normalized in $\ca F$. As $\ca F$ is inductive there exists an 
$\ca F$-homomorphism $\phi:N_S(U)\to N_S(V)$ with $U\phi=V$. Set $Q=P\phi$. Then $Q\in\ca F^{cr}$, and 
so $C_T(V)\leq V$ by the result of the preceding paragraph. Let $V'\in V^{\ca E}$. Then $V'\in V^{\ca F}$, 
and so there exists an $\ca F$-homomorphism $\psi:N_S(V')\to N_S(V)$ with $V'=V\psi$. Then 
$N_T(V')\psi\leq N_T(V)$ as $T$ strongly closed in $\ca F$, and thus $|N_T(V')|\leq |N_T(V)|$. This shows 
that $V$ is fully normalized in $\ca E$. As $\ca E$ is inductive, $V$ is fully centralized in $\ca E$ by 
1.13. As $C_T(V)\leq V$, $V$ is then centric in $\ca E$ by 1.10. 
\qed 
\enddemo 

\proclaim {Lemma 6.16} Let $(\ca L,\D,S)$ be a proper locality on $\ca F$, let $X$ be fully normalized in 
$\ca F$, and let $X\leq P\leq N_S(X)$ with $P$ fully normalized in $N_{\ca F}(X)$. Then there exists an 
$\ca F$-homomorphism $\phi:N_S(X)\cap N_S(P)\to S$ such that $Q:=P\phi$ is fully normalized in $\ca F$ 
and $Y:=X\phi$ is fully normalized in $N_{\ca F}(Q)$. Moreover, $\phi$ then induces an isomorphism 
$N_{N_{\ca F}(X)}(P)\to N_{N_{\ca F}(Q)}(Y)$ of fusion systems. 
\endproclaim 

\demo {Proof} First, choose an $\ca F$-homomorphism $\psi:N_S(X)\to S$ so that $Q:=P\psi$ is fully 
normalized in $\ca F$, and set $Y'=X\psi$. By 6.8 $N_{\ca F}(Q)$ is saturated, 
so $N_{\ca F}(Q)$ is inductive, and so there exists an $N_{\ca F}(Q)$-homomorphism 
$\xi:N_S(Y')\cap N_S(Q)\to N_S(Q)$ such that $Y:=Y'\xi$ is fully normalized in $N_{\ca F}(Q)$. Let 
$\psi_0$ be the restriction of $\psi$ to $N_S(X)\cap N_S(P)$. Then $\psi_0\circ\xi$ sends 
$N_S(X)\cap N_S(P)$ into $N_S(Y)\cap N_S(Q)$, and maps $X$ to $Y$ and $P$ to $Q$. 

Now let $\phi$ be any $\ca F$-homomorphism $N_S(X)\cap N_S(P)\to S$ such that $P\phi$ is fully normalized 
in $\ca F$ and such that $X\phi$ is fully normalized in $N_{\ca F}(P\phi)$, and write $Q=P\phi$ and 
$Y=X\phi$. Set $U=N_S(X)\cap N_S(P)$ and $V=N_S(Y)\cap N_S(Q)$. Then $U\phi\leq V$. The 
argument of the preceding paragraph shows that there exists an $\ca F$-homomorphism $V\to U$, so 
$U\phi V$, and thus $\phi$ induces an $\ca F$-isomorphism $\a:U\to V$. It follows from two applications of 
1.5 that $\a$ is fusion-preserving. 
\qed 
\enddemo 

\proclaim {Lemma 6.17} Let $\ca F$ be a saturated fusion system on $S$, let $X$ be fully normalized in 
$\ca F$, and let $P$ be a subgroup of $N_S(X)$ containing $X$. Set $\ca F_X=N_{\ca F}(X)$. Then: 
\roster 

\item "{(a)}" $P\in(\ca F_X)^c$ if and only if $P\in\ca F^c$. 

\item "{(b)}" $P\in(\ca F_X)^s$ if and only if $P\in\ca F^s$. 

\endroster 
\endproclaim 

\demo {Proof} As $\ca F_X$ is saturated by 6.8, we may appeal to 1.18 to obtain the existence of an 
$\ca F_X$-conjugate $Q$ of $P$ such that both $Q$ and $R:=O_p(N_{\ca F_X}(Q))$ 
are fully normalized in $\ca F_X$. Then $Q$ and $R$ are fully centralized in $\ca F$, by 1.13. If 
$P\in(\ca F_X)^c$ then $C_S(Q)=C_{C_S(X)}(Q)\leq Q$, so $Q\in\ca F^c$, and then $P\in\ca F^c$. Similarly, 
if $R\in(\ca F_X)^c$ then $R\in\ca F^c$, whence $Q\in\ca F^s$ and $P\in\ca F^s$. This establishes 
the forward direction of the ``if and only if" in both (a) and (b). As $\ca F_X$ is a subsystem of 
$\ca F$ the reverse direction in each of the two cases is obvious.  
\qed 
\enddemo

\proclaim {Lemma 6.18} Let $(\ca L_1,\D_1,S)$ and $(\ca L_2,\D_2,S)$ be proper localities on the same fusion 
system $\ca F$. For $i=1$ and $2$ let $\ca L_i^+$ be the expansion of $\ca L_i$ to $\D_1\cup\D_2$. Then 
there exists an isomorphism $\a:\ca L_1^+\to\ca L_2^+$ such that $\a$ restricts to the identity map on $S$. 
\endproclaim 

\demo {Proof} By the main theorem of [Ch1] there exists a proper locality $(\ca L,\D,S)$ on $\ca F$ with 
$\D=\ca F^c$, and if $(\ca L',\D,S)$ is any other such locality on $\ca F$ then there is an 
isomorphism $\b:\ca L\to\ca L'$ which restricts to the identity map on $S$. This result, in combination 
with Theorem A1, yields the lemma. 
\qed 
\enddemo

\proclaim {Lemma 6.19} Let $(\ca L,\D,S)$ be a proper locality on $\ca F$, with subcentric closure 
$\ca L^s$. Let $H$ be a subgroup of $\ca L^s$. Then there exists $g\in \ca L^s$ such 
that conjugation by $g$ is defined on $H$, $H^g$ is a subgroup of $\ca L$, and conjugation by $g$ 
is an isomorphism $H\to H^g$. 
\endproclaim 

\demo {Proof} For each subgroup $P\leq S$ set  $P_0=P$, and then for each $k\geq 0$ define $P_{k+1}$ 
recursively, by $P_{k+1}=O_p(N_{\ca F}(P_k))$. We first show: 
\roster 

\item "{(*)}" For each $P\leq S$ there exists $Q\in P^{\ca F}$ such that each of the groups $Q_k$ is fully 
normalized in $\ca F$. Moreover, if $P\in\ca F^s$ then the sequence $(Q_k)_{k\geq 0}$ stabilizes at a 
member of $\ca F^{cr}$. 

\endroster 
For the proof of (*): Suppose that we are given $Q\in P^{\ca F}$ and $k\geq 0$ such that $Q_i<Q_{i+1}$ for 
all $i$ with $0\leq i\leq k$, and such that each $Q_i$ is fully normalized in $\ca F$. Further, let 
$Q$ be chosen so that $k$ is as large as possible. Set $R=Q_{k+1}$ and let $R'\in R^{\ca F}$ be fully 
normalized in $\ca F$. Thus there exists an $\ca F$-homomorphism $\phi:N_S(R)\to N_S(R')$ with 
$R\phi=R'$. Each of the groups $N_S(Q_i)$ $(0\leq i\leq k)$ normalizes $R$, and then each of the groups 
$N_S(Q_i)\phi$ is fully normalized in $\ca F$. Replacing $Q$ by $Q\phi$ thus yields a 
contradiction to the maximality of $k$. Thus all of the groups $Q_k$ are fully normalized. 
Suppose now that $P\in\ca F^s$. Then $Q_1$ is $\ca F$-centric, and since there exists $k$ with 
$Q_k=Q_{k+1}$ it follows that $Q_k\in\ca F^{cr}$. Thus (*) holds. 

Now let $H$ be a subgroup of $\ca L^s$. By I.2.11(a) there exists $P\in\ca F^s$ with $H\leq N_{\ca L^s}(P)$. 
Let $Q\in P^{\ca F}$ satisfy the conditions in (*). Then $Q=P^g$ for some $g\in\ca L^s$, and I.2.3(b) 
shows that there is an isomorphism $c_g:N_{\ca L^s}(P)\to N_{\ca L^s}(Q)$ given by conjugation by $g$. 
Notice that $N_{\ca L^s}(Q_k)\leq N_{\ca L^s}(Q_{k+1})$ for all $k$. Then (*) shows that $H^g$ normalizes 
a member of $\ca F^{cr}$, and thus $H^g$ is a subgroup of any restriction of $\ca L^s$ to a proper 
locality on $\ca F$. In particular, $H^g$ is a subgroup of $\ca L$. 
\qed 
\enddemo

\vskip .2in 
\noindent 
{\bf Section 7: $O^p_{\ca L}(\ca N)$ and $O^{p'}_{\ca L}(\ca N)$} 
\vskip .1in 

\vskip .1in 
\definition {Definition 7.1} Let $(\ca L,\D,S)$ be a locality, let $\ca N\norm\ca L$ be a partial 
normal subgroup of $\ca L$, and set $T=S\cap\ca N$. Set 
$$ 
\Bbb K=\{\ca K\norm\ca L\mid\ca KT=\ca N\}\ \ \text{and}\ \ \Bbb K'=\{\ca K'\norm\ca L\mid T\leq\ca K\},  
$$ 
and set 
$$ 
O^p_{\ca L}(\ca N)=\bigcap\Bbb K\ \ \text{and}\ \ O^{p'}_{\ca L}(\ca N)=\bigcap\Bbb K'. 
$$ 
Write $O^p(\ca L)$ for $O^p_{\ca L}(\ca L)$, and $O^{p'}(\ca L)$ for $O^{p'}_{\ca L}(\ca L)$. 
\enddefinition

\proclaim {Proposition 7.2} $O^p_{\ca L}(\ca N)\in\Bbb K$ and $O^{p'}_{\ca L}(\ca N)\in\Bbb K'$
\endproclaim 

\demo {Proof} Let $\ca K_1,\ca K_2\in\Bbb K$, set $\ca K=\ca K_1\cap\ca K_2$, and set $T_i=S\cap\ca K_i$ 
($i=1,2$). Let $x\in\ca K_1$. Then $x\in\ca N=\ca K_2T$, so we may write $x=yt$ with $y\in\ca K_2$ and 
$t\in T$. Then $S_x=S_y$, so $(y\i,x)\in\bold D$ and $y\i x=t$. Thus $t\in S\cap\ca K_1\ca K_2$, and so 
$t\in T_2T_1$ by [He1, Theorem A]. This shows that $\ca K_1\leq\ca K_2(T_2T_1)$. As 
$$ 
\ca K_2(T_2T_1)=(\ca K_2T_2)T_1=\ca K_2T_1 
$$ 
by I.2.9, we obtain $\ca K_1\leq\ca K_2T_1$. The Dedekind lemma (I.1.10) then yields  
$$ 
\ca K_1=\ca K_1\cap\ca K_2T_1=(\ca K_1\cap\ca K_2)T_1=\ca KT_1\leq\ca KT, 
$$ 
and so $\ca K_1T\leq\ca KT$. As $\ca K_1T=\ca N$ we conclude that $\ca KT=\ca N$, and thus $\ca K\in\Bbb K$. 
As $\Bbb K$ is finite, iteration of this procedure yields $O^p_{\ca L}(\ca N)\in\Bbb K$. 
The proof that $O^{p'}_{\ca L}(\ca N)\in\Bbb K'$ is simply the observation that $T\leq\bigcap\Bbb K'$. 
\qed 
\enddemo

\proclaim {Lemma 7.3} Let ``$*$" be either of the symbols ``$p$" or ``$p'$, let $(\ca L,\D,S)$ be a 
proper locality, and let $(\ca L^+,\D^+,S)$ be an expansion of $\ca L$. Let $\ca N\norm\ca L$, and for 
any partial normal subgroup $\ca K\norm\ca L$ let $\ca K^+$ be the corresponding partial normal subgroup 
of $\ca L^+$ given by Theorem A2. Then $O^{*}_{\ca L}(\ca N)^+=O^{*}_{\ca L^+}(\ca N^+)$. 
\endproclaim 

\demo {Proof} Write $\Bbb K^+$ for the set of all $\ca K^+$ with $\ca K\in\Bbb K$. Then 
$$ 
(\bigcap\Bbb K)^+\leq\bigcap(\Bbb K^+)
$$ 
as $\bigcap(\Bbb K^+)$ is a partial normal subgroup of $\ca L^+$ containing $\bigcap\Bbb K$. The 
reverse inclusion is given by Theorem A2, along with the observation that 
$$ 
\ca L\cap(\bigcap(\Bbb K^+))=\bigcap\{\ca L\cap\ca K^+\}_{\ca K\in\Bbb K}=\bigcap\Bbb K. 
$$
Thus the lemma holds for ``$p$", and the same argument applies to ``$p'$". 
\qed 
\enddemo 

\proclaim {Lemma 7.4} Let $\ca N$ and $\ca M$ be partial normal subgroups of $\ca L$, with $\ca N\leq\ca M$. 
Let ``$*$" be either of the symbols ``$p$" or ``$p'$. Then $O^*_{\ca L}(\ca N)\leq O^*_{\ca L}(\ca M)$.   
\endproclaim 

\demo {Proof} Since $O^{p'}_{\ca L}(\ca M)$ is a partial normal subgroup of $\ca L$ containg $S\cap\ca N$, 
we have $O^{p'}_{\ca L}(\ca N)\leq O^{p'}_{\ca L}(\ca M)$ by definition. Set $\ca K=O^p_{\ca L}(\ca M)$, 
set $\bar{\ca L}=\ca L/\ca K$, and let $\r:\ca L\to\bar{\ca L}$ be the canonical projection. Then 
$(S\cap\ca M)\r=\ca M\r\geq\ca N\r$, and hence $\ca N\r=(S\cap\ca N)\r$. Subgroup correspondence (I.4.7) 
then yields $\ca N\leq\ca K(S\cap\ca N)$, and then $O^p_{\ca L}(\ca N)\leq\ca K$ by definition. 
\qed 
\enddemo

\enddocument